\documentclass[11pt]{amsart} 
\usepackage{verbatim, latexsym, amssymb, amsmath,color}
\usepackage{epsfig}
\def\R{\mathbb R}
\def\RP{\mathbb {RP}}
\def\N{\mathbb N}
\def\Z{\mathbb Z}

\def\dmn{\mathrm{dmn}}

\theoremstyle{remark}
\newtheorem*{rmk}{Remark}

\theoremstyle{definition}

\title{Min-Max theory and the Willmore conjecture}
\author{Fernando C. Marques and Andr\'e Neves}
\address{Instituto de Matem\'atica Pura e Aplicada (IMPA) \\ Estrada Dona Castorina 110 \\ 22460-320 Rio de Janeiro \\ Brazil}
\email{coda@impa.br}
\address{Imperial College London\\ Huxley Building \\ 180 Queen's Gate \\ London SW7 2RH \\ United Kingdom}
\email{a.neves@imperial.ac.uk}
\thanks{The first author was partly supported by CNPq-Brazil, FAPERJ, Math-Amsud and the Stanford Department of Mathematics. The second author was partly supported by Marie Curie IRG Grant and ERC Start Grant.}

\begin{document}

\maketitle

\begin{abstract}
{ {In 1965,  T. J. Willmore  conjectured that the integral of the square of the mean curvature of a torus immersed in $\mathbb{R}^3$ is at least $2\pi^2$. We prove this conjecture using the min-max theory of minimal surfaces. }}
\end{abstract}

\setcounter{tocdepth}{1}
\tableofcontents

\section{Introduction}

The most basic geometric invariants of a closed surface $\Sigma$ immersed in Euclidean three-space are the
Gauss curvature $K$ and the mean curvature $H$. These invariants have been studied in differential geometry since its very beginning. The total integral of the Gauss curvature is a topological
invariant by the Gauss-Bonnet theorem. The integral of the square of the mean curvature, known as the Willmore energy, is specially interesting because it has the 
remarkable property of being invariant under conformal transformations of $\mathbb{R}^3$ \cite{blaschke,whitej}.  This
fact was already known to Blaschke \cite{blaschke} and Thomsen \cite{thomsen} in the 1920s (see also \cite{whitej}).  

 {Sometimes called bending energy, the Willmore energy appears naturally in some physical contexts}.  {For instance, it had been proposed in 1812 by Poisson} \cite{poisson}   {and later by Germain} \cite{germain}  {to describe elastic shells}.  {In mathematical biology it appears  in the Helfrich model} \cite{helfrich}  {as one of the terms that contribute to the energy of cell membranes}.

If we fix the topological type of $\Sigma$ and ask the question of what is the optimal immersion of $\Sigma$ in $\mathbb{R}^3$, it is natural to search among solutions to geometric variational problems. 
 It is not difficult to
see that the Willmore energy is minimized, among the class of all closed surfaces, precisely by the round spheres with value $4\pi$. The global problem of minimizing the Willmore energy among the class of immersed tori was proposed by T. J. Willmore  \cite{willmore}.

The main purpose of this paper is to prove the {\it Willmore Conjecture}:
\subsection{Willmore Conjecture (1965, \cite{willmore})}\label{conjecture} \textit{The integral of the square of the mean curvature of a torus immersed in $\mathbb{R}^3$ is at least $2\pi^2$.}

\medskip

The equality is achieved by the torus of revolution whose generating circle has radius 1 and center at distance $\sqrt{2}$ from  the axis of revolution:
$$
(u,v) \mapsto \big( (\sqrt{2} +\cos\, u) \cos\, v, (\sqrt{2}+\cos\,u)\sin\,v, \sin\,u) \in \mathbb{R}^3.
$$
This torus can also be seen as a stereographic projection of the Clifford torus $S^1(\frac{1}{\sqrt{2}}) \times S^1(\frac{1}{\sqrt{2}}) \subset S^3.$

The Willmore conjecture can be reformulated as a question about surfaces in the three-sphere because  if $\pi:S^3\setminus \{(0,0,0,1)\} \rightarrow \mathbb{R}^3$ denotes the stereographic projection
and $\Sigma \subset S^3\setminus \{(0,0,0,1)\}$ is a closed surface, then
\begin{equation}\label{conformal.invariance}
\int_{\tilde{\Sigma}}\tilde{H}^2 d \tilde{\Sigma}=\int_\Sigma (1+H^2) d\Sigma.
\end{equation}
Here $H$ and $\tilde{H}$ are the mean curvature functions of $\Sigma \subset S^3$ and $\tilde{\Sigma}=\pi(\Sigma) \subset \mathbb{R}^3$, respectively.

The conformal invariance of \eqref{conformal.invariance} motivates the following definition.  {Unless otherwise stated,} we will assume throughout the paper that surfaces are smooth and connected.

\subsection{Definition} The {\bf Willmore energy} of a closed surface $\Sigma \subset S^3$ is the quantity:
$$
\mathcal{W}(\Sigma) = \int_{\Sigma} (1+H^2) \, d\Sigma.
$$
Here $H$ denotes the mean curvature of $\Sigma$, i.e., $H=\frac{k_1+k_2}{2}$ where $k_1$ and $k_2$ are the principal curvatures. 

Note that if $F:S^3\rightarrow S^3$ is a conformal map, then $\mathcal{W}(F(\Sigma))=\mathcal{W}(\Sigma)$.

The Willmore conjecture follows as a consequence of our main theorem:
 
 \subsection*{Theorem A}\label{willmore.conjecture.theorem}
\textit{Let $\Sigma \subset S^3$ be an   embedded closed  surface of genus $g\geq 1$. Then
$$
\mathcal{W}(\Sigma) \geq 2\pi^2,
$$
and the equality holds if and only if $\Sigma$ is the Clifford torus up to conformal transformations of $S^3$.}

Theorem A indeed implies the Willmore conjecture because Li and Yau  \cite{li-yau}  proved  that if an immersion $f:\Sigma \rightarrow S^3$ covers a point $x\in S^3$ at least $k$ times,  then $\mathcal{W}(\Sigma) \geq 4\pi k$ .  {Therefore  a  non-embedded surface $\Sigma$ has   $\mathcal{W}(\Sigma) \geq 8\pi>2\pi^2$.}

If $\Sigma$ is a  critical point for the functional $\mathcal{W}$, we say that $\Sigma$ is a {\em Willmore surface}. The Euler-Lagrange equation 
for this variational problem, attributed by Thomsen \cite{thomsen} to Schadow, is
$$\Delta H + 2(H^2-K)H=0,$$
where $K$ denotes the Gauss curvature.  {Hence}  the image of a minimal surface under a conformal transformation of $S^3$  is {a Willmore surface}. (Minimal surfaces in $S^3$ with arbitrary genus were constructed by Lawson \cite{lawson70}). These are the simplest examples of Willmore surfaces, but not the only ones.  Bryant \cite{bryant2} found and classified  {immersed Willmore spheres} and Pinkall \cite{pinkall} constructed infinitely many embedded Willmore tori in $S^3$ which are not conformal to a minimal surface. Weiner \cite{weiner} checked that the second variation of $\mathcal{W}$ at the Clifford torus is nonnegative.

The existence of a  torus that minimizes the Willmore energy was established by Simon \cite{simon93}. His work
was later extended to surfaces of higher genus by Bauer and Kuwert \cite{bauer-kuwert} (see also \cite{kusner96}). We note that the existence of minimizers among higher genus surfaces in three-space also follows from our work, since Theorem A immediately implies the Douglas-type condition of \cite{simon93}. The minimum Willmore energy among all  orientable closed surfaces of genus $g$ is less than $8\pi$ \cite{lawson70, kuhnel-pinkall, kusner89}, and converges to $8\pi$ as $g\rightarrow \infty$ \cite{kuwert-li-schatzle}. The minimum Willmore energy among all immersed projective planes in $\mathbb{R}^3$
is known to be $12\pi$ \cite{bryant, kusner}.

Conjecture \ref{conjecture} was known to be true in some particular cases.  Willmore himself \cite{willmore71}, and independently
 Shiohama and Takagi \cite{shiohama-takagi}, proved it when the torus is a tube of constant radius around a space curve in $\mathbb{R}^3$. Chen \cite{chen} proved it for conformal images of flat tori in $S^3$  {(see \cite{topping} and \cite{ammann} for related results).} Langer and Singer \cite{langer-singer} proved it for tori of revolution (see also \cite{pinkall2} for a generalization). Langevin and Rosenberg  \cite{langevin-rosenberg} proved that any embedded knotted torus  $\Sigma$ in $\mathbb{R}^3$ satisfies
$\int_\Sigma |K|\, d\Sigma \geq 16\pi$. (Recall that a torus is knotted if it is not isotopic to the standard embedding.) Since 
$\int_\Sigma H^2 \, d\Sigma \geq \frac12 \int_\Sigma |K|\, d\Sigma$ for any torus $\Sigma \subset \mathbb{R}^3$, we conclude that  $\mathcal{W}(\Sigma) \geq 8\pi$ if $\Sigma$ is knotted. Li and Yau \cite{li-yau} 
introduced the notion of conformal volume and proved the conjecture for a class of conformal structures on $T^2$ that includes that of the square torus. The family of  conformal structures for which their method applies was later enlarged by Montiel and Ros \cite{montiel-ros}. Ros \cite{ros} proved the conjecture for tori $\Sigma \subset S^3$ that are invariant under the antipodal map. This result also follows from the work of  Topping \cite{topping, topping2005} on integral geometry. The conjecture was also known to be true for tori in $\mathbb{R}^3$ that are symmetric with respect to a point (Ros \cite{ros2001}).

Due to its connection to mathematical biology, {evidence for} the  fact that the Clifford torus and its {Dupin} cyclides minimize the Willmore energy was {experimentally} observed in membranes with the aide of a microscope by Mutz and Bensimon \cite{bensimon} (see also \cite{bensimon2}).

Finally,   {our understanding of the analytical aspects of the Willmore equation has  been greatly improved in recent years thanks to the work of Kuwert-Sch\"atzle (e.g. \cite{kuwert-schatzle}) and Rivi\`ere (e.g. \cite{riviere}).}

 % In particular,  $\mathcal{W}(\Sigma) \geq 8\pi$ if $\Sigma$ is not embedded.

 The next result is a corollary of Theorem A, but in fact we will prove it first. This theorem rules out the existence of a minimal surface of higher genus 
 {in $S^3$}  with area less than $2\pi^2$:

\subsection*{Theorem B}\label{minimal.area.theorem}
\textit{Let $\Sigma\subset S^3$ be an embedded closed minimal surface of genus $g \geq 1$. Then ${\rm area}(\Sigma) \geq 2\pi^2$,
and ${\rm area}(\Sigma)=2\pi^2$ if and only if  $\Sigma$ is the Clifford torus up to  isometries of $S^3$.}

\subsection{Remark} 
We  note that a   closed minimal surface $\Sigma \subset S^3$ of genus zero has to be totally geodesic
(Almgren \cite{almgren66}), and so its area is $4\pi$.
If $g\geq 1$ and $\Sigma$ is not embedded then ${\rm area}(\Sigma)=\mathcal{W}(\Sigma)\geq 8\pi > 2\pi^2$, by Li and Yau \cite{li-yau}.
\medskip

Finally, Theorem B will follow from the min-max theorem below.  The relevant definitions are in Sections \ref{gmt.definitions} and \ref{almgren.pitts.section}.
\subsection*{Theorem C}\label{existence.theorem}
\textit{Let $\Sigma \subset S^3$ be an embedded   closed  surface of genus $g\geq 1$, and let $\Pi$ be the  homotopy class associated with $\Sigma$ (see Definition \ref{homotopy.class.sigma}) with width ${\bf L}(\Pi).$
Then there exists an embedded closed minimal surface $\tilde{\Sigma} \subset S^3$ such that 
$$4\pi<{\rm area}(\tilde{\Sigma}) = {\bf L}(\Pi) \leq \mathcal{W}(\Sigma).$$
}

 {
Theorems B and C together immediately imply the next corollary. The corollary presents the Clifford torus as the min-max surface of a 5-dimensional family in $S^3$.

\subsection*{Corollary D}\label{clifford.minmax}
\textit{Let $\widehat \Pi$ be the homotopy class associated with the Clifford torus $\widehat \Sigma=S^1(\frac{1}{\sqrt{2}}) \times S^1(\frac{1}{\sqrt{2}}) \subset S^3$. Then
$$
{\bf L}(\widehat \Pi) = {\rm area}(\widehat \Sigma) = 2\pi^2.
$$
}
}

We  give an outline of our proof in the next section. Very briefly, to each embedded closed surface  $\Sigma$ in $S^3$, we associate a continuous $5$-parameter family  of surfaces (integral 2-currents with boundary zero, to be more precise)  in $S^3$ such that the area of each surface in the family is bounded above
by $\mathcal{W}(\Sigma)$. This family is parametrized by a map $\Phi$ defined on the 5-cube  $I^5$, and is constructed so that
\begin{itemize}
\item $\Phi(x,0)=\Phi(x,1) = 0$ (trivial surface) for any $x\in I^4$,
\item $\Phi(x,t)$ is an oriented round sphere in $S^3$ for any $x\in \partial I^4$, $t\in [0,1]$,
\item $\{\Phi(x,t)\}_{t\in [0,1]}$ is a homotopically nontrivial sweepout of $S^3$ for any $x\in \partial I^4$,
%\item $\sup\{{\rm area}(\Phi(x,t)):(x,t) \in I^5\} \leq E(\gamma_1,\gamma_2)$.
\end{itemize}
If ${\rm genus}(\Sigma)\geq 1$, this map $\Phi$ has the crucial property that its restriction to $\partial I^4 \times \{1/2\}$ is a homotopically nontrivial map
into  the space of oriented great spheres, which is homeomorphic to $S^3$. The min-max theory developed  in this paper  shows that for any such family  $\Phi$ there must exist  $y\in I^5$ such that 
${\rm area}(\Phi(y)) \geq 2\pi^2$. 

\medskip

{\bf Acknowledgements:} The authors would like to thank Brian White for his constant availability and helpful discussions. The authors are also thankful to Richard Schoen for his friendliness and encouragement while this
work was being completed.   Finally we would like to thank Rob Kusner for his interest and useful comments. Part of this work was done while the authors were visiting Stanford University.

%%%%%%%%%%%%%%%%%%%%%%%%%%%%%%%%%%%%%%%%%%%%%%%%%%%%%%%%%%
%%%%%%%%%%%%%%%%%%%%%%%%%%%%%%%%%%%%%%%%%%%%%%%%%%%%%%%%%%%%%%%

\section{Main ideas and organization}

 {We outline our proof of the Willmore conjecture. For the purpose of this discussion, we will ignore several  technical issues until Section \ref{technique}. Until then, we will appeal mainly to intuition in order to explain
the principal ideas behind our approach. 

\subsection{Min-Max Theory}\label{ideas.almgren}
We begin by describing the min-max theory of minimal surfaces in an informal way. We restrict our discussion to the case of 2-surfaces in a compact Riemannian three-manifold $M$.

Let $I^n=[0,1]^n$, and suppose we have a continuous map $\Phi$ defined  on $I^n$ such that $\Phi(x)$ is a compact surface with no boundary in $M$ for each $x\in I^n$. Two such maps $\Phi$ and $\Phi'$ are homotopic to each other relatively to 
$\partial I^n$ if there exists a continuous map $\Psi$, defined on $I^{n+1}$, such that:
\begin{itemize}
\item $\Psi(y)$ is a compact surface with no boundary in $M$ for each $y \in I^{n+1}$;
\item $\Psi(0,x)=\Phi(x)$ and $\Psi(1,x) = \Phi'(x)$ for each $x \in I^n$;
\item $\Psi(t,x)=\Phi(x)=\Phi'(x)$ for every $t \in I$, $x \in \partial I^n$.
\end{itemize}
The set $\Pi$ of all maps $\Phi'$ that are homotopic to $\Phi$ is called the homotopy class of $\Phi$.  The {\it width}
of $\Pi$ is then defined to be the min-max invariant:
$$
{\bf L}(\Pi) = \inf_{\Phi'\in \Pi}\, \sup_{x\in I^n} \, {\rm area}(\Phi'(x)).
$$

For instance, we could define $\Phi(s)= \{x_4=2s-1\}\subset S^3$ for $s\in [0,1]$. If $\Pi_1$ denotes its homotopy class, one should have ${\bf L}(\Pi_1)=4\pi$. Informally, $\Phi$ can be thought of as an element of $\pi_1(\mathcal{S}, \{0\})$, where $\mathcal{S}$ denotes the space of 2-surfaces in $S^3$ (0 means the trivial surface, of area zero).

The main goal of what we call min-max theory is to {realize the width as} the area of a minimal surface.
The prototypical result is:
\subsection*{Min-Max Theorem}\textit{If
$${\bf L}(\Pi)>\sup_{x\in \partial I^n}{\rm area}(\Phi(x)),$$
then there exists a smooth embedded closed minimal surface $\Sigma\subset M$ (possibly disconnected, with multiplicities) whose area  is equal to  ${\bf L}(\Pi).$ Moreover, if $\{\Phi_i\}$ is a sequence of maps in $\Pi$ such that 
$$
\lim_{i\rightarrow\infty} \sup_{x\in I^n} {\rm area}(\Phi_i(x)) = {\bf L}(\Pi),
$$
then we can choose $\Sigma$ to be the limit, as $i \rightarrow \infty$, of $\Phi_i(x_i)$ for some $x_i\in I^n$.
}

\subsection{Remark}\label{index.conjecture}
By analogy with standard Morse theory, and since $n$ is the number of parameters, one should expect that the index of 
$\Sigma$ as a minimal surface is at most $n$.  In general verifying this could be a delicate issue.

\subsection{Canonical family}\label{ideas.canonical}

Let  $B^4$ be the unit ball.  For every $v\in B^4$  we consider   the conformal map
$$F_v:S^3 \rightarrow S^3, \quad F_v(x) = \frac{(1-|v|^2)}{|x-v|^2}(x-v) -v.$$ 
Note that  if $v\neq 0$ then $F_v$ is a centered dilation of $S^3$ that fixes $v/|v|$ and $-v/|v|$.
To each smooth embedded closed surface $\Sigma \subset S^3$, we 
 associate a canonical  five-dimensional family of surfaces:
$$
\Sigma_{(v,t)}= \partial \left\{x\in S^3 :d_v(x)<\,t\right\},\quad (v,t) \in B^4 \times [-\pi,\pi].
$$
Here $d_v:S^3 \rightarrow S^3$ denotes the signed distance function to the oriented surface $\Sigma_v=F_v(\Sigma)$, {which becomes well defined }after we choose a unit normal vector field $N$ to $\Sigma$ in $S^3$. The distance
is computed with respect to the standard metric of $S^3$.   Note that
$\Sigma(v,\pi)=\Sigma(v,-\pi)=\emptyset$ for every $v\in B^4$.

The fundamental relation between the canonical family and the Willmore energy is given by Ros  \cite{ros} (see also \cite{heintze-karcher}):
 \begin{equation}\label{ideas.heinz}
 {\rm area}(\Sigma_{(v,t)}) \leq \mathcal{W}(\Sigma_v)= \mathcal{W}(\Sigma)\quad\mbox{for all }(v,t) \in B^4 \times [-\pi,\pi],
 \end{equation}
 where the last equality follows from the conformal invariance of the Willmore energy.

 \subsection{Boundary blow-up}\label{ideas. boundary.blowup2}
In view of \eqref{ideas.heinz}, we would like to apply the min-max method to the $5$-dimensional family $$\{\Sigma_{(v,t)}\}_{(v,t)\in B^4\times [-\pi,\pi]}.$$
Unfortunately  this family is not continuous in any reasonable sense, if we try to extend it to 
$\overline B^4\times [-\pi,\pi]  \approx  I^5$. As $v\in B^4$ converges to $p\in \Sigma$, we will see that the
limit depends on the angle of convergence.
In fact, if $$v_n=|v_n|(\cos(s_n)p+\sin(s_n)N(p))$$ is a sequence in $B^4$ converging to $p\in\Sigma$, i.e., $|v_n|$ tends to one, $|v_n|<1$,  and $s_n$ tends to zero, the limit of $\Sigma_{(v_n,t)}$ is  the geodesic sphere
$$\partial B_{\frac{\pi}{2}-\theta+t}(-\sin(\theta)p-\cos(\theta)N(p)),$$
where
$$\theta = \lim_{n\to\infty}\arctan\frac{s_n}{1-|v_n|}\in\left[-\frac{\pi}{2}, \frac{\pi}{2}\right].$$

\subsection{Remark} As $v\in B^4$ converges to $p\in S^3 \setminus \Sigma$, $\Sigma_{(v,t)}$ converges to
$$\partial B_{\pi+t}(p)\quad\mbox{or}\quad\partial B_{t}(-p),$$
depending on which connected component of $S^3 \setminus \Sigma$ contains $p$.

\medskip

In order to fix the failure of continuity, and after computing every boundary limit, we reparametrize the canonical family to
make it continuous on $\overline B^4 \times [\-\pi,\pi]$. This is done by ``blowing-up'' $\overline B^4$ along the
surface $\Sigma$, {a procedure which we describe now}. 

We first choose $\varepsilon>0$  to be small, and $\Omega_\varepsilon$ to be a   tubular neighborhood of radius $\varepsilon$ around $\Sigma$ in $\overline B^4$:
$$
\Omega_\varepsilon =\{(1-s_1)( \cos(s_2)p + \sin(s_2)N(p)): |(s_1,s_2)| < \varepsilon,  s_1\geq 0\}.
$$
Then we construct a continuous map $T:\overline B^4\rightarrow \overline B^4$  such that:
\begin{itemize}
\item $T$ maps $B^4\setminus \overline \Omega_\varepsilon$ homeomorphically onto $B^4$;
\item $T$ maps $\overline{\Omega}_\varepsilon$ onto $\Sigma$ by nearest point projection;
\item the map $$C(v,t)=\Sigma_{(T(v),t)},\quad (v,t) \in (B^4\setminus \overline\Omega_\varepsilon) \times [-\pi,\pi],$$ admits a continuous extension to  $\overline{(B^4\setminus \Omega_\varepsilon)} \times[-\pi,\pi]$, which we still denote by $C$. 
\end{itemize}

Finally  we extend $C$ to $\Omega_\varepsilon$ so that $C$ is constant along the radial directions. The resulting map $C$, defined on $\overline{B}^4 \times [-\pi,\pi]$,   satisfies  the following  properties:
 \begin{itemize}
 \item [(i)]${\rm area}(C(v,\pi))={\rm area}(C(v,-\pi))=0$ for every $v\in \overline{B}^4$;
  \item [(ii)] $C(v,t)$ is a geodesic sphere whenever $v \in S^3 \cup \overline{\Omega}_\varepsilon$;
  \item [(iii)]  for each $v\in S^3$, there exists a unique $s(v) \in [-\pi/2,\pi/2]$ such that $C(v,s(v))$ is a great sphere, i.e.,
  such that
  $$
  C(v,s(v)) = \partial B_{\pi/2}(\overline{Q}(v))
  $$  
  for some $\overline{Q}(v) \in S^3$.
  \end{itemize}
  If we {take into account the} orientation, $\partial B_{\pi/2}(p) \neq \partial B_{\pi/2}(-p)$. Hence $\overline{Q}(v)$ is also unique.
  
   In particular,
  \begin{equation}\label{sup.boundary.area}
  \sup_{(v,t) \in \partial (\overline{B}^4 \times [-\pi,\pi])} {\rm area}(C(v,t)) = 4\pi.
  \end{equation}
  Because of condition (i), we can extend $C$ to be zero (trivial surface) on $\overline{B}^4 \times (\R\setminus [-\pi,\pi])$.  
  
\subsection{Min-max family}\label{ideas.phi}  To apply the min-max theory described earlier, we  will reparametrize  $C$ to get a map $\Phi$ defined on $I^5$.  The min-max family is given by 
$$
\Phi(x,t) = C(f(x), 2\pi(2t-1) + \hat s(f(x))), \quad  x\in I^4, \, t \in I,
$$ 
for some choice of homeomorphism $f:I^4 \rightarrow \overline{B}^4$ and some extension $\hat s:\overline{B}^4 \rightarrow [-\pi/2,\pi/2]$ of the function $s$ to $\overline{B}^4$.
Note that this reparametrization is chosen so that  when $x\in \partial I^4$, we have that
\begin{equation}\label{level.one.half}
\Phi(x,t) \quad \mbox{is a great sphere if and only if}\quad t=1/2.  
\end{equation}

The estimate (\ref{ideas.heinz}) becomes
\begin{equation}\label{relation.area.willmore}
\sup_{x\in I^5} {\rm area}(\Phi(x)) \leq  \mathcal{W}(\Sigma).
 \end{equation}
 
 From (\ref{sup.boundary.area}), we also get 
   \begin{equation}\label{sup.boundary.area.phi}
  \sup_{x \in \partial I^5} {\rm area}(\Phi(x)) = 4\pi.
  \end{equation}
  
  Informally, the min-max family $\Phi$ can be thought of as an element of the relative homotopy group $\pi_5(\mathcal{S}, \mathcal{G})$, where
  $\mathcal{S}$ denotes the space of 2-surfaces in $S^3$ as before and $\mathcal{G}$ denotes the space of geodesic spheres.

  \subsection{Degree of $\overline Q$}  The map $\Phi$ is continuous and defined on $I^5$, so let $\Pi$ be its homotopy class.    Because of (\ref{sup.boundary.area.phi}),   we have that $\sup_{x\in \partial I^5}{\rm area}(\Phi(x))=4\pi.$ Therefore we need to check that 
  $L(\Pi)>4\pi$, in order to apply the Min-Max Theorem to this class.  Of course this might not be the case if
  $\Sigma$ is a topological sphere, but we will prove that $L(\Pi)>4\pi$  whenever $g={\rm genus}(\Sigma) \geq 1$.
  
  The main topological ingredient in the proof of this fact is:
  \begin{center}
  {\it $\overline{Q}:S^3 \rightarrow S^3$ is a continuous map with degree equal to $g$.}
  \end{center}
  This means that the canonical family detects the genus of $\Sigma$, and this is what will make the min-max approach work.
  The above fact, derived in Section \ref{associated}, is a consequence of the Gauss-Bonnet Theorem.
  
  This has an important homological implication as follows.  First note that 
  \begin{equation}\label{phi.half}
  \Phi(x,1/2) =  \partial B_{\pi/2}(\overline{Q}(f(x)))
  \end{equation}
  for every $x \in \partial I^4$. Now let $\mathcal{T}$ denote the set of all unoriented great spheres 
  in $S^3$.  By associating to each sphere in $\mathcal{T}$ the line generated by its center, we see that $\mathcal{T}$ is naturally homeomorphic to $\RP^3$. If $|\Phi|(x)=|\Phi(x)|$ denotes the surface $\Phi(x)$ after forgetting orientations (the reason we introduce this will be explained in Section \ref{technique}), then $|\Phi|$ maps $\partial I^4 \times \{1/2\}$ into $\mathcal{T}$. The fact that ${\rm deg}(\overline{Q})=g$ and (\ref{phi.half}) then implies 
  \begin{equation}\label{2g.homology}
  |\Phi|_*(\partial I^4 \times \{1/2\}) = 2g \in H_3(\RP^3,\Z).
 \end{equation}
  This will play a crucial role in the proof that ${\bf L}(\Pi) > 4\pi$.

  \subsection{${\bf L}(\Pi)>4\pi$}\label{width.outline}  
  
  Here we assume $g\geq 1$. The proof is by contradiction, therefore assume we can find a sequence of maps $\{\phi_i\}_{i\in \N}$ in $\Pi$ such that
$$\sup_{x\in I^5}{\rm area}(\phi_i(x))\leq 4\pi+\frac{1}{i}.$$ 
Note that $\phi_i=\Phi$ on $\partial I^5$. 

First we summarize the argument. We will construct a 4-dimensional submanifold $R(i) \subset I^5$, with $\partial R(i) \subset \partial I^4 \times I$, that separates $I^4 \times \{0\}$ from $I^4 \times \{1\}$. We construct $R(i)$ so that  for every $x \in R(i)$, the surface $|\phi_i(x)|$ is close to a great sphere in $\mathcal{T}$. This can be used to
produce by approximation a continuous function {$$f_i: R(i) \rightarrow \mathcal{T}\quad\mbox{such that}\quad f_i((x,t))=|\Phi(x,1/2)|\quad\mbox{for}\quad   (x,t) \in \partial R(i).$$}
Since we prove that $\partial R(i)$ is homologous to
$\partial I^4 \times \{1/2\}$ in $\partial I^4 \times I$, the existence of $f_i$ implies that $|\Phi|_*(\partial I^4 \times \{1/2\}) =0$ in $H_3(\RP^3,\Z)$. This is in contradiction with (\ref{2g.homology}). 

We now give more details. In what follows $\varepsilon>0 $ is a fixed small number. We denote by $\overline{A}(i)$ the set of all $x\in I^5$ such that the distance of the surface $|\phi_i(x)|$ to $\mathcal{T}$  (in an appropriate sense) is at least $\varepsilon.$ Since $\phi_i$, like $\Phi$, vanishes on
$I^4 \times \{0\}$ and $I^4 \times \{1\}$, these sets are both contained in $\overline{A}(i)$. 

We define $A(i)$ to be  the connected component of $\overline{A}(i)$ that contains $I^4 \times \{0\}$. For the purpose of this  discussion, we assume $\overline{A}(i)$ and $A(i)$ are compact manifolds with boundary.

We claim that $A(i)$  does not intersect $I^4\times\{1\}$  if  $i$ is  sufficiently large. Suppose this is false. Then we  find, 
after passing to a subsequence, a sequence of continuous paths $$\gamma_i:[0,1]\rightarrow A(i) \subset \overline A(i)\quad\mbox{with} \quad\gamma_i(0)\in I^4\times\{0\},\quad \gamma_i(1)\in I^4\times\{1\}.$$
The maps $\sigma_i =\phi_i\circ\gamma_i$, defined on $I=[0,1]$, are all homotopic to each other. Their homotopy class $\Pi_1$, just like in the one-dimensional example {in Section \ref{ideas.almgren},} satisfies ${\bf L}(\Pi_1) = 4\pi$. 
Moreover, we have
$$4\pi={\bf L}(\Pi_1)\leq \sup_{t\in I}{\rm area}(\sigma_i(t))\leq \sup_{x\in I^5}{\rm area}(\phi_i(x)) \leq 4\pi+\frac{1}{i}.$$Therefore, by the Min-Max Theorem, we can find $t_i\in I$ such that $\sigma_i(t_i)$ converges to an embedded minimal  surface $S$ with area $4\pi$. We must have that $S$ is a great sphere,  but this contradicts the fact that the distance of $|\sigma_i(t_i)|=|\phi_i(\gamma_i(t_i))|$ to $\mathcal{T}$ is at least $\varepsilon$.

One immediate consequence of the claim is that 
 $$\partial A(i)\cap \partial I^5 \subset (\partial I^4 \times I) \cup (I^4 \times \{0\}).$$

 Let $R(i)$ be the closure of $\partial A(i)\cap {\rm int}(I^5)$. It follows from the definition of $A(i)$ that
 \begin{equation}\label{distance.phi}
 d(|\phi_i(x)|, \mathcal{T})    \leq \varepsilon\quad {\mbox{for every}\quad x \in R(i).}
 \end{equation}
 In particular, 
 $\partial R(i) \subset  \partial I^4 \times I$. In fact it follows from (\ref{level.one.half}) that, given any $\delta>0$, we can choose $\varepsilon>0$ sufficiently small so that 
 \begin{equation}\label{Jdelta}
 \partial R(i) \subset  \partial I^4 \times [1/2-\delta,1/2+\delta].
 \end{equation}

 Let $C(i) = \partial A(i) \cap (\partial I^4 \times I)$. Since $\partial A(i)$ has no boundary,  we get that
 $$
 \partial C(i) = \partial R(i) \cup \partial (I^4 \times \{0\}).
 $$
 Therefore, since $C(i) \subset \partial I^4 \times I$, we have that  $\partial R(i)$ is homologous to $\partial I^4 \times\{0\}$ in $\partial I^4 \times I$. Consequently, $\partial R(i)$ is also  homologous to $\partial I^4 \times\{1/2\}$ in $\partial I^4 \times I$.
 
 Now let $\hat{\Phi}(x,t)=|\Phi(x,1/2)| \in \mathcal{T}$ for $x\in \partial I^4$. Because $\phi_i=\Phi$ on $\partial I^5$, we get from (\ref{Jdelta}) that ${|\phi_i|}_{|\partial R(i)}$ is close to $\hat{\Phi}_{|\partial R(i)}$. We use this, together with (\ref{distance.phi}), to approximate $|\phi_i|$ on $R(i)$ by a continuous map $f_i:R(i) \rightarrow \mathcal{T}$ 
 such that $f_i=\hat{\Phi}$ on $\partial R(i)$. This implies in homology that 
$$ \hat \Phi_{*}[\partial R(i)]={f_i}_{*}[\partial R(i)]=[{f_i}_{\#}\partial (R(i))]=[\partial {f_i}_{\#}(R(i))]=0.$$
On the other hand, we have 
 $$\hat \Phi_{*}[\partial R(i)]=\hat \Phi_{*}[\partial I^4\times\{1/2\}]=|\Phi|_{*}([\partial I^4\times\{1/2\}])=2g\in H_3(\RP^3,\Z).$$ We have reached a contradiction. 
 
 \subsection{Proof of Theorem B}\label{thmB.outline}
Let $\Sigma$ be the minimal surface with least area among all minimal surfaces in $S^3$ with genus greater than or equal to 1.  (The existence of $\Sigma$ follows from standard arguments in Geometric Measure Theory. This is explained in Appendix \ref{F1.section}.) The  area of $\Sigma$ is of course bounded above by $2\pi^2$, the area of the Clifford torus.

We claim that ${\rm index}(\Sigma)\leq 5$. This claim implies, by  a theorem of Urbano \cite{urbano}, that $\Sigma$ must be the Clifford torus up to isometries of $S^3$. 

Suppose, by contradiction, that ${\rm index}(\Sigma)>6$. If $\{\Sigma_{(v,t)}\}_{(v,t) \in B^4 \times [-\pi,\pi]}$ denotes the canonical family, then (\ref{ideas.heinz}) gives
$$
\sup_{(v,t) \in B^4 \times [-\pi,\pi]} {\rm area}(\Sigma_{(v,t)}) \leq \mathcal{W}(\Sigma) = {\rm area}(\Sigma). 
$$ 
The last equality follows from the fact that  $\Sigma$ is a minimal surface. The fact that $\Sigma$ is minimal also implies that  the function $(v,t) \mapsto {\rm area}(\Sigma_{(v,t)})$ has an isolated global maximum point at $(0,0)$. Since we are assuming that the
index is strictly bigger than the dimension of the parameter space, we can slightly perturb $\{\Sigma_{(v,t)}\}$ in a  neighborhood of $(0,0)$ to produce a new family $\{\Sigma'_{(v,t)}\}$ with 
\begin{equation}\label{ideas.heinz.modified}
\sup_{(v,t) \in B^4 \times [-\pi,\pi]} {\rm area}(\Sigma'_{(v,t)}) < {\rm area}(\Sigma).
\end{equation}

Let $\Phi'$  be the min-max family produced out of  $\{\Sigma'_{(v,t)}\}$, just like we constructed $\Phi$ out of $\{\Sigma_{(v,t)}\}$. Let $\Pi'$ be the homotopy class of $\Phi'$. Since $\Phi'$ agrees with $\Phi$ on $\partial I^5$, and since   $g={\rm genus}(\Sigma) \geq 1$, we can argue similarly as in \ref{width.outline} to get
  ${\bf L}(\Pi') > 4\pi$.  Therefore, because of (\ref{sup.boundary.area.phi}),  we can apply the Min-Max Theorem to $\Pi'$ in order to find an embedded minimal surface $\widehat \Sigma$ (with possible multiplicities) in $S^3$ such that
$${\rm area}(\widehat \Sigma)={\bf L}(\Pi')>4\pi.$$
But it follows from (\ref{ideas.heinz.modified}) that
$${\bf L}(\Pi')\leq \sup_{x\in I^5}{\rm area}(\Phi'(x))<{\rm area}(\Sigma)\leq 2\pi^2.$$
Thus ${\rm area}(\widehat \Sigma)< {\rm area}(\Sigma) \leq 2\pi^2$.

The area of any embedded minimal  surface in $S^3$ is at least  $4\pi$. It follows that the multiplicity of $\widehat \Sigma$ must be equal to one (otherwise ${\rm area}(\widehat \Sigma)\geq 8\pi$). Moreover, since ${\rm area}(\widehat \Sigma) >4\pi$ we get that  ${\rm genus}(\widehat \Sigma) \geq 1.$ Since ${\rm area}(\widehat \Sigma)< {\rm area}(\Sigma)$, we obtain a contradiction with the least-area property of $ \Sigma$.  Therefore ${\rm index}(\Sigma)\leq 5$, and $\Sigma$ is the Clifford torus up to isometries of $S^3$.

\subsection{Proof of Theorem A}  Let $\Sigma$ be an embedded closed surface in $S^3$, not necessarily minimal, with genus $g\geq 1$. We can suppose $\mathcal{W}(\Sigma) < 8\pi$ (otherwise the theorem follows immediately). Let $\Phi$ be the min-max family associated with $\Sigma$, and let $\Pi$ be its homotopy class.  From \ref{width.outline}, we get that ${\bf L}(\Pi)>4\pi$. Because of (\ref{sup.boundary.area.phi}), we can apply the Min-Max Theorem to $\Pi$ in order to find an embedded minimal surface $\widehat \Sigma$ (with possible multiplicities) in $S^3$ such that
$${\rm area}(\widehat \Sigma)={\bf L}(\Pi)>4\pi.$$

But it follows from (\ref{ideas.heinz}) that
$${\bf L}(\Pi)\leq \sup_{x\in I^5}{\rm area}(\Phi(x))\leq \mathcal{W}(\Sigma) < 8\pi.$$
Thus $4\pi<{\rm area}(\widehat \Sigma)\leq \mathcal{W}(\Sigma)<8\pi$. As in \ref{thmB.outline}, this implies
that the multiplicity of $\widehat \Sigma$ is equal to one and that ${\rm genus}(\widehat \Sigma)\geq 1$. It follows from Theorem B that ${\rm area}(\widehat \Sigma) \geq 2\pi^2$. Hence $\mathcal{W}(\Sigma) \geq 2\pi^2$, and the Willmore Conjecture holds. The rigidity statement follows by a perturbation argument similar to the one in \ref{thmB.outline}. 

\subsection{The technique}\label{technique} We  discuss the technical work that is necessary  to 
rigorously implement the min-max argument described above. In this subsection we assume the reader is  familiar with some
concepts of Geometric Measure Theory (see Section \ref{gmt.defi} for definitions).

In 1981, building on the work of Almgren \cite{almgren-varifolds}, Pitts \cite{pitts} succeeded in proving by min-max methods that any compact Riemannian manifold of dimension $n\leq 7$ contains a smooth embedded closed  minimal hypersurface, {where the regularity for the case $n=7$ was provided by  Schoen and Simon in \cite{schoen-simon}.} The methods of \cite{almgren-varifolds} and  \cite{pitts} are based in tools from Geometric Measure Theory, and comprise what we refer to in this paper as the Almgren-Pitts Min-Max Theory. The surfaces of a min-max family in this theory are integral currents, while the convergence to the min-max minimal hypersurface is in the sense of varifolds. 

There are other treatments of the min-max theory,  such as \cite{smyth, colding-delellis}. These impose stronger regularity 
 and convergence conditions on the surfaces of a min-max family.  These conditions are not satisfied
 by our sets $\Sigma_{(v,t)}$. {In particular, the family  $\{\Sigma_{(v,t)}\}$ can exhibit the well-known phenomenon of cancellation of mass: the possibility that two pieces of the surface match with opposite orientations and cancel out.}

 %The canonical family of Section \ref{ideas.canonical} and the min-max family of Section \ref{ideas.phi} can be interpreted as maps taking values in $\mathcal{Z}_2(S^3)$, the space of 2-dimensional integral  currents with boundary zero.  These maps are  continuous in the flat topology of currents. Unfortunately, if
 %we see them as maps into the space of varifolds (by forgetting orientations) they no longer have to be continuous. This is due to the well-known phenomenon of cancellation of mass: the possibility that two pieces of the surface match with opposite orientations and cancel out.  In particular the mass functional is not continuous. This is the  main  technical difficulty that we have to deal with.
 
 In Section \ref{ideas.almgren}, we considered families of surfaces parametrized by  the $n$-cube. In reality, Almgren and Pitts work with a discretized version: the maps are defined on the vertices of grids in $I^n$ that become finer and finer. The notion of continuity is replaced by the concept of fineness of a map, and appropriate discretized notions of homotopy have to be provided. Pitts chooses to work with families of currents that are fine in the mass norm ${\bf M}$.  The advantage of using the  ${\bf M}$-norm in $\mathcal{Z}_2(S^3)$ is that it can easily be localized (unlike the ${\bf F}$-metric), making it ideal for area comparisons, cut and paste arguments, and thus, regularity theory. The other advantage is that the mass functional is continuous in the ${\bf M}$-norm, as in the ${\bf F}$-metric (but not in the flat topology). 
 
 The disadvantage is that even the simplest family, like the 1-dimensional family  $\{x_4=s\}$ described in Section   \ref{ideas.almgren}, is not continuous with respect to the mass norm. This issue is addressed  by discretizing  the family $\{x_4=s\}$, and then interpolating, which  means adding  currents to the family or grid so that it becomes fine in the ${\bf M}$-norm. This is done is a way that both the original and the new families  represent, under a suitable homomorphism, the same element in $H_3(S^3,\Z)$. The min-max procedure is then applied to the interpolated family.

In this work we deal with the technical difficulties mentioned above by  following Almgren-Pitts approach. The min-max family $\Phi$ is defined on $I^5$ (as in Section \ref{ideas.phi}), takes values in  $\mathcal{Z}_2(S^3)$, and is continuous in the flat topology. By discretizing and interpolating, we construct a sequence of discrete maps $\phi_i$ that are fine in the mass norm and approximate $\Phi$ in the flat topology. Since the original map $\Phi$ is already continuous in varifold sense when restricted to $\partial I^5$, we can take $\phi_i$ to approximate $\Phi$ on $\partial I^5$ in the {\bf F}-metric. We also need to keep the fact that the width is bounded by the Willmore energy of $\Sigma$. Therefore the interpolation has to be carried out in such a way that the supremum of ${\bf M}(\phi_i)$ is not much bigger than the supremum of ${\bf M}(\Phi)$.  

The sets $A(i)$ and $R(i)$ that appear in Section \ref{width.outline}  will be replaced by  cubical singular chains in the rigorous argument. This is more appropriate for the homological conclusions, and  fits nicely with  the discrete nature of $\phi_i$.  The reason we sometimes need to forget orientations and work with $|\phi_i|, |\Phi|$ instead of $\phi_i, \Phi$, as in Section \ref{width.outline}, is that  the convergence to the
minimal surface in the Min-Max Theorem, using the Almgren-Pitts Min-Max Theory, is in the sense of varifolds. 
Later $|T|$ will denote the varifold associated with the integral current $T$.

The construction of the interpolating maps $\phi_i$ follows basic ideas of Almgren and Pitts, but is quite lengthy and technical. We dedicate a considerable part of the paper to carry it out. 

\subsection{Organization} The remaining material of this paper is organized as follows. 

{The main work needed to prove the Willmore conjecture is in Part I. This {contains}  Sections  \ref{associated},   \ref{gmt.defi},  \ref{canonical.boundary.section},  \ref{minmax.family.section}, \ref{gmt.definitions},  \ref{almgren.pitts.section},  \ref{bound.width.section}, \ref{thmb.section}, and \ref{thma.section}.}

In Section \ref{associated}, we define the five-dimensional canonical family $\{\Sigma_{(v,t)}\}$ associated with an embedded closed surface $\Sigma$ in $S^3$. We prove that the  area of $\Sigma_{(v,t)}$  is bounded above by $\mathcal{W}(\Sigma)$, and we compute  the degree of the map $\overline{Q}:S^3 \rightarrow S^3$.

In Section \ref{gmt.defi}, we collect the notation and the definitions from Geometric Measure Theory that are
relevant in this paper.

In Section \ref{canonical.boundary.section}, we reparametrize the canonical family and then we extend it to obtain the continuous map $C$ (in the sense of currents).

In Section \ref{minmax.family.section}, we define the min-max family $\Phi$ to which we will apply the Almgren-Pitts Min-Max Theory.   We collect all of its relevant properties.

In Section \ref{gmt.definitions}, we  give the basic definitions of the Almgren-Pitts min-max theory, adapted to our setting.  

In Section \ref{almgren.pitts.section}, we state a theorem that produces a discrete sequence of maps, needed by Almgren-Pitts min-max theory, out of the min-max family $\Phi$.  We also discuss Pitts Min-Max Theorem, adapted to our setting. 

In Section \ref{bound.width.section}, we show that the width is strictly bigger than $4\pi$ if the genus of $\Sigma$ is at least one. 

In Section \ref{thmb.section}, we prove Theorem B.

In Section \ref{thma.section}, we prove Theorem A. 

The technical machinery that makes the min-max argument work is done  {in Part II}. This {contains}  Sections     \ref{concentration.section}, \ref{continuous.discrete},  \ref{discrete.continuous},  and \ref{proof.pulltight}.

 In Section \ref{concentration.section}, we prove  that the canonical family has no concentration
 of area. 

In Section \ref{continuous.discrete}, we construct the discrete sequence of maps mentioned in Section
\ref{almgren.pitts.section}. This is done by discretizing $\Phi$ and then interpolating. 

In Section \ref{discrete.continuous}, we prove an interpolation theorem that associates to a discrete map a continuous map in the mass norm. This is needed in the pull-tight argument of Section \ref{proof.pulltight}.

In Section \ref{proof.pulltight}, we adapt the pull-tight procedure of Almgren and Pitts to our setting.

In Appendix \ref{F1.section}, we use standard arguments of Geometric Measure Theory to show that 
there exists  a minimal surface  with least area among  all embedded closed minimal surfaces with genus $g\geq 1$  in $S^3$.

In Appendix \ref{conformal.images}, we compute the conformal images of geodesic spheres in $S^3$.  

In Appendix \ref{appendix.map}, we construct the  map ${\bf r}_m(j)$ used in Section \ref{continuous.discrete}.
}

\part*{Part I. Proof of the Willmore  conjecture}

\medskip

%%%%%%%%%%%%%%%%%%%%%%%%%%%%%%%%%%%%%%%%%%%%%%%%%%%%%%%%%%
%%%%%%%%%%%%%%%%%%%%%%%%%%%%%%%%%%%%%%%%%%%%%%%%%%%%%%%%%%%%%%%

\section{Canonical family: First properties}\label{associated}

Before we construct the canonical family we need to introduce some notation.

\subsection{Notation and definitions}\label{associated.t} 
We use the following notation:
\begin{itemize}
\item $B^4 \subset \mathbb{R}^4$  the open unit ball and  $S^3=\partial B^4$ the unit sphere.
\item
$$B^4_R(Q)=\{x \in \mathbb{R}^4:|x-Q|<R\}\quad\mbox{and}\quad B_r(p)=\{x\in S^3:d(x,p)<r\},$$ where $Q \in \mathbb{R}^4$, $p\in S^3$, $R,r>0$, and $d$ is the spherical geodesic distance.
\end{itemize}

For each $v \in B^4$, we consider the conformal map $$F_v :S^3 \rightarrow S^3,\quad F_v(x) = \frac{(1-|v|^2)}{|x-v|^2}(x-v) -v.
$$

Consider $\Sigma \subset S^3$ an embedded   closed surface of genus $g$. We make  several definitions regarding the geometry of a tubular neighborhood of $\Sigma$ in $\overline B^4$.

\begin{itemize}
\item $A$ and $A^*$  denote the disjoint connected components of $S^3\setminus \Sigma=A \cup A^*$. 
\item  $N$  denotes  the  unit normal to $\Sigma$ that points into $A^*$.
\item Denote
 $$D_+^2(r) = \{s=(s_1,s_2) \in \mathbb{R}^2 : |s| <r, s_1 \geq 0\}.$$ 
 \item If  $\varepsilon >0$ is sufficiently small,  the map $\Lambda:\Sigma \times D_+^2(3\varepsilon)\rightarrow \overline{B}^4$ given by
\begin{equation}\label{tubular.neighborhood}
\Lambda(p,s) = (1-s_1)( \cos(s_2)p + \sin(s_2)N(p))
\end{equation}
is a diffeomorphism onto a neighborhood of $\Sigma$ in $\overline{B}^4$. 
\item  Let $\Omega_r=\Lambda(\Sigma\times D_+^2(r))$ for all $r\leq 3\varepsilon.$
\end{itemize}

Consider the continuous map $T:\overline B^4\rightarrow \overline B^4$ such that
\begin{itemize}
\item $T$ is the identity on $\overline B^4\setminus \Omega_{3\varepsilon};$
\item on $\Omega_{3\varepsilon}$ we have $$T(\Lambda(p,s)))=\Lambda(p,\phi(|s|)s),$$ 
where $\phi$ is smooth, zero on $[0,\varepsilon]$, strictly increasing  on $[\varepsilon,2\varepsilon],$ and one on $[2\varepsilon,3\varepsilon].$
\end{itemize}
The map $T$ collapses a tubular neighborhood of $\Sigma$ onto $\Sigma$.

Define $$A_{v}= F_v(A),\quad A^*_v=F_v(A^*),\quad\mbox{and }\Sigma_v=F_v(\Sigma)= \partial A_{v}$$ and let
 $d_v:S^3 \rightarrow \mathbb{R}$ be the signed distance to $\Sigma_v\subset S^3$: 
\begin{eqnarray*}
d_v(x)=\left\{
\begin{array}{rl}
d(x,\Sigma_v) & \mbox{if \ } x\notin A_{v},\\
-d(x,\Sigma_v)& \mbox{if \ } x\in A_{v}.
\end{array}
\right.
\end{eqnarray*}

\subsection{Definition}\label{associated.family}
The {\em canonical family} of  $\Sigma$ is the five-dimensional family of $2$-rectifiable subsets of $S^3$ given by 
$$
\Sigma_{(v,t)} = \partial A_{(v,t)}{,}\quad\mbox{where}\quad A_{(v,t)}= \{x\in S^3 :d_v(x)<\,t{\,\}}
$$
 {and} $(v,t) \in B^4 \times [-\pi,\pi]$.
\vskip0.03in

\subsection{Remark}\label{canonical.remark}
\begin{enumerate}
\item Let $N_v$ be the normal vector to $\Sigma_v$ given by $N_v=D{F_v}(N)/|D{F_v}(N)|$ and 
consider   the smooth map  $$\psi_{(v,t)}:\Sigma_v  \rightarrow S^3, \quad \psi_{(v,t)}(y)=\exp_{y}(tN_v(y))=\cos t \, y + \sin t\,  N_v(y).$$
We have
$$
\Sigma_{(v,t)} \subset \psi_{(v,t)} (\{ {\rm Jac \ } \psi_{(v,t)} \geq 0{\})}
$$
and so $\Sigma_{(v,t)}$ is indeed a $2$-rectifiable set.
\item Notice that $A_{(v,0)}=A_{v}$, $A_{(v,\pi)}=S^3$, and $A_{(v,-\pi)}=\emptyset$, which means that $$\Sigma_{(v,0)}=\Sigma_v, \quad \Sigma_{(v,\pi)}=\emptyset,\quad\mbox{and}\quad\Sigma_{(v,-\pi)}=\emptyset.$$
\end{enumerate}
The importance of this family is described in the next theorem. A related result appears in Proposition 1 of \cite{ros}.

\subsection{Theorem}\label{heintze.karcher}
{\em  We have, for every $(v,t) \in B^4 \times (-\pi,\pi),$
$${\rm area}\left(\Sigma_{(v,t)}\right) \leq \mathcal{W}(\Sigma).$$%-\frac{\sin^2 t}{2}\int_{\Sigma}{|\mathring{A}|^2}d\Sigma,
Moreover, if $\Sigma$ is  not a geodesic sphere and 
$$ {\rm area}\left(\Sigma_{(v,t)}\right) =\mathcal{W}(\Sigma),$$
then $t=0$ and $\Sigma_v$ is a minimal surface.
}

\begin{proof}
The following calculation can be found in \cite{ros}: 
\subsection{Lemma}\label{jacobian.calculation}
\textit{
We have
$$
{\rm Jac \ }\psi_{(v,t)}(y) = (1+H(v)^2) - (\sin t + H(v) \cos t)^2 -\frac{(k_{1}(v)-k_{2}(v))^2}{4}\sin^2 t,
$$
where $k_{1}(v)$ and $k_{2}(v)$ are the principal curvatures of $\Sigma_v$ at $y$, and $ H(v)=\frac{k_{1}(v)+k_{2}(v)}{2}$ is the mean curvature.
}

\begin{proof}
Let $\{e_1,e_2\} \subset T_y\Sigma_v$ be an orthonormal basis of principal directions, with principal curvatures $k_{1}(v)$ and $k_{2}(v)$, respectively. Hence
$$D{\psi_{(v,t)}}_{|y}e_i = (\cos t - k_{i}(v) \sin t) e_i,$$ from which we conclude that $${\rm Jac \ }\psi_{(v,t)}(y) = (\cos t - k_{1}(v)\sin t) (\cos t - k_{2}(v) \sin t). $$
The lemma follows by expanding this  out.
\end{proof}

Using this lemma we can finish the proof.
From Lemma \ref{jacobian.calculation}, the area formula,  and conformal invariance of the Willmore energy we obtain
 \begin{multline*}
 {\rm area}(\Sigma_{(v,t)}) 
 \leq {\rm area} \big(\psi_{(v,t)}(\{{\rm Jac \ }\psi_{(v,t)}(p)\geq0\})\big)\\
\leq \int_{\{{\rm Jac \ }\psi_{(v,t)}\geq 0\}} ({\rm Jac \ }\psi_{(v,t)}) \ d\Sigma_v\\
\leq \int_{\{{\rm Jac \ }\psi_{(v,t)}(p)\geq0\}} (1+ H({v})^2) -\sin^2 t\frac{(k_{1}(v)-k_{2}(v))^2}{4} \,d\Sigma_v\\
 {\leq  \int_{\Sigma} (1+ H({v})^2)\ d\Sigma_v=\mathcal{W}(\Sigma)}.
 \end{multline*}
If equality holds for some $(v,t)\in B^4\times  (-\pi,\pi)$, we obtain from the set of inequalities above that $\{{\rm Jac \ }\psi_{(v,t)}\geq 0\}=\Sigma$ and
$$\frac{\sin^2 t}{2}\int_{\Sigma_v}{|\mathring{A}|^2}d\Sigma_v=\frac{\sin^2 t}{2}\int_{\Sigma}{|\mathring{A}|^2}d\Sigma=0,$$
where $\mathring{A}$ denotes the trace-free part of the second fundamental form. This implies the rigidity statement.

%$$0\neq\int_{\Sigma}\frac{(k_{1}-k_{2})^2}{4}\,d\Sigma=\int_{\Sigma_v}\frac{(k_{1}(v)-k_{2}(v))^2}{4}\,d\Sigma_v\quad\mbox{for all }v\in B^4$$
%and so $t=0$, $t=\pi$, or $t=-\pi$. The last two cases are impossible by Remark \ref{canonical.remark} (2) and thus $t=0$. This means ${\rm area}(\Sigma_v)={\mathcal W}(\Sigma_v)$ and so $\Sigma_v$ is a minimal surface.

\end{proof}

\subsection{Extended Gauss map}\label{extended.map.gauss}
For every $p\in\Sigma$ and $k\in [-\infty,+\infty]$ consider
\begin{equation}\label{Q.R.}
\overline Q_{p,k}=-\frac{k}{\sqrt{1+k^2}}p-\frac{1}{\sqrt{1+k^2}}N(p)\in S^3.
\end{equation}
This induces a function
$\overline Q: \overline\Omega_{\varepsilon}\rightarrow S^3$ such that $$\overline Q\left(\Lambda(p,s)\right)=\overline Q_{p,k}\quad\mbox{where }k=\frac{s_2}{\sqrt{\varepsilon^2-s_2^2}}.$$ 
We extend this map  in the following way:
  \begin{eqnarray}\label{Qmap}
  \overline Q: S^3\cup {\overline\Omega_{\varepsilon}} \rightarrow S^3,\quad 
\overline Q(v)=\left\{
\begin{array}{rl}
-T(v) & \mbox{if  } v\in A^* \setminus \overline{\Omega}_\varepsilon,\\
T(v) & \mbox{if  } v\in A\setminus \overline{\Omega}_\varepsilon,\\
\overline Q(v) & \mbox{if  } v\in \overline{\Omega}_\varepsilon.\\
\end{array}
\right.
\end{eqnarray}
\begin{rmk}
If $p\in \Sigma$, i.e, $p=\Lambda(p,(0,0))$, then $\overline Q(p)=-N(p)$ is the classical Gauss map for surfaces in $S^3$.
\end{rmk}

The next theorem is  absolutely crucial to the proof of the Willmore conjecture.
 
 \subsection{Theorem}\label{degree.gauss.map}\textit{The map $\overline Q$ is continuous and $$\overline Q:S^3\rightarrow S^3$$has degree g.}

\begin{proof}
We start by showing that $\overline Q:S^3\rightarrow S^3$ is continuous. Clearly $Q$ is continuous on $S^3 \cap \overline{\Omega}_\varepsilon$, $A^* \setminus \overline{\Omega}_\varepsilon$, and $A \setminus \overline{\Omega}_\varepsilon$. Assume $$v=\Lambda(p,(0,t))=\cos t\, p+ \sin t\, N(p)\in \Omega_{2\varepsilon}.$$
If $v\in S^3\cap \overline {\Omega}_\varepsilon$ we see from \eqref{Qmap} that
$$\lim_{t\to\varepsilon_{-}}\overline Q(v)=\overline Q_{p,+\infty}=-p\quad\mbox{and}\quad\lim_{t\to -\varepsilon_{+}} \overline Q(v)=\overline Q_{p,-\infty}=p.$$
If $v\in A\setminus  \overline {\Omega}_\varepsilon$ we see from the definition  of $T$ and  \eqref{Qmap}  that
$$\lim_{t\to-\varepsilon_{-}}\overline Q(v)= \lim_{t\to-\varepsilon_{-}} T(v)=-\lim_{t\to 0_{-}}\Lambda(p,(0,t))=p.$$
If $v\in A^*\setminus  \overline {\Omega}_\varepsilon$ we see from the definition of $T$ and  \eqref{Qmap}  that
$$\lim_{t\to\varepsilon_{+}}\overline Q(v)= -\lim_{t\to\varepsilon_{+}} T(v)=-\lim_{t\to 0_{+}}\Lambda(p,(0,t))=-p.$$
Hence $\overline Q:S^3\rightarrow S^3$ is continuous.

\subsection{Lemma}\label{degree}\textit{
The degree of $\overline Q:S^3 \rightarrow S^3$ is $g$.
}
\begin{proof}

We will use the fact that $\overline{Q}$ is piecewise smooth. Let  $dV$ denote the volume form of $S^3$ and $ \nabla$ the induced connection on $S^3$.

Since $\overline Q=-T$ on $A^*\setminus \overline\Omega_\varepsilon$, we have from the definition of $T$ that $\overline Q$ is an orientation-preserving  diffeomorphism of $A^* \setminus \Omega_\varepsilon$ onto $-\overline{A^*}$.
Therefore 
\begin{equation}\label{grau1}
\int_{A^*\setminus \overline\Omega_\varepsilon} \overline Q^*(dV)=\int_{-A^*} dV= {\rm vol}(A^*).
\end{equation}
Since $\overline Q=T$ on $A\setminus\overline \Omega_\varepsilon$, we have from the definition of $T$ that $\overline Q$ is an orientation-preserving  diffeomorphism of $A \setminus \overline\Omega_\varepsilon$ onto $\overline{A}$.
Therefore 
\begin{equation}\label{grau2}
\int_{A\setminus \overline\Omega_\varepsilon} \overline Q^*(dV)=\int_{A} dV= {\rm vol}(A).
\end{equation}

Recall that $\{e_1,e_2,e_3\} \in T_pS^3$ is a positive basis if $\{e_1,e_2,e_3,p\}$ is a positive basis of $\mathbb{R}^4$, and $\{e_1,e_2\}\in T_p\Sigma$
is a positive basis if $\{e_1,e_2,N(p)\}$ is a positive basis of $T_pS^3$.

Consider the diffeomorphism
$
G:\Sigma\times [-\varepsilon,\varepsilon] \rightarrow S^3 \cap\overline{\Omega}_\varepsilon
$
defined by $$G(p,t)=\Lambda(p,(0,t))=\cos\, t \, p+ \sin\, t\, N(p).$$
The orientation of $\Sigma \times [-\varepsilon,\varepsilon]$ is chosen so that $\{e_1,e_2,\partial_t\}$ is a positive basis whenever $\{e_1,e_2\}$ is
a positive basis of $T\Sigma$. We have $${G_{*}}(e_1\wedge e_2\wedge\partial_t)_{|(p,0)}=e_1\wedge e_2\wedge N(p)$$ and thus $G$ is orientation preserving.

Consider ${Q}=\overline Q \circ G:\Sigma\times [-\varepsilon,\varepsilon]  \rightarrow S^3$ which is given by
$$
{Q}(p,t)=-\frac{t}{\varepsilon}\, p-\frac{\sqrt{\varepsilon^2-t^2}}{\varepsilon}\,N(p).
$$
Hence
\begin{equation*}
\int_{S^3 \cap\overline{\Omega}_\varepsilon} \overline Q^*(dV) =\int_{\Sigma\times [-\varepsilon,\varepsilon]} G^*(\overline Q^*(dV))
= \int_{\Sigma\times [-\varepsilon,\varepsilon]} {Q}^*(dV).
\end{equation*}

Let $\{e_1,e_2\}$ be a positive orthonormal basis of $T_p\Sigma$ which diagonalizes the second fundamental form:
$$\nabla_{e_i} N=-k_i e_i\quad\mbox{for }i=1,2.$$ 
We have
$$
DQ_{|(p,t)}(\partial_t ) = -\frac{1}{\varepsilon}\,p+\frac{t}{\varepsilon\sqrt{\varepsilon^2-t^2}}\,N(p),
$$
and
$$
DQ_{|(p,t)}(e_i)= \left(-\frac{t}{\varepsilon}+\frac{\sqrt{\varepsilon^2-t^2}}{\varepsilon}k_i\right)e_i\quad\mbox{for }i=1,2,
$$
and thus, denoting by $vol_{\mathbb{R}^4}$ the standard volume form of $\mathbb{R}^4$, we have
\begin{align*}
{Q}^*(dV)_{|(p,t)}(e_1,e_2,\partial_t)& =
dV_{|Q(p,t)}(D{Q}(e_1), DQ(e_2), DQ(\partial_t))\\
&={vol_{\mathbb{R}^4}}_{|Q(p,t)}(DQ(e_1), DQ(e_2), DQ(\partial_t),{Q}(p,t))\\
&=\left(-\frac{t}{\varepsilon}+\frac{\sqrt{\varepsilon^2-t^2}}{\varepsilon}k_1\right)\left(-\frac{t}{\varepsilon}+\frac{\sqrt{\varepsilon^2-t^2}}{\varepsilon}k_2\right)\frac{(-1)}{\sqrt{\varepsilon^2-t^2}}
\end{align*}
since
\begin{multline*}
 DQ_{|(p,t)}(\partial_t)\wedge Q(p,t)=\\
\left(-\frac{1}{\varepsilon}\,p+\frac{t}{\varepsilon\sqrt{\varepsilon^2-t^2}}\,N(p)\right) \wedge \left(-\frac{t}{\varepsilon}\, p-\frac{\sqrt{\varepsilon^2-t^2}}{\varepsilon}\,N(p)\right)\\
=-\frac{1}{\sqrt{\varepsilon^2-t^2}}\, N(p) \wedge p.
\end{multline*}

The Gauss equation implies that $K = 1 + k_1k_2$, where $K$ denotes the Gauss curvature of $\Sigma$ and so 
we conclude that
\begin{multline}\label{grau3}
\int_{\Sigma\times [-\varepsilon,\varepsilon]} {Q}^*(dV)\\
=- \int_{\Sigma}\int_{-\varepsilon}^\varepsilon \frac{1}{\varepsilon^2}\left(k_1k_2\sqrt{\varepsilon^2-t^2}-(k_1+k_2)t+\frac{t^2}{\sqrt{\varepsilon^2-t^2}}\right)dt \, d\Sigma\\
=-\frac{\pi}{2}\int_\Sigma (K-1) \, d\Sigma - \frac{\pi}{2} \int_\Sigma d\Sigma
=-\pi^2 \chi(\Sigma)
=\pi^2(2g-2).
\end{multline}
In the calculation above we have used that
\begin{itemize}
\item $\int_{-\varepsilon}^\varepsilon \sqrt{\varepsilon^2-t^2} dt =\varepsilon^2 \int_{-\pi/2}^{\pi/2} \cos^2 \theta\, d\theta =\frac{\pi \varepsilon^2}{2}$,
\item $\int_{-\varepsilon}^\varepsilon t\, dt=0$, 
\item $\int_{-\varepsilon}^\varepsilon \frac{t^2}{\sqrt{\varepsilon^2-t^2}} dt =\varepsilon^2 \int_{-\pi/2}^{\pi/2} \sin^2 \theta\, d\theta =\frac{\pi \varepsilon^2}{2}$.
\end{itemize}

Finally, since ${\rm vol}(S^3)=2\pi^2$, we combine \eqref{grau1}, \eqref{grau2}, and \eqref{grau3} to obtain
\begin{multline*}
\int_{S^3} \overline Q^*(dV)=\int_{A^*\setminus \overline\Omega_\varepsilon} \overline Q^*(dV)+\int_{A\setminus  \overline\Omega_\varepsilon}  \overline Q^*(dV)+\int_{S^3\cap \overline{\Omega}_\varepsilon}  \overline Q^*(dV)\\
= {\rm vol}(A^*)+{\rm vol}(A)+\int_{\Sigma\times [-\varepsilon,\varepsilon]} {Q}^*(dV)\\
=2\pi^2 + \pi^2(2g-2)
=2\pi^2g 
= g \cdot \int_{S^3} dV.
\end{multline*}

It follows that ${\rm deg}(Q)=g$.

\end{proof}
{This lemma finishes the proof of Theorem \ref{degree.gauss.map}.}
\end{proof}

{For technical reasons that will be relevant later, we need to ensure that the areas of the sets $\Sigma_{(v,t)}$ cannot concentrate at a point:

\subsection{Theorem}\label{no.concentration.mass}\textit{
For every  $\delta>0$, there exists $r>0$ such that
$$
{\rm area}(\Sigma_{(v,t)} \cap B_r(q)) \leq \delta\quad\mbox{for every }q\in S^3\mbox{ and } (v,t) \in B^4 \times [-\pi,\pi].
$$
}

The proof of Theorem \ref{no.concentration.mass} will be postponed to Section \ref{concentration.section}.
}

%%%%%%%%%%%%%%%%%%%
%%%%%%%%%%%%%%%%%%%%%%%%%%%%%%

\section{Definitions from Geometric Measure Theory}\label{gmt.defi}

{
 {In this section we recall some definitions and notation from Geometric Measure Theory. A standard reference   is the book of Simon \cite{simon}. Sometimes we will also follow the notation of  Pitts book \cite{pitts}.}
 
 Let $(M,g)$ be an orientable compact Riemannian three-manifold. We  assume $M$ is
isometrically embedded in $\mathbb{R}^L$. We denote by $B_r(p)$ the open geodesic
ball in $M$ of radius $r$ and center $p\in M$.

 The spaces we will work with in this paper are:
\begin{itemize}
\item the space ${\bf I}_k(M)$ of $k$-dimensional integral currents in $\mathbb{R}^L$ with support contained  in $M$;
\item the space ${\mathcal Z}_k(M)$ of integral currents $T \in {\bf I}_k(M)$ with  $\partial T=0$;
\item the closure $\mathcal{V}_k(M)$, in the weak topology, of the space of $k$-dimensional rectifiable varifolds in $\mathbb{R}^L$ with support contained in $M$.
\end{itemize}

Given $T\in {\bf I}_k(M)$,  we denote by $|T|$ and $||T||$ the integral varifold   and Radon measure in $M$ associated with $T$, respectively;
 given $V\in \mathcal{V}_k(M)$, $||V||$ denotes the Radon measure in $M$ associated with $V$.  
 If $U\subset M$ is an open set of finite perimeter, the associated current in ${\bf I}_{3}(M)$ is denoted by $[|U|]$.

The above spaces come with several relevant metrics. The   {\it mass} of $T \in {\bf I}_k(M)$,   defined by
$${\bf M}(T)=\sup\{T(\phi): \phi \in \mathcal{D}^k(\mathbb{R}^L), ||\phi|| \leq 1\},$$
induces the metric ${\bf M}(S,T)={\bf M}(S-T)$ on ${\bf I}_k(M)$.
 Here $\mathcal{D}^k(\mathbb{R}^L)$ denotes  the space of smooth $k$-forms in $\mathbb{R}^L$ with compact support, and $||\phi||$ denotes the comass norm of $\phi$.  
 
The {\it flat metric} is defined by 
$$\mathcal{F}(S,T)=\inf\{{\bf M}(P)+{\bf M}(Q): S-T=P+\partial Q, P\in{\bf I}_k(M), Q \in {\bf I}_{k+1}(M)\},$$
for $S, T\in {\bf I}_k(M)$. We also use $\mathcal{F}(T)=\mathcal{F}(T,0)$. Note that
$$\mathcal{F}(T)\leq  {\bf M}(T)\quad\mbox{ for all }T \in {\bf I}_k(M).$$

The  ${\bf F}$-{\it metric} on $\mathcal{V}_k(M)$ is defined in  {Pitts book} \cite[page 66]{pitts} as:
\begin{eqnarray*}
&&{\bf F}(V,W) = \sup \{V(f)-W(f) : f \in C_c(G_k(\mathbb{R}^L)),\\
&&\hspace{5.5cm} |f|\leq 1, {\rm Lip}(f) \leq 1\},
\end{eqnarray*}
for $V,W \in \mathcal{V}_k(M)$. Here $C_c(G_k(\mathbb{R}^L))$ denotes the space of all  real-valued continuous functions with compact support defined on $G_k(\mathbb{R}^L)$ - the $k$-dimensional {Grassmannian} bundle over $\mathbb{R}^L$. 
 The ${\bf F}$-metric  induces the varifold weak topology on $\mathcal{V}_k(M)$,  and it satisfies
$$ {\bf F}(|S|,|T|)\leq {\bf M}(S-T)\quad\mbox{ for all }S,T \in {\bf I}_k(M).$$

Finally,  the ${\bf F}$-{\it metric} on ${\bf I}_k(M)$ is defined by
$$ {\bf F}(S,T)=\mathcal{F}(S-T)+{\bf F}(|S|,|T|).$$

We assume that  ${\bf I}_k(M)$ and ${\mathcal Z}_k(M)$ both have the topology induced by the flat metric. When endowed with
the topology of the mass norm, these spaces will be denoted by  ${\bf I}_k(M;{\bf M})$ and ${\mathcal Z}_k(M;{\bf M})$, respectively. If endowed with the ${\bf F}$-metric, we will denote them by  ${\bf I}_k(M;{\bf F})$ and ${\mathcal Z}_k(M;{\bf F})$, respectively.   The space $\mathcal{V}_k(M)$ is considered with the weak topology of varifolds.

If ${\bf \nu}$ is either the flat, mass, or ${\bf F}$-metric, then
  $${\bf B}^{\bf \nu}_r(T)=\{S\in {\mathcal Z}_k(M): {\bf \nu}(T,S)< r\}.$$
  Given $\mathcal{A,B}\subset \mathcal{V}_k(M)$, we also define
  $${\bf F}(\mathcal{A},\mathcal{B})=\inf\{{\bf F}(V,W):V\in \mathcal{A}, W\in \mathcal{B}\}.$$
  
  The mass {\bf M} is continuous in the topology induced by the ${\bf F}$-metric, but not in the flat topology. In the flat topology the mass functional is only lower semicontinuous.
 Keep in mind that 
 $$\mathcal{F}(S-T)\leq  {\bf F}(S,T)\leq 2{\bf M}(S-T)$$
 for every $S,T\in{\bf I}_k(M)$.

The following lemma will be useful.
\subsection{Lemma}\label{flat+mass=f} \textit{Let $\mathcal{S}$ be a compact subset of ${\mathcal Z}_k(M;{\bf F})$. For every $\varepsilon>0$ there is $\delta>0$ so that for every $S\in \mathcal{S}$ and $T\in { {\mathcal Z}_k(M)}$
$$
{\bf M}(T)<{\bf M}(S)+\delta\mbox{ and }\mathcal{F}(T-S)\leq \delta\Rightarrow {\bf F}(S,T)\leq \varepsilon.
$$
}

\begin{proof}
 In \cite[page 68]{pitts}, it is observed that $\lim_{i\rightarrow \infty} {\bf F}(S,T_i)=0$ if and only if $\lim_{i\rightarrow \infty}{\bf M}(T_i)={\bf M}(S)$ and $\lim_{i\rightarrow \infty}\mathcal{F}(S-T_i)=0$, for $T_i,S \in {\mathcal Z}_k(M)$. The lemma then follows from the continuity properties of the mass functional and the compactness of $\mathcal{S}$ in ${\mathcal Z}_k(M;{\bf F})$, via a standard finite covering argument.
\end{proof}

%The mass {\bf M} is continuous in the topology induced by the ${\bf F}$-metric, but not in the flat %topology. In the flat topology the mass functional is only lower semicontinuous.
 %Keep in mind that 
 %$$\mathcal{F}(S-T)\leq  {\bf F}(S,T)\leq 2{\bf M}(S-T)$$
 %for every $S,T\in{\bf I}_k(M)$.
 
 Given a $C^1$-map $F:M\rightarrow M$, the push-forwards of $V\in \mathcal{V}_k(M)$ and $T\in {\bf I}_k(M)$ are denoted by  ${F}_{\#}(V)$ and ${F}_{\#}(T)$, respectively. Denote by $\mathcal{X}(M)$ the space of smooth vector fields of $M$ with the $C^1$-topology. 
The first variation 
$$\delta:\mathcal{V}_k(M)\times \mathcal{X}(M)\rightarrow \R $$
is defined as
$$\delta V(X)=\frac{d}{dt}_{|t=0}||{F_{t}}_{\#}(V)||(M),\quad\mbox{where } \frac{dF_t}{dt}_{|t=0}=X.$$
The first variation is continuous  with respect to the product topology of $\mathcal{V}_k(M)\times \mathcal{X}(M)$. 
Recall that a varifold $V$ is said to be {\it stationary} if $\delta V(X)=0$ for every $X\in \mathcal{X}(M)$.

We will also need the following definition. $I^n$ denotes the $n$-dimensional cube.
\subsection{Definition}\label{mass.def} Given a continuous map $\Phi:I^n\rightarrow {\mathcal Z}_2(M)$, with respect to the flat topology, we define
$${\bf m}(\Phi,r)=\sup\{||\Phi(x)||(B_r(p)):x\in I^n, p\in M\}.$$

}

%%%%%%%%%%%%%%%%%%
%%%%%%%%%%%%%%%%%%%

\section{Canonical family: Boundary blow-up}\label{canonical.boundary.section}

Following the discussion in Section \ref{ideas. boundary.blowup2}, we want to reparametrize and extend the canonical family to be defined on all of $\overline B^4\times [-\pi,\pi]$. {The resulting family will be continuous in the sense of currents.}

%Since we want to apply the Almgren-Pitts Min-Max theory,  we will adopt in this section the language of %Geometric Measure Theory and work with integral currents. {More precisely, we use 
%\begin{itemize}
%\item ${\bf I}_k(S^3)$ the space of  $k$-dimensional integral currents in $S^3$;
%\item ${\mathcal Z}_k(S^3)$ the space of $k$-dimensional integral currents $T\in {\bf I}_k(S^3)$ with $\partial %T=0$;
%\item $[|U|]\in {\bf I}_3(S^3)$ denotes the integral current associated to an open set of finite perimeter $U
%\subset S^3$;
%\item ${\bf M}(T)$ denotes the mass of $T\in {\bf I}_k(S^3)$ (see Section \ref{gmt.defi} for precise definition);
%\item  given $S,T \in {\bf I}_k(S^3)$, $\mathcal{F}(S,T)$ denotes the flat metric (see Section \ref{gmt.defi} for %precise definition).
%\end{itemize}
%}  

{The goal is to produce}, out of the
canonical family, a five-dimensional family  of integral currents of boundary zero that is continuous
in the flat topology of currents (Theorem \ref{canonical.family.continuous}).

For every $k\in [-\infty,+\infty]$, consider $$\overline r_k=\frac{\pi}{2}-\arctan k \in [0,\pi].$$ We note that
\begin{equation}\label{extrinsic.ball}
B^4_{\overline R_k}(\overline Q_{p,k})\cap S^3=B_{\overline r_k}(\overline Q_{p,k})\quad\mbox{where}\quad \overline R_k=\sqrt{2\left(1-\frac{k}{\sqrt{1+k^2}}\right)}.
\end{equation}
Consider also $\overline r: \overline \Omega_{\varepsilon}\rightarrow [0,\pi]$ given by
$$\bar r\left(\Lambda(p,s)\right)=\bar r_k,\quad\mbox{where }k=\frac{s_2}{\sqrt{\varepsilon^2-s_2^2}}.$$
We extend this function in the following way:
  \begin{eqnarray}\label{r.definition}
\overline r: S^3\cup \overline {\Omega_{\varepsilon}} \rightarrow [0,\pi],\quad \overline r(v)=\left\{
\begin{array}{rl}
0  & \mbox{if  } v\in A^* \setminus \overline{\Omega}_\varepsilon,\\
\pi & \mbox{if  } v\in A\setminus \overline{\Omega}_\varepsilon,\\
\overline r(v) & \mbox{if  } v\in \overline{\Omega}_\varepsilon.\\
\end{array}
\right.
\end{eqnarray}

The goal of this section is to prove the following result.
\subsection{Theorem}\label{canonical.family.continuous}\textit{The map below is well defined and continuous in the flat topology:
$$C: \overline B^4\times[-\pi,\pi]\rightarrow \mathcal{Z}_2(S^3),$$
 \begin{eqnarray*}
 C(v,t)=\left\{
\begin{array}{rl}
\partial [|A_{(T(v),t)}|]  & \mbox{if  } v\in B^4\setminus \overline\Omega_{\varepsilon},\\[0.5em] 
%\mbox{ } &\mbox{ }\\
\partial [|B_{\overline r(v)+t}(\overline Q(v))|] & \mbox{if  } v\in S^3\cup\overline \Omega_{\varepsilon}.\\
\end{array}
\right.
\end{eqnarray*}
{Furthermore,
$${\bf M}(C(v,t))\leq \mathcal{W}(\Sigma)\quad\mbox{for all}\quad(v,t)\in \overline B^4\times[-\pi,\pi]$$
and $C(v,\pi)=C(v,-\pi)=0$ for every $v\in \overline B^4$.}
}

\subsection{Preliminary results}

Given sets $A,B$ of $\R^4$, the symmetric difference is denoted by
$$A\,\Delta \,B=(A\setminus B)\cup (B\setminus A).$$

Recall the definition of the map $\Lambda$ in \eqref{tubular.neighborhood}. 
If $v_n\in B^4$ is a sequence converging to $p\in \Sigma$, then for all $n$ sufficiently large there are  unique $p_n\in \Sigma$ and $s_n\in D^2_+(3\varepsilon)$ so that $v_n=\Lambda(p_n,s_n)$. Necessarily, $p_n$ tends to $p$ and $s_n$ tends to zero. By passing to a subsequence, we can also assume that
$$\lim_{n\to\infty} \frac{{s_n}_2}{{{s_n}_1}}= k \in [-\infty, +\infty].$$

\subsection{Proposition}\label{convergence.sets}\textit{Consider a sequence $(v_n,t_n)\in B^4\times [-\pi,\pi]$ converging to $(v,t)\in \overline B^4\times [-\pi,\pi]$.
\begin{itemize}
\item[(i)]  If $v\in B^4$ then
$$\lim_{n\to\infty} {\rm vol}\left(A_{(v_n,t_n)}\,\Delta\, A_{(v,t)}\right)=0.$$
\item[(ii)] If $v\in A$ then
$$\lim_{n\to\infty} {\rm vol} \left(A_{(v_n,t_n)}\, \Delta\, B_{\pi+t}(v)\right)=0$$
and,  {given any $\delta>0$,
$$ \Sigma_{(v_n,t_n)}\subset \overline B_{\pi+t+\delta}(v)\setminus B_{\pi+t-\delta}(v)\quad\mbox{for all $n$ sufficiently large}.$$}
\item[(iii)] If $v\in A^*$ then
$$\lim_{n\to\infty} {\rm vol} \left(A_{(v_n,t_n)}\, \Delta\, B_{t}(-v)\right)=0$$
and,  {given any $\delta>0$,
$$ \Sigma_{(v_n,t_n)}\subset \overline B_{t+\delta}(-v)\setminus B_{t-\delta}(-v)\quad\mbox{for all $n$ sufficiently large}.$$}
\item[(iv)] If $v=p\in \Sigma$ and 
$$v_n=\Lambda(p_n,({s_n}_1,{s_n}_2))\quad\mbox{with}\quad \lim_{n\to\infty} \frac{{s_n}_2}{{{s_n}_1}}=k \in [-\infty, +\infty],$$
 then
$$
\lim_{n\to\infty} {\rm vol} \left(A_{(v_n,t_n)}\,\Delta\,B_{\overline{r}_{k}+t}(\overline{Q}_{p,k}) \right)=0
$$
and,  {given any $\delta>0$,
$$ \Sigma_{(v_n,t_n)}\subset \overline B_{\overline{r}_{k}+t+\delta}(\overline{Q}_{p,k})\setminus B_{\overline{r}_{k}+t-\delta}(\overline{Q}_{p,k})\quad\mbox{for all $n$ sufficiently large}.$$}
\end{itemize}
}
\begin{proof}
We denote by $N_v$ the normal vector to $\Sigma_v$ with the same direction as $DF_v(N)$. Consider the normal exponential map of $\Sigma_v$ given by
$$\exp_v:\Sigma_v\times\R\rightarrow S^3, \quad \exp_v(y,t)=\cos t \, y + \sin t\,  N_v(y).$$ 
 For every $x\in S^3$, there exists $y\in \Sigma_v$ such that $x = {\rm exp}_v(y,d_v(x)).$ In particular,
\begin{equation}\label{closest.point}
\left(A_{(v,t)} \setminus A_{(v,s)}\right)  \subset {\rm exp}_v(\Sigma_v \times [s,t))\quad\mbox{for }s\leq t.
\end{equation}

We now prove Proposition {\ref{convergence.sets}} (i).
Let $\delta>0$ and  choose $\eta>0$ such that
$$
{\rm vol}\left({\rm exp}_v(\Sigma_v \times [t-\eta,t+\eta])\right) \leq \delta.
$$

 The sequence of surfaces  $\Sigma_{v_n}$ converges smoothly to $\Sigma_v$, since $v_n$ tends to $ v \in B^4$. 
This, together with the triangle inequality and the fact that $t_n$ tends to $t$, implies that
we can choose $n_0$ such that 
$$
A_{(v,t-\eta)} \subset A_{(v_n,t_n)} \subset A_{(v,t+\eta)}\quad\mbox{for all }n\geq n_0.
$$
Hence, for $n\geq n_0$, we have
$$
A_{(v_n,t_n)} \,\Delta\, A_{(v,t)} 
 \subset\left( A_{(v,t+\eta)} \setminus A_{(v,t-\eta)}\right).
$$
From  \eqref{closest.point} we have 
$$
\left(A_{(v,t+\eta)} \setminus A_{(v,t-\eta)}\right)  \subset {\rm exp}_v(\Sigma \times [t-\eta,t+\eta])
$$
and thus
$$
 {\rm vol}\left(A_{(v_n,t_n)}\,\Delta\, A_{(v,t)}\right) \leq \delta
$$
for each $n\geq n_0$.

We now prove Proposition {\ref{convergence.sets}} (ii). Let $r>0$ be such that $B_r(v) \subset A$. Given $\delta >0$, there exists $n_0  \in \mathbb{N}$ such that for all $n\geq n_0$
\begin{equation}\label{inclusion.ii.sets}
\overline{B}_{\pi-\delta/2}(v) \subset F_{v_n}(B_r(v)) \subset  F_{v_n}(A)=A_{(v_n,0)}\quad\mbox{and}\quad|t_n-t|\leq \frac{\delta}{2}.
\end{equation}
In particular,
\begin{equation}\label{inclusion.ii.sets2}
\Sigma_{v_n} \subset B_{\delta/2}(-v)\quad\mbox{for all } n\geq n_0. 
\end{equation}

If $t\geq 0$, then from \eqref{inclusion.ii.sets}  and the triangle inequality we have, for all $n\geq n_0$,
$$S^3 \setminus B_{\delta}(-v)=\overline{B}_{\pi-\delta}(v) \subset A_{(v_n,-\delta/2)} \subset A_{(v_n,t_n)}.$$
Hence, because  $\overline{B}_{\pi+t}(v)=S^3$, 
\begin{equation}\label{inclusion.ii.sets3}
{\rm vol} \left(A_{(v_n,t_n)}\,\Delta\, {B}_{\pi+t}(v)\right) \leq {\rm vol}\,(B_\delta(-v))\quad\mbox{and}\quad \Sigma_{(v_n,t_n)}\subset B_\delta(-v).
\end{equation}
Notice that if $t>0$, (\ref{inclusion.ii.sets}) implies that $A_{(v_n,t_n)} = S^3$ and hence $\Sigma_{(v_n,t_n)}=\emptyset$ for any sufficiently large $n$.

If  $t <0$,  choose $n_1 \geq n_0$ such that $t_n <0$ for each $n \geq n_1$. We have
$$A_{(v_n,t_n)}\subset B_{\pi+t+\delta}(v) \quad\mbox{for all }n\geq n_1$$
because,  picking $x\in A_{(v_n,t_n)}$ and $y\in \Sigma_{v_n}$ with $d_{v_n}(x)=-d(x,y)$, we obtain from \eqref{inclusion.ii.sets2} and the triangle inequality
\begin{multline*}
d(x,-v)\geq d(x,y)-d(y,-v)= -d_{v_n}(x)-d(y,-v)
\geq -t_n-\frac{\delta}{2}\geq -t-\delta.
\end{multline*}
Also
$$B_{\pi+t-\delta}(v)\subset A_{(v_n,t_n)} \quad\mbox{for all }n\geq n_1$$
because if $x\in B_{\pi+t-\delta}(v)$, then $x\notin B_{\delta-t}(-v)$ and we obtain from \eqref{inclusion.ii.sets2}
$$ d(x,\Sigma_{v_n})>d(x, \partial B_{\delta/2}(-v))> -t+\frac{\delta}{2}\geq -t_n.$$
  Hence, for all $n\geq n_1$,
\begin{equation}\label{inclusion.ii.sets4}
\left(A_{(v_n,t_n)}\,\Delta\, B_{\pi+t}(v)\right)\cup  \Sigma_{(v_n,t_n)}\subset  \overline{B}_{\pi+t+\delta}(v)\setminus B_{\pi+t-\delta}(v)
\end{equation}
In any case Proposition \ref{convergence.sets}(ii) follows from \eqref{inclusion.ii.sets3} and  \eqref{inclusion.ii.sets4}, since we can choose $\delta$ arbitrarily small.

Proposition {\ref{convergence.sets}} (iii) is proven exactly in the same way as Proposition {\ref{convergence.sets}} (ii).

We now prove Proposition {\ref{convergence.sets}} (iv).

\subsection{Lemma}\label{two.balls}\textit{
There exists $r_0>0$ such that for every $p \in \Sigma$  we have
\begin{multline*}
B_{r_0}\left((\cos  r_0) p -(\sin r_0) N(p)\right)  \subset A,\\
\mbox{and}\quad \overline{A} \subset S^3 \setminus B_{r_0}\left((\cos  r_0) p + (\sin r_0) N(p)\right).
\end{multline*}
}

\begin{proof}
Choose $r_0>0$ sufficiently small such that for every $x \in S^3$ with $d(x,\Sigma) \leq r_0$, there exists a unique $q \in \Sigma$ such that the shortest
geodesic segment joining $x$ and $q$ is  orthogonal to $\Sigma$ at $q$. We must have  $d(x,q)=d(x,\Sigma)$. 

If  $x_1 = (\cos  r_0) p -(\sin r_0) N(p)$ and $x_2=(\cos  r_0) p + (\sin r_0) N(p)$, then  $d(x_1,\Sigma)=d(x_2,\Sigma)=r_0$. Therefore
$B_{r_0}(x_1) \cap \Sigma = B_{r_0}(x_2) \cap \Sigma =\emptyset$. The result follows since  $x_1\in A$ and   $x_2 \in A^*$.
\end{proof}

Write $v_n=\Lambda(p_n,({s_n}_1,{s_n}_2))$, where $k_n={s_n}_2/{s_n}_1$ tends to $k$ and $p_n$ tends to $p$.

Set $B_q = B_{\pi/2}(-N(q))=B^4_{\sqrt{2}}(-N(q))\cap S^3$ for $q \in \Sigma$. It follows from Lemma \ref{two.balls} that
\begin{multline*}
A\,\Delta\,B_{p_n} \subset S^3 \setminus \left(B_{r_0}\left((\cos  r_0) p_n + (\sin r_0) N(p_n)\right)\right.\\
 \cup B_{r_0}\left.\left((\cos  r_0) p_n - (\sin r_0) N(p_n)\right) \right).
\end{multline*} 

From Proposition \ref{sphere.image} of the Appendix \ref{conformal.images}, we obtain the existence of $C>0$ and $n_0\in \mathbb{N}$ such that 
$$S^3 \cap B^4_{\overline{R}_{k_n}-C\sqrt{a_n}}(\overline{Q}_{p_n,k_n}) \subset F_{v_n}(A)\subset S^3 \cap B^4_{\overline{R}_{k_n}+C\sqrt{a_n}}(\overline{Q}_{p_n,k_n})$$
for all $n \geq n_0$, where $a_n=\sqrt{1+k_n^2}\,{s_n}_1$. Notice that $a_n \rightarrow 0$.

Therefore, from \eqref{extrinsic.ball}, we see that for each $\delta >0$ there exists $n_1 \geq n_0$ such that for every $n\geq n_1$ we have
\begin{equation}\label{delta.close.inclusion}
B_{\overline{r}_{k_n}-\delta/2}(\overline{Q}_{p_n,k_n}) \subset F_{v_n}(A)\subset B_{\overline{r}_{k_n}+\delta/2}(\overline{Q}_{p_n,k_n})
\end{equation}
and
\begin{equation}\label{delta.close.inclusion2}
\Sigma_{v_n} \subset \overline{B}_{\overline{r}_{k_n}+\delta/2}(\overline{Q}_{p_n,k_n}) \setminus B_{\overline{r}_{k_n}-\delta/2}(\overline{Q}_{p_n,k_n}).
\end{equation}

Assume $\overline{r}_{k} \in (0,\pi)$ and $0<\delta<\min\{\overline{r}_{k}, \pi-\overline{r}_{k}\}$. The cases $\overline{r}_{k}=0$ and $\overline{r}_{k}=\pi$ can be dealt with as in the proof of Proposition {\ref{convergence.sets}} (ii). 

We can find $n_2 \geq n_1$ such that for each $n\geq n_2$ we have  $$|t_n-t|+ d(\overline{Q}_{p_n,k_n},\overline{Q}_{p,k})+|\overline{r}_{k_n}-\overline{r}_{k}|\leq \delta/2.$$ 
Thus, from \eqref{delta.close.inclusion}, we have $\overline{Q}_{p,k} \in F_{v_n}(A)$ and $-\overline{Q}_{p,k} \notin F_{v_n}(A)$ for $n\geq n_2$.

We claim
$$A_{(v_n,t_n)} \subset B_{\overline{r}_{k}+t+\delta}(\overline{Q}_{p,k})\quad\mbox{for all }n\geq n_2.$$

Let $n \geq n_2$, and $x \in A_{(v_n,t_n)}$. Then $d_{v_n}(x) <t_n$, and $x={\rm exp}_{v_n}(y,d_{v_n}(x))$ for some $y\in \Sigma_{v_n}$.

If $d_{v_n}(x) \geq 0$, we obtain from \eqref{delta.close.inclusion}
\begin{eqnarray*}d(x,\overline{Q}_{p,k}) &\leq& d(x,y) + d(y,\overline{Q}_{p,k})\\
&\leq& d(x,y) + d(y,\overline{Q}_{p_n,k_n})+d(\overline{Q}_{p_n,k_n},\overline{Q}_{p,k})\\
&\leq& d_{v_n}(x) + \overline{r}_{k_n}+\delta/2+d(\overline{Q}_{p_n,k_n},\overline{Q}_{p,k})\\
&<&   t_n + \overline{r}_{k_n}+\delta/2+d(\overline{Q}_{p_n,k_n},\overline{Q}_{p,k})\\
&<& \overline{r}_{k} +t + \delta.
\end{eqnarray*}

If $d_{v_n}(x) <0$, then $x \in F_{v_n}(A)$.  Thus, from \eqref{delta.close.inclusion},  any continuous path joining $x$ to   $-\overline{Q}_{p_n,k_n}$ must intersect $\Sigma_{v_n}$
and  using \eqref{delta.close.inclusion2} we obtain
\begin{align*}
d(x,-\overline{Q}_{p,k})&\geq  d(x,-\overline{Q}_{p_n,k_n})-d(\overline{Q}_{p_n,k_n},\overline{Q}_{p,k})\\
&\geq  d(x, \Sigma_{v_n}) +d(\Sigma_{v_n},-\overline{Q}_{p_n,k_n})-d(\overline{Q}_{p_n,k_n},\overline{Q}_{p,k})\\
&\geq   d(x, \Sigma_{v_n}) + \pi -\overline{r}_{k_n}-\delta/2-d(\overline{Q}_{p_n,k_n},\overline{Q}_{p,k})\\
&=  -d_{v_n}(x) + \pi -\overline{r}_{k_n}-\delta/2-d(\overline{Q}_{p_n,k_n},\overline{Q}_{p,k})\\
&>  -t_n + \pi -\overline{r}_{k_n}-\delta/2-d(\overline{Q}_{p_n,k_n},\overline{Q}_{p,k})\\
&> -t + \pi -\overline{r}_{k}-\delta.
\end{align*}
In any case, $d(x,\overline{Q}_{p,k}) <\overline{r}_{k} +t+ \delta$ and the claim follows.  

Arguing in the very same way, one can also show that
$$B_{\overline{r}_{k}+t-\delta}(\overline{Q}_{p,k})\subset A_{(v_n,t_n)} \quad\mbox{for all }n\geq n_2.$$
\begin{comment}
Suppose $x' \notin A_{(v_n,t_n)}$. Then $d_{v_n}(x') \geq t_n$, and $x'={\rm exp}_{v_n}(y',d_{v_n}(x'))$ for some $y'\in \Sigma_{v_n}$.

If $d_{v_n}(x') \geq 0$, then $x'\notin F_{v_n}(A)$. Hence  any continuous path joining $x'$ to   $\overline{Q}_{p,k}$ must intersect $\Sigma_{v_n}$. Thus we obtain from \eqref{delta.close.inclusion2}
\begin{eqnarray*}d(x',\overline{Q}_{p,k}) &\geq& d(x',\Sigma_{v_n}) -d(\Sigma_{v_n},\overline{Q}_{p_n,k_n})-\delta/2-d(\overline{Q}_{p_n,k_n},\overline{Q}_{p,k})\\
&\geq& d(x',\Sigma_{v_n}) + \overline{r}_{k_n}-\delta/2-d(\overline{Q}_{p_n,k_n},\overline{Q}_{p,k})\\
&=& d_{v_n}(x') + \overline{r}_{k_n}-\delta/2-d(\overline{Q}_{p_n,k_n},\overline{Q}_{p,k})\\
&\geq&   t_n + \overline{r}_{k_n}-\delta/2-d(\overline{Q}_{p_n,k_n},\overline{Q}_{p,k})\\
&>& \overline{r}_{k}+t - \delta.
\end{eqnarray*}

If $d_{v_n}(x') <0$, we obtain from \eqref{delta.close.inclusion}
\begin{eqnarray*}
d(x',\overline{Q}_{p,k})&\geq& d(y',\overline{Q}_{p,k}) -d(y',x') \\
&\geq&\overline{r}_{k_n}-\delta/2 +d_{v_n}(x') -d(\overline{Q}_{p_n,k_n},\overline{Q}_{p,k}) \\
&\geq& \overline{r}_{k_n}+ t_n-\delta/2-d(\overline{Q}_{p_n,k_n},\overline{Q}_{p,k})\\
&>&  \overline{r}_{k}+ t-\delta.
\end{eqnarray*}
In any case, $d(x',\overline{Q}_{p,k}) > \overline{r}_{k}+ t-\delta$ and the claim follows.
\end{comment}
Hence, for $n \geq n_2$,
$$
B_{\overline{r}_{k}+t-\delta}(\overline{Q}_{p,k})\subset A_{(v_n,t_n)}\subset B_{\overline{r}_{k}+t+\delta}(\overline{Q}_{p,k}).
$$
This implies  
$$
\left(A_{(v_n,t_n)}\,\Delta\, B_{\overline{r}_{k}+t}(\overline{Q}_{p,k})\right)\cup\Sigma_{(v_n,t_n)} \subset \overline B_{\overline{r}_{k}+t+\delta}(\overline{Q}_{p,k})\setminus B_{\overline{r}_{k}+t-\delta}(\overline{Q}_{p,k}).
$$
The result follows since $\delta>0$ can be chosen arbitrarily small.
\end{proof}

\subsection{Proof of Theorem \ref{canonical.family.continuous}} We start by arguing that the function $\overline r$ defined on $S^3 \cup \overline\Omega_{\varepsilon}$ is continuous. Clearly $\overline r$ is continuous on $S^3 \cap \overline{\Omega}_\varepsilon$, $A^* \setminus \overline{\Omega}_\varepsilon$, and $A \setminus \overline{\Omega}_\varepsilon$. Assume $$v=\Lambda(p,(0,t))=\cos t \,p+ \sin t \, N(p)\in \Omega_{2\varepsilon}.$$ The continuity follows at once from
$$\lim_{t\to\varepsilon_{-}}\overline r(v)=\lim_{k\to\infty}\left(\frac{\pi}{2}-\arctan(k)\right)=0$$
and
$$\lim_{t\to-\varepsilon_+}\overline r(v)=\lim_{k\to-\infty}\left(\frac{\pi}{2}-\arctan(k)\right)=\pi.$$

Consider the map
$$U: \overline B^4\times[-\pi,\pi]\rightarrow {\bf I}_3(S^3)$$
 \begin{eqnarray}\label{mapu}
 U(v,t)=\left\{
\begin{array}{rl}
[|A_{(T(v),t)}|]  & \mbox{if  } v\in B^4\setminus \overline\Omega_{\varepsilon},\\[0.5em] 
%\mbox{ } &\mbox{ }\\
 [|B_{\overline r(v)+t}(\overline Q(v))|] & \mbox{if  } v\in S^3\cup\overline \Omega_{\varepsilon}.\\
\end{array}
\right.
\end{eqnarray}
{Note that, by Theorem \ref{heintze.karcher},   $A_{(T(v),t)}$ has finite perimeter and so indeed $U(v,t)\in {\bf I}_3(S^3)$ for all $(v,t)\in\overline B^4\times[-\pi,\pi]$.}

\subsection{Lemma}\label{mass.continuity}\textit{
The map $U$ is continuous with respect to the mass topology of currents.}
\begin{proof}
We will use the fact that if $V_1,V_2 \subset S^3$ are open sets, then
$$
{\bf M}([|V_1|]-[|V2|]) = {\rm vol}\left(V_1\,\Delta\,V_2\right).
$$

Let $(v_n,t_n)$ tend to $(v,t)$ with $v_n,v\in B^4 \setminus \overline{\Omega}_\varepsilon$. 
Hence $(T(v_n),t_n)$ tends to $(T(v),t)$ and we obtain from Proposition {\ref{convergence.sets}} (i) that
$$\lim_{n\to\infty}{\bf M}(U(v_n,t_n)-U(v,t))=0.$$

Suppose now that $(v_n,t_n)$ tends to $(v,t)$ with $v\in A \setminus \overline{\Omega}_\varepsilon$.  We have $T(v_n)$ converging to  $T(v) \in A\setminus \Sigma$ and $\overline r(v)=\pi$. Thus $U(v,t)= [| B_{\pi+t}(T(v))|]$. For every $n$ sufficiently large
$$U(v_n,t_n)= [| B_{\pi+t_n}(T(v_n))|]\quad \mbox{if }v_n \in S^3,$$
or
$$U(v_n,t_n)=[|A_{(T(v_n),t_n)}|]\quad \mbox{if }v_n \in B^4.$$   
In any case, using Proposition {\ref{convergence.sets}} (ii), we get that 
$$\lim_{n\to\infty}{\bf M}(U(v_n,t_n)-U(v,t))=0.$$

The case $(v_n,t_n)$ tending to $(v,t)$ with $v\in A^* \setminus \overline{\Omega}_\varepsilon$ follows similarly, using
Proposition {\ref{convergence.sets}} (iii).

The restriction of $U$ to $\overline{\Omega}_\varepsilon$ is clearly continuous in the mass topology because $\overline Q$ and $\overline r$ are continuous functions.

It remains to consider the case $(v_n,t_n)$ converging to $(v,t)$ with $v_n \in B^4\setminus \overline{\Omega}_\varepsilon$ and $v\in \partial \Omega_\varepsilon$. We write $$v_n =\Lambda(p_n,s_n)\quad\mbox{and}\quad v =\Lambda(p,s),$$
where $\varepsilon=|s|<|s_n|,$  $s,s_n\in D^2_+(2\varepsilon),$
and set
$$\lim_{n\to\infty}\frac{{s_n}_2}{{s_n}_1}=\frac{s_2}{s_1}=k\in [-\infty,+\infty].$$
Recalling the definition of $T$ in Section \ref{associated.t} we have
 $$T(v_n)=\Lambda(p_n,u_n),\quad\mbox{where}\quad u_n=\phi(|s_n|)s_n,$$
 and so
 $$ \lim_{n\to\infty}\frac{{u_n}_2}{{u_n}_1}=k.$$
 Therefore, Proposition {\ref{convergence.sets}} (iv) implies that
\begin{multline*}
\lim_{n\to\infty}{\bf M}(U(v_n,t_n)-[|B_{\overline{r}_{k}+t}(\overline{Q}_{p,k})|])\\
=\lim_{n\to\infty}{\bf M}([|A_{(T(v_n),t_n)}|]-[|B_{\overline{r}_{k}+t}(\overline{Q}_{p,k})|]) =0.
\end{multline*}
We claim that $U(v,t)=[|B_{\overline{r}_{k}+t}(\overline{Q}_{p,k})|]$ and this implies the desired continuity at once.

Indeed, since $$|s|=\varepsilon\implies k=\frac{s_2}{\sqrt{\varepsilon^2-s_2^2}},$$ we see from the definition of $\overline Q$ in  \eqref{Q.R.} and $\overline r$ in \eqref{r.definition} that
$$\overline Q(v)=\overline Q_{p,k}\quad\mbox{and}\quad\overline r(v)=\overline r_k{.}$$
{This} implies $U(v,t)=[|B_{\overline{r}_{k}+t}(\overline{Q}_{p,k})|]$.
\end{proof}

{From the  Boundary Rectifiability Theorem (Theorem 30.3 of \cite{simon}) we know that  $C(v,t)=\partial U(v,t)\in \mathcal{Z}_2(S^3)$ and}  Lemma \ref{mass.continuity}  implies at once that  $C$ is continuous in the flat  topology. 

We now argue that 
$${\bf M}(C(v,t))\leq \mathcal{W}(\Sigma)\quad\mbox{for all}\quad(v,t)\in \overline B^4\times[-\pi,\pi].$$
This only needs justification if $v\in B^4\setminus \overline\Omega_{\varepsilon}$.

{If $\partial^*A_{(T(v),t)}\subset \Sigma_{(T(v),t)}$ denotes the reduced boundary  of $A_{(T(v),t)}$ (see \cite[Section 5.7]{evans} for the definition), we have from  \cite[Remark 27.7]{simon} and the Structure Theorem in \cite[p. 205]{evans} that
 $${\bf M}(C(v,t))=\mathcal{H}^2(\partial^*A_{(T(v),t)})\leq {\rm area}(\Sigma_{(T(v),t)}).$$
Theorem \ref{heintze.karcher}  then proves the desired inequality.}

We are left to prove the final statement of  Theorem \ref{canonical.family.continuous}. If $v\in S^3 \cup \overline{\Omega}_\varepsilon$, it is clear from \eqref{mapu} that  $U(v,\pi) = [|S^3|]$
and $U(v,-\pi) = 0$ and thus $C(v,\pm\pi) = \partial U(v,\pm\pi) = 0$.

If $v \in B^4\setminus \overline{\Omega}_\varepsilon$, set $w=T(v) \in B^4$. Since $\Sigma_w$ is 
a smooth surface, there can be no point $p\in S^3$ with $d(p,\Sigma_w)\geq\pi$ (otherwise $\Sigma_w\subset \{-p\}$). Therefore 
$A_{(w,\pi)} = S^3$, and $A_{(w,-\pi)}=\emptyset$. Since in this case $U(v,t) =[|A_{(w,t)}|]$, we again have 
$C(v,\pm\pi) = \partial U(v,\pm\pi) = 0$.

%%%%%%%%%%%%%%%%%%%%%%%%%%%%%%%%%%%%%%%%%%%%%%%%%%%%%%%%%%
%%%%%%%%%%%%%%%%%%%%%%%%%%%%%%%%%%%%%%%%%%%%%%%%%%%%%%%%%%%%%%

%%%%%%%%%%%%%%%%%%%%%%%%%%%%%%%%%%%%%%%%%%%%%%%%%%%%%
%%%%%%%%%%%%%%%%%%%%%%%%%%%%%%%%%%%%%%%%%%%%%%%%%%%%%%%%%%%

%%%%%%%%%%%%%%%%%%%%%%%%%%%%%%%%%%%%%%%%%%%%%%%%%%%%%%%%%%
%%%%%%%%%%%%%%%%%%%%%%%%%%%%%%%%%%%%%%%%%%%%%%%%%%%%%%%%%%%%%%%

\section{Min-max family}\label{minmax.family.section}

In this section we construct the  continuous map $\Phi$ into $\mathcal{Z}_2(S^3)$ to which we apply Almgren-Pitts Min-Max Theory. 

Recall the definition of the map $C$ in Section \ref{canonical.boundary.section}. From Theorem \ref{canonical.family.continuous} we can extend $C$ continuously to $\overline B^4\times \R$ by  defining  $C(v,t)=0$ when $|t|\geq \pi$. We denote this extension by $C$ as well.

We also choose a continuous extension  of $\overline r$, defined in \eqref{r.definition}, to a function $\overline r:\overline B^4\rightarrow [0,\pi]$.  

Choose an orientation preserving homeomorphism  $f:I^4 \rightarrow \overline{B}^4$ (hence   {$f_{|\partial I^4}$ is a homeomorphism from $\partial I^4$ onto $S^3$}) and consider $$\gamma:\R \rightarrow \R, \quad \gamma(s)=0\mbox{ if } s\leq \frac{1}{2}, \quad \gamma(s)=2s-1 \mbox{ if }s\geq \frac 1 2.$$

\subsection{Definition}\label{Fi.family} The {\bf min-max family} of $\Sigma$ is the map $\Phi:I^5 \rightarrow \mathcal{Z}_2(S^3)$ given by 
$$
 {\Phi(x,t) = C\left(f(x),2\pi\, (2t-1)+ \gamma (|f(x)|)\left(\frac \pi 2-\overline r\left(f(x)\right)\right)\right).}
$$
\subsection{Remark}\label{1/2.geodesic} The motivation for this definition is that if $x\in \partial I^4$ then we see from the definition of the map $C$ in Theorem \ref{canonical.family.continuous} that
\begin{multline*}
\Phi(x,t)=C\left(f(x),2\pi\, (2t-1)+\frac{\pi}{2}-\overline r(f(x))\right)\\
=\partial[|B_{2\pi\, (2t-1)+\pi/2}(\overline Q(f(x))|]
\end{multline*}

The properties of $\Phi$ that are important for our proof are collected in the next theorem. We denote by $\mathcal{T}\subset \mathcal{V}_2(S^3)$ the set of all  (unoriented) great spheres, which is  {homeomorphic} to $\RP^3$.  {The quantity  ${\bf m}(\Phi,r)$ appears in  Definition \ref{mass.def}.}
\subsection{Theorem}\label{modified.family}
{\em Let $\Sigma \subset S^3$ be an embedded closed surface of genus $g$. The map $$\Phi:I^5 \rightarrow \mathcal{Z}_2(S^3)$$ satisfies  the following properties:
\begin{enumerate}
\item [(i)] $\Phi$ is continuous with respect to the flat topology of currents;
\item [(ii)] $\Phi(I^4 \times \{0\})=\Phi(I^4 \times \{1\}) = \{0\};$
\item [(iii)] $$\sup\{{\bf M}(\Phi(x)):x\in I^5\}\leq \mathcal{W}(\Sigma);$$
\item[(iv)] the restriction $\Phi: \partial I^4 \times I \rightarrow\mathcal{Z}_2(S^3)$ is continuous in the ${\bf F}$-metric;
\item[(v)] for every $c\in I^4$, the map $\gamma:I\rightarrow\mathcal{Z}_2(S^3),\quad \gamma(t)=\Phi(c,t)$ is such that
\begin{itemize}
\item $\gamma(t)=\partial [|U(t)|]\quad\mbox{for all }0\leq t\leq 1$, where $U(t)$ are open sets of finite perimeter of $S^3$;
\item $U(0)=\emptyset$  and $U(1)=S^3$;
\item the map $t\to [|U(t)|]$ is continuous in the mass norm.
\end{itemize}
\item[(vi)]  {$$\max\{{\bf M}(\Phi(x)):x\in \partial I^5\}{=} 4\pi$$ and $$x\in \partial I^5\quad\mbox{and}\quad {{\bf M}(\Phi(x))=4\pi \Rightarrow |\Phi(x)|\in \mathcal{T}};$$} 
\item[(vii)]  {for every $\delta>0$ there is $\varepsilon>0$ so that, for all $(x,t)\in \partial I^5$},
$${\bf F}(|\Phi {(x,t)}|,\mathcal{T})\leq \varepsilon\implies |t-1/2|\leq \delta;$$
%\item[vii)] for every $\delta$ there is $\eta$ such that  the continuous map 
%$$\widehat \Phi:\partial I^4\times[1/2-\eta,1/2+\eta]\rightarrow \RP^3,\quad \widehat \Phi(x,t)=\partial B_{\pi/2}(\overline Q(f(x))$$
%satisfies
%$${\bf F}(|\Phi(x,t)|,\widehat \Phi(x,t))\leq \delta\quad\mbox{for all }(x,t)\in \partial I^4\times[1/2-\eta,1/2+\eta];$$
\item[(viii)] {the map $|\Phi|:\partial I^4\times\{1/2\}\rightarrow \mathcal{T}$ defined by
$$|\Phi|(x,1/2)= |\Phi(x,1/2)|=|\partial B_{\pi/2}(\overline Q(f(x)) | $$
has
$$|\Phi|_{*}([\partial I^4\times\{1/2\}])|=2g\in H_3(\RP^3,\Z);$$ }
\item [(ix)] {$\lim_{r\to 0}{\bf m}(\Phi,r)=0.$}
%\item[x)] if a sequence $\{y_j\}_{j\in \N}$ in $I^5$ tends to $\partial I^4\times I$, then
%$$
%\liminf_{j\to\infty} {\bf M}(\Phi(y_j))>0 
%\implies \liminf_{j\to\infty} {\mathcal F}(\Phi(y_j))>0.$$
\end{enumerate}
}

\begin{proof}
Property (i) comes from the fact that $\Phi$ is a composition of continuous functions. 

Because $0\leq \overline r\leq \pi$, we have from Theorem \ref{canonical.family.continuous} that 
$$\Phi(x,1)=C(f(x),2\pi+\gamma(|f(x)|)(\frac \pi 2-\overline r(f(x))))=0\quad\mbox{for all }x\in I^4.$$
Likewise, $\Phi(x,0)=0$ and this shows property (ii).

Property (iii) follows from the {mass estimate in Theorem \ref{canonical.family.continuous}.}

From the definition of $C$ it is clear that  $C$ restricted to $S^3\times[-\pi,\pi]$ is continuous in the ${\bf F}$-metric and thus $\Phi$ restricted to $\partial I^4 \times I$ is also continuous in the ${\bf F}$-metric. This proves property (iv).

Property (v) follows at once from the fact that $C(v,t)=\partial U(v,t)$, where the map $U$ is defined in \eqref{mapu}, and from Lemma \ref{mass.continuity}.

 %{From the definition of $\Phi$ we have that $\Phi(x,1)=\Phi(x,0)=0$ for all $x\in \partial I^4$.}
 Property (vi) follows  from Remark  \ref{1/2.geodesic}.
Property (vii) follows from property (iv) and  {the fact that, from Remark \ref{1/2.geodesic}, for every $x\in \partial I^4$
$$|\Phi(x,t)|\in \mathcal{T} \iff  t=1/2.$$}

Consider the $2$-fold cover of $\mathcal{T}$ given by $$\pi:S^3\rightarrow \mathcal{T},\quad \pi(p)= {|\partial B_{\pi/2}(p)|}.$$
We have $ {|\Phi|(x ,1/2)}=\pi\circ \overline Q\circ f(x)$, and so Theorem \ref{degree.gauss.map} implies the degree of $x\mapsto  {|\Phi|(x ,1/2)}$ is $2g$. This implies property (viii).

Property (ix) is a consequence of Theorem \ref{no.concentration.mass}.

%Finally, we prove property (x). Assume we have a sequence $\{(y_n)\}_{n\in\N}$ in $I^5$ tending to $y\in \partial I^4\times I$ and with
%\begin{equation}\label{area.big.delta}
%\liminf_{n\to\infty}{\bf M}(\Phi(y_n))\geq 2\delta>0.
%\end{equation}
%Suppose that
%\begin{equation}\label{flat.zero.phi}
%\liminf_{n\to\infty}{\mathcal F}(\Phi(y_n)=0.
%\end{equation}
%We have $\Phi(y_n)=C(v_n,t_n)$, where $(v_n,t_n)$ tends to $S^3\times [-\pi,\pi]$. From \eqref{flat.zero.phi} we obtain $\mathcal{F}(C(v,t))=0$ and thus, because $v\in S^3$,  Proposition \ref{convergence.sets} implies that  $C(v_n,t_n)$ converges in Hausdorff distance to a point. Therefore, we obtain from \eqref{area.big.delta} that $\liminf_{r\to 0}{\bf m}(\Phi,r)\geq 2\delta$. This contradicts property (ix).

\end{proof}

%%%%%%%%%%%%%%%%%%%%%%%%%%%%%%%%%%%%%%%%%%%%%%
%%%%%%%%%%%%%%%%%%%%%%%%%%%%%%%%%%%%%%%%%%%%%%%%%%%

\section{Almgren-Pitts Min-Max Theory I }\label{gmt.definitions}

{We will set up the notation needed to apply the Almgren-Pitts Min-Max Theory to our setting.
$(M,g)$ will denote an orientable compact Riemannian three-manifold.}

%We  assume $M$ is
%isometrically embedded in $\mathbb{R}^L$. We denote by $B_r(p)$ the open geodesic
%ball in $M$ of radius $r$ and center $p\in M$.

\subsection{Cell complexes} 

We denote by  $I^n=[0,1]^n\subset \R^n$  the $n$-dimensional cube, with boundary $I^n_0=\partial I^n=I^n\setminus (0,1)^n$. 

For each $j\in \N$, $I(1,j)$ denotes the cell complex on $I^1$  whose $1$-cells and $0$-cells (those are sometimes called vertices) are, respectively,  
$$[0,3^{-j}], [3^{-j},2 \cdot 3^{-j}],\ldots,[1-3^{-j}, 1]\quad\mbox{and}\quad [0], [3^{-j}],\ldots,[1-3^{-j}], [1].$$

We consider the $n$-dimensional cell complex on $I^n$: 
$$I(n,j)=I(1,j)\otimes\ldots \otimes I(1,j)\quad (\mbox{$n$ times}).$$
$\alpha=\alpha_1 \otimes \cdots\otimes \alpha_n$ is a $p$-cell of $I(n,j)$ if and only if $\alpha_i$ is a cell
of $I(1,j)$ for each $i$, and $\sum_{i=1}^n {\rm dim}(\alpha_i) =p$. We often abuse notation by identifying  a $p$-cell $\alpha$ with its support: $\alpha_1 \times \cdots \times \alpha_n \subset I^n$. 
 
  {We use the following notation: 
\begin{itemize}
\item $I(n,j)_p$ denotes the set of all $p$-cells in $I(n,j)$;
\item $I_0(n,j)_p$ denotes the set of all $p$-cells  of $I(n,j)$ which are contained in the boundary  $I^n_0$;
\item $I_0(n,j)$ is the subcomplex of $I(n,j)$ generated by  all cells that are contained in the boundary $I^n_0$.
\end{itemize}
Given  a $p$-cell $\alpha\in I(n,j)_p$ we use the notation:
\begin{itemize}
\item$\alpha(0)$ denotes the $p$-dimensional subcomplex of $I(n,j)$ whose  cells are those with support contained in $\alpha$;
\item  $\alpha(k)$ denotes the $p$-dimensional subcomplex of $I(n,j+k)$ formed by  all cells 
that are contained in $\alpha$;
\item $\alpha(k)_q$, with $q\leq p$, denotes the set of all $q$-dimensional cells of $\alpha(k)$;
\item $\alpha_0(k)_q$, with $q\leq p$, denotes  the set of all $q$-dimensional cells of $\alpha(k)$ whose support is contained in the boundary of $\alpha$.
\item $\alpha_q=\alpha(0)_q$ denotes the $q$-dimensional faces of $\alpha$.
\end{itemize} 
 }

We also define the following cell subcomplexes of $I(n,j)$:
\begin{align*}
(\mbox{top})&\quad T(n,j)=I(n-1,j) \otimes \langle[1]\rangle,\\
  (\mbox{bottom})&\quad B(n,j)=I(n-1,j)\otimes\langle[0]\rangle,\\
(\mbox{side})&\quad S(n,j)=I_0(n-1,j) \otimes I(1,j).
\end{align*}
(Here $\langle[x]\rangle$ is the cell complex whose only cell is  $[x]$.)  Let $T(n,j)_p,  B(n,j)_p$, and $S(n,j)_p$  be the corresponding sets of $p$-cells. Note that
$T(1,j)=\langle[1]\rangle$, and $B(1,j)=\langle[0]\rangle$.

The boundary homomorphism
$$\partial:I(n,j)\rightarrow I(n,j)$$
is defined  by 
$$\partial (\theta^1\otimes\ldots\otimes\theta^n)=\sum_{i=1}^n(-1)^{\sigma(i)}\theta^1\otimes\ldots\otimes\partial\theta^i\otimes\ldots\otimes\theta^n,$$
where
$$\sigma(i)=\sum_{p<i}\mathrm{dim}(\theta^p)$$
and $$\partial([a,b])=[b]-[a]\mbox{ if }[a,b]\in I(1,j)_1, \quad \partial([a])=0\mbox{ if }[a]\in I(1,j)_0.$$

The distance between two vertices of $I(n,j)$ is defined by
$${\bf d}: I(n,j)_0\times I(n,j)_0 \rightarrow \mathbb{Z}_+, \quad {\bf d}(x,y)=3^j\sum_{i=1}^n|x_i-y_i|.$$
It has the  property  that two vertices $x,y$ satisfy ${\bf d}(x,y)=1$ if and only if $[x,y]$  is a 1-cell of $I(n,j)$. 

We will also need the map ${\bf n}(i,j):I(n,i)_0\rightarrow I(n,j)_0$, defined as follows: for each  $x\in I(n,i)_0$, ${\bf n}(i,j)(x)$ is the unique element of $I(n,j)_0$ such that
$${\bf d}(x,{\bf n}(i,j)(x))=\inf\{{\bf d}(x,y):y\in I(n,j)_0\}.$$
Note that ${\bf n}(i,j)(x)=x$ if $i\leq j$, and ${\bf n}(k,i)={\bf n}(j,i) \circ {\bf n}(k,j)$ if $i\leq j\leq k$.

\subsection{Maps into currents}\label{fineness.defi}
Given a map $\phi:I(n,j)_0\rightarrow  \mathcal{Z}_2(M)$, we define the {\em fineness} of $\phi$ to be
$${\bf f}(\phi)=\sup\left\{\frac{{\bf M}(\phi(x)-\phi(y))}{{\bf d}(x,y)}: x,y\in I(n,j)_0, x\neq y\right\}.$$

The reader should think of the notion of fineness as being a discrete  measure of continuity with respect to the mass norm. The following lemma is useful for computational purposes.
 
\subsection{Lemma}\label{lemma_adjacent}\textit{
 ${\bf f}(\phi)<\delta$ if and only if  ${\bf M}(\phi(x)-\phi(y))< \delta$ whenever ${\bf d}(x,y)=1$.
}

\begin{proof} If ${\bf f}(\phi)<\delta$, then it follows directly from the definition of fineness that  ${\bf M}(\phi(x)-\phi(y))< \delta$ whenever ${\bf d}(x,y)=1$. Suppose now that ${\bf M}(\phi(x)-\phi(y))< \delta$ if ${\bf d}(x,y)=1$. Given any $x,y\in I(n,j)_0$ with ${\bf d}(x,y)=k$,  we can find a  sequence $\{y_i\}_{i=0}^k$ in $I(n,j)_0$ so that $y_0=y$, $y_k=x$, and $[y_i,y_{i+1}]$ is a 1-cell of $I(n,j)$.   Thus
$$\frac{{\bf M}(\phi(x)-\phi(y))}{{\bf d}(x,y)}\leq \frac 1 k \sum_{i=1}^{k}{\bf M}(\phi(y_{i})-\phi(y_{i-1}))<\frac 1 k k\delta=\delta.$$
\end{proof}

%It is also convenient to measure the local concentration of mass of continuous maps $\Phi:I^n
%\rightarrow {\mathcal Z}_2(M)$:

\subsection{Homotopy notions}\label{homotopy} 
Suppose we have  a map
$$\Phi_0:\partial I^n\rightarrow {\mathcal Z}_2(M)$$
that satisfies
\begin{itemize}
\item $\Phi_0$ is continuous in the ${\bf F}$-metric;
\item 
$$\Phi_0(I^{n-1}\times\{0\})=\Phi_0(I^{n-1}\times\{1\})=0.$$
\end{itemize}

Let $\phi_i:I(n,k_i)_0\rightarrow  \mathcal{Z}_2(M)$, $i=1,2$. We say that $\phi_1$ is {\it $n$-homotopic to $\phi_2$ in $(\mathcal{Z}_2(M;{\bf M}),\Phi_0)$ with fineness $\delta$}  if we can find  $k\in \N$ and a map
$$\psi: I(1,k)_0\times I(n,k)_0\rightarrow  \mathcal{Z}_2(M)$$
such that 
\begin{itemize}
\item[(i)] ${\bf f}(\psi)<\delta;$
\item[(ii)] if $i=1,2$ and $x\in I(n,k)_0$, then
$$\psi([i-1],x)=\phi_i({\bf n}(k,k_i)(x));$$
\item[(iii)] $$\psi(I(1,k)_0\times T(n,k)_0)=\psi(I(1,k)_0\times B(n,k)_0)=\{0\};$$
\item[(iv)] 
\begin{eqnarray*}
&&\sup\,\{{\mathcal{F}}(\psi(t,x)-\Phi_0(x)):(t,x)\in I(1,k)_0\times S(n,k)_0\} \leq \delta,\\
&&{\bf M}(\psi(t,x))\leq {\bf M}(\Phi_0(x))+\delta  {\rm \ for \  any \ } (t,x)\in I(1,k)_0\times S(n,k)_0.
%\sup\,\{{\bf M}(\psi(t,x))-{\bf M}(\Phi_0(x)):(t,x)\in I(1,k)_0\times S(n,k)_0\} &\leq& \delta.
\end{eqnarray*}
\end{itemize}
In particular we must have that $\phi_i=0$  on $T(n,k_i)_0\cup B(n,k_i)_0,$ 
$$\sup\{{\mathcal{F}}(\phi_i(x)-\Phi_0(x)):x\in  S(n,k_i)_0\} \leq \delta,$$
 and
$$\sup\{{\bf M}(\phi_i(x))-{\bf M}(\Phi_0(x)):x\in  S(n,k_i)_0\} \leq \delta,$$
for each $i=1,2$.

We note that if $\phi_1$ is homotopic to $\phi_2$ with fineness $\delta_1$, and $\phi_2$ is homotopic to $\phi_3$ with fineness $\delta_2$, then $\phi_1$ is homotopic to $\phi_3$ with fineness $\max\{\delta_1,\delta_2\}.$

\subsection{Remark}
It is convenient to compare the above definition  with a related definition used by  Pitts \cite[Section 4.1]{pitts}:
 $\phi_1$ is $n$-homotopic to $\phi_2$ in $(\mathcal{Z}_2(M;{\bf M}),\{0\})$ with fineness $\delta$, according to Pitts,  if 
 we can find  $k\in \N$ and a map
$$\psi: I(1,k)_0\times I(n,k)_0\rightarrow  \mathcal{Z}_2(M)$$
such that 
\begin{itemize}
\item[(i)] ${\bf f}(\psi)<\delta;$
\item[(ii)] if $i=1,2$ and $x\in I(n,k)_0$, then
$$\psi([i-1],x)=\phi_i({\bf n}(k,k_i)(x));$$
\item[(iii)]$$ \psi(S(n+1,k)_0)=\{0\}.$$
\end{itemize}
Note that for the definition of Pitts to make sense, it is required that  $\phi_i(I_0(n,k_i)_0)=\{0\}$ for each $i=1,2$. In the one-dimensional case ($n=1, \Phi_0 =0$), our notion is equivalent to the
definition of Pitts. 

\medskip

Instead of considering continuous maps from $I^n$ into $\mathcal{Z}_2(M;{\bf M})$, Almgren-Pitts consider sequences of discrete maps into  $\mathcal{Z}_2(M)$ with fineness tending to zero. 

\subsection{Definition}\label{homotopy.sequence.phi}
 An $$\mbox{{\it $(n,{\bf M})$-homotopy sequence of mappings into $(\mathcal{Z}_2(M;{\bf M}),\Phi_0)$}}$$ is a sequence of mappings $\{\phi_i\}_{i\in \N}$,
$$\phi_i:I(n,k_i)_0\rightarrow \mathcal{Z}_2(M),$$
such that $\phi_i$ is $n$-homotopic to $\phi_{i+1}$ in  $(\mathcal{Z}_2(M;{\bf M}),\Phi_0)$ with fineness $\delta_i$ and
\begin{itemize}
\item[(i)] $\lim_{i\to\infty} \delta_i=0$;
\item[(ii)]$\sup\{{\bf M}(\phi_i(x)):x\in I(n,k_i)_0, i\in \N\}<+\infty.$
\end{itemize}

\subsection{Remark}
This is similar to the notion of an $$\mbox{{\em $(n,{\bf M})$-homotopy sequence of mappings into $(\mathcal{Z}_2(M;{\bf M}),\{0\})$}}$$ in  \cite[Section 4.1]{pitts}. Both notions coincide in the one-dimensional case ($n=1, \Phi_0=0$).

 {The next lemma says that $\phi_i$  restricted to the boundary of its domain tends to $\Phi_0$ in the ${\bf F}$-metric.}

\subsection{Lemma} \label{homotopy.sequence.boundary}\textit{Let $S=\{\phi_i\}_{i\in \N}$ be an $(n,{\bf M})$-homotopy sequence of mappings into  $(\mathcal{Z}_2(M;{\bf M}),\Phi_0)$. If $I(n,k_i)_0$ denotes the domain of $\phi_i$, then
$$
\lim_{i\rightarrow \infty} \sup \{{\bf F}(\phi_i(x),\Phi_0(x)): x \in I_0(n, k_i)_0\} = 0.
$$
}
\begin{proof}
First note that $\phi_i(x)=\Phi_0(x)=0$  for $x \in T(n,k_i)_0\cup B(n,k_i)_0.$ Since $\Phi_0$ is continuous in the {\bf F}-metric,  $\Phi_0(I^n_0)$ is a compact subset of $\mathcal{Z}_2(M,{\bf F})$. The lemma follows from condition (iv) in the definition of ``homotopic to'', by using Lemma \ref{flat+mass=f}.
\end{proof}

The  next definition explains what  it means for two distinct homotopy sequences  of mappings into  $(\mathcal{Z}_2(M;{\bf M}),\Phi_0)$ to  be homotopic.

\subsection{Definition}
 Given $S^1=\{\phi^1_i\}_{i\in \N}$ and $S^2=\{\phi^2_i\}_{i\in \N}$  $(n,{\bf M})$-homotopy sequences  of mappings into  $(\mathcal{Z}_2(M;{\bf M}),\Phi_0)$, we say  that {\it $S^1$ is homotopic with $S^2$} if there exists $\{\delta_i\}_{i\in \N}$ such that
 \begin{itemize}
\item $\phi^1_i$  is $n$-homotopic to $\phi^2_i$ in  $(\mathcal{Z}_2(M;{\bf M}),\Phi_0)$ with fineness $\delta_i$;
\item $\lim_{i\to\infty} \delta_i=0.$
 \end{itemize}

\subsection{Remark}
There is a similar definition for $(n,{\bf M})$-homotopy sequences  of mappings into 
$(\mathcal{Z}_2(M;{\bf M}),\{0\})$ \cite[Section 4.1]{pitts}. Once again these definitions coincide
in the one-dimensional case  ($n=1, \Phi_0=0$).

\medskip

The relation ``is homotopic with'' is an equivalence relation on the set of all $(n,{\bf M})$-homotopy sequences  of mappings into  $(\mathcal{Z}_2(M;{\bf M}),\Phi_0)$. We call  the equivalence class 
of any such sequence  an {\it $(n,{\bf M})$-homotopy class of mappings into $(\mathcal{Z}_2(M;{\bf M}),\Phi_0)$}. We denote 
by $\pi_n^{\#}(\mathcal{Z}_2(M;{\bf M}),\Phi_0)$  the set of all equivalence classes. 

 {Finally, a  $$\mbox{{\em $(n,\mathcal{F})$-homotopy sequence  (or class) of mappings into $(\mathcal{Z}_2(M;\mathcal{F})),\{0\})$}}$$  is defined similarly to what  we just did   but with the mass  ${\bf M}$ in the definition of ${\bf f}$  being replaced by the flat metric $\mathcal{F}$. The set of all equivalence classes is denoted by $\pi_{n}^{\#}(\mathcal{Z}_2(M;\mathcal{F}),\{0\})$.}
  In \cite[Section 4.1]{pitts} (see also \cite[Section 3]{almgren}) it is also considered $\pi_1(\mathcal{Z}_2(M;\mathcal{F}),\{0\})$ to be the usual homotopy group  of equivalence classes of continuous mappings $(I,I_0)\rightarrow (\mathcal{Z}_2(M;\mathcal{F}),\{0\}$).
\subsection{Min-Max definitions}
Given  $\Pi \in  \pi_n^{\#}(\mathcal{Z}_2(M;{\bf M}),\Phi_0)$, let
$$ {\bf L}: \Pi\rightarrow [0,+\infty]$$
be defined by
$${\bf L}(S)=\limsup_{i\to\infty}\max\{{\bf M}(\phi_i(x)):x\in \mathrm{dmn}(\phi_i)\},\quad\mbox{where }S=\{\phi_i\}_{i\in \N}.$$
Note that ${\bf L}(S)$ is the discrete replacement for  the maximum area of a continuous map into $\mathcal{Z}_2(M;{\bf M})$.

\subsection{Definition}The {\it width} of $\Pi$ is defined by
$${\bf L}(\Pi)=\inf\{{\bf L}(S):S\in \Pi\}.$$

We also consider 
$${\bf K}:\Pi\rightarrow \{K:K \mbox{ compact subset of }\mathcal{V}_2(M)\},$$ given by  
\begin{multline*}
{\bf K}(S)=\{V:V=\lim_{j\to\infty}|\phi_{i_j}(x_j)|\mbox{ as varifolds, for some increasing}\\
\mbox{sequence }\{i_j\}_{j\in \N} \mbox{ and }x_j\in \mathrm{dmn}(\phi_{i_j})\}
\end{multline*}
 for $S=\{\phi_i\}_{i\in \N}\in \Pi$.

We say  $S\in \Pi$ is a {\it critical sequence} for $\Pi$ if $${\bf L}(S)={\bf L}(\Pi).$$  The {\em critical set}  ${\bf C}(S)$ of a critical sequence   $S\in \Pi$ is given by 
$${\bf C}(S)={\bf K}(S)\cap\{V: ||V||(M)={\bf L}(S)\}.$$
The set ${\bf C}(S)\subset \mathcal{V}_2(M)$ is nonempty and compact.

%%%%%%%%%%%%%%%%%%%%%%%%%%%%%%%%%%%%%%%%%%%%%%%%%%%%%%%%
%%%%%%%%%%%%%%%%%%%%%%%%%%%%%%%%%%%%%%%%%%%%%%%%%%%%%%%%%%%

%%%%%%%%%%%%%%%%%%%%%%%%%%%%%%%%%%%%%%%%%%%%%%%%%%%%%%%%%%%%
%%%%%%%%%%%%%%%%%%%%%%%%%%%%%%%%%%%%%%%%%%%%%%%%%%%%%%%%%%%%

\section{Almgren-Pitts Min-Max Theory II}\label{almgren.pitts.section}

{In our setting, the Almgren-Pitts Min-Max Theory applies to elements of  $\pi_n^{\#}(\mathcal{Z}_2(M;{\bf M}),\Phi_0)$. Therefore it is important to generate an $(n,{\bf M})$-homotopy sequence of mappings
 into $(\mathcal{Z}_2(M;{\bf M}),\Phi_0)$ out of    a continuous map 
$\Phi:I^{n}\rightarrow   \mathcal{Z}_2(M)$  in the flat topology.   This is the content of Theorem \ref{discrete.sweepout} below. In this section we also discuss Pitts Min-Max Theorem.}

Let $$c=\frac{1}{3}(1,\ldots,1,0)\in I^{n-1}\times\{0\},$$ and $e_{n}$ be the coordinate vector corresponding to the  $x_{n}$-axis.

 {We consider the following hypotheses for the continuous map in the flat topology $\Phi:I^{n}\rightarrow   \mathcal{Z}_2(M)$.
 
\begin{enumerate}
 \item[($A_0$)] $\Phi_{|I_0^n} \mbox{ is continuous in the ${\bf F}$-metric.}$
  \item[($A_1$)] $\Phi(I^{n-1}\times\{0\})=\Phi(I^{n-1}\times\{1\})=0.$
\item[($A_2$)] ${\bf L}(\Phi)=\sup\{{\bf M}(\Phi(x)):x\in I^{n}\}<+\infty.$
\item[($A_3$)] $\lim_{r\to 0}{\bf m}(\Phi,r)=0$ (recall Definition \ref{mass.def}).
\item[($A_4$)] The map $t\mapsto \Phi(c+tx_n)$, $0\leq t\leq 1$, defines a non-trivial class in $\pi_1(\mathcal{Z}_2(M;\mathcal{F}),\{0\}).$
% there exists a family $\{U(t)\}_{0\leq t\leq 1}$ of open sets of $M$ of finite perimeter with
%\begin{itemize}
%\item $U(0)=\emptyset$  and $U(1)=M$;
%\item $ \Phi(c+tx_n)=\partial [|U(t)|]\quad\mbox{for all }0\leq t\leq 1;$
%where $c=\frac{1}{3}(1,\ldots,1,0)\in I^{n-1}\times\{0\};$
%\item the map $t\to [|U(t)|]$ is continuous in the mass norm.
%\end{itemize}
\end{enumerate}
}

The next lemma assures that the min-max family $\Phi$  associated to  an embedded closed surface $\Sigma$ of $S^3$ satisfies the conditions above.

{
\subsection{Lemma} \em{Let $\Phi$ be the min-max family defined in Definition \ref{Fi.family}. Then $\Phi$ satisfies hypotheses $(A_0)$--$(A_4)$.}
\begin{proof} From Theorem \ref{modified.family} it is clear that hypotheses $(A_0)$--$(A_3)$ are satisfied.

Let $[\gamma]\in \pi_1(\mathcal{Z}_2(S^3;\mathcal{F}),\{0\})$ be the class generated by the map $\gamma(t)=\Phi(c+tx_n)$, $0\leq t\leq 1$.

For each $i$ is sufficiently large, Corollary 1.14 of Almgren  \cite{almgren} guarantees  the existence, for each $x\in I(1,i)_0\setminus \{[1]\}$, of $A_i(x)\in {\bf I}_3(S^3)$  so that
\begin{equation}\label{mass.small.homotopy}
 \partial A_i(x)=\gamma(x+3^{-i})-\gamma(x)\quad\mbox{and }\quad{\bf M}(A_i(x))={\mathcal{F}}(\partial A_i(x)).
 \end{equation}
 
If $F:\pi_1(\mathcal{Z}_2(S^3;\mathcal{F}),\{0\}) \rightarrow H_3(S^3,\mathbb{Z})$ is the natural isomorphism constructed by Almgren in Section 3 of \cite{almgren} (see also Theorem 13.4 of \cite{almgren-varifolds}), then   
$$F[\gamma]=\left[\sum_{j=0}^{3^{i}-1}A_i(j3^{-i})\right] \in H_3(S^3,\Z),$$
for every $i$ is sufficiently large. 

We now argue that $F[\gamma]=\big[S^3\big] \in H_3(S^3,\Z)$ and so condition  $(A_4)$ is also satisfied. 

From Theorem \ref{modified.family} (v) We know that 
$$\gamma(x+3^{-i})-\gamma_i(x)=\partial ([|U(x+3^{-i})|]- [|U(x)|]).$$
Thus $$B(x)=[|U(x+3^{-k_i})|]- [|U(x)|]-A_i(x)\in {\bf I}_3(S^3)$$ satisfies $\partial B(x)=0$. The Constancy Theorem (see \cite{simon}) then implies that $B(x)=k[|M|]$ for some  $k=k(x)\in \mathbb{Z}$. On the other hand, the continuity of $t\to [|U(t)|]$ in the mass norm, together with continuity of $\gamma$ and \eqref{mass.small.homotopy}, implies
that the mass of $B(x)$ becomes uniformly and arbitrarily small as   $i\rightarrow \infty.$ We conclude that if $i$ is sufficiently large then  $B(x)=0$ for all $x\in I(1,k_i)_0$. 

Therefore, for large $i$,
\begin{eqnarray*}
F[\gamma]&=&\left[\sum_{j=0}^{3^{i}-1}\left([|U((j+1)3^{-i})|]- [|U(j3^{-i})|]\right)\right]\\
&=&\big[[|U(1)|]-[|U(0)|]\big]=\big[S^3\big] \in H_3(S^3,\Z).
\end{eqnarray*}
\end{proof}
}

Then:
\subsection{Theorem }\label{discrete.sweepout}{\em  {Assume $\Phi$ satisfies hypotheses $(A_0)$--$(A_4)$}.

There exists an $(n,{\bf M})$-homotopy sequence of mappings into  $(\mathcal{Z}_2(M;{\bf M}),\Phi_{|I_0^n})$ $$\tilde{\phi}_i:I(n,k_i)_{0}\rightarrow  \mathcal{Z}_2(M),$$
 with the following properties:
\begin{itemize}
\item[(i)]  {There is a sequence $\{l_i\}_{i\in \N}$ tending to infinity such that for every sequence $x_i\in I(n,k_i)_0$ we have
$$\limsup_{i\to\infty}{\bf M}(\tilde\phi_i(x_i))\leq \limsup_{i\to\infty}\{{\bf M}(\Phi(x)):\alpha\in I(n,l_i)_n, x,x_i\in \alpha\}.$$
 In particular
  $${\bf L}(\{\tilde{\phi}_i\}_{i\in\N})\leq \sup\{{\bf M}(\Phi(x)):x\in I^{n}\};$$

 }
 \item[(ii)] {$$\lim_{i\to\infty}\sup\{{\mathcal{F}}(\tilde\phi_i(x)-\Phi(x))\,|\, x\in I(n,k_i)_0\}=0;$$}
\item[(iii)] The sequence of mappings 
$$v_i:I(1,k_i)_0\rightarrow  \mathcal{Z}_2(M;{\bf M}),\quad v_i(x)=\tilde{\phi}_i(c+xe_{n}),$$
 is a $(1,{\bf M})$-homotopy sequence of mappings into $(\mathcal{Z}_2(M;{\bf M}),\{0\})$  that belongs
 to a non-trivial element of $\pi_1^{\#}(\mathcal{Z}_2(M;{\bf M}),\{0\})$.
\end{itemize}
}
{The proof of Theorem \ref{discrete.sweepout} is {postponed to}  Section \ref{continuous.discrete}.}

%We note that if $\Phi$ is the min-max family associated to  an embedded closed surface $\Sigma$ of $S^3$, we have from Theorem  \ref{modified.family} that {$\Phi$  hypotheses $(A_0)$--$(A_4)$ and thus we can apply Theorem \ref{discrete.sweepout}.}

\subsection{Definition}\label{homotopy.class.sigma} Let $\Sigma$ be an embedded closed surface in $S^3$ and let $\Phi$ be the min-max family associated to $\Sigma$ constructed in Section \ref{minmax.family.section}.  The {\em homotopy class associated with $\Sigma$} is defined to  be the homotopy class of $S=\{\tilde \phi_{i}\}_{i\in\N}$ given by Theorem \ref{discrete.sweepout} applied to $\Phi$.

\subsection{Min-Max Theorem}
{We now} adapt the celebrated  {Pitts  Min-Max Theorem} to our setting.   
Assume we have a continuous map in the flat topology
$$\Phi:I^{n}\rightarrow   \mathcal{Z}_2(M)$$
 which satisfies the hypotheses $(A_0)-(A_1)$.
We denote by   $|\Phi|:I^n \rightarrow \mathcal{V}_2(M)$ the map given by $|\Phi|(x)=|\Phi(x)|$ for $x\in I^n$.

Consider $\Pi \in  \pi_n^{\#}(\mathcal{Z}_2(M;{\bf M}),\Phi_{|I_0^n})$.  

\subsection{Proposition}\label{pulltight}\textit{
There exists  a critical sequence $S^* \in \Pi$. For each  critical sequence $S^*$, there exists a critical sequence $S\in\Pi$ such that
\begin{itemize}
\item ${\bf C}(S)\subset {\bf C}(S^*)$;
\item every $\Sigma\in {\bf C}(S)$ is either a stationary varifold or belongs to $|\Phi|(I^n_0)$.
\end{itemize}
}
{The sequence $S$ is obtained from a pull-tight procedure applied to  $S^*$. The proof follows very closely Theorem 4.3 of \cite{pitts}  and is postponed to Section \ref{proof.pulltight}}.

One  consequence of Proposition \ref{pulltight} is the following  theorem,  established by Pitts  \cite{pitts} when $\Pi$ is a non-trivial element of $\pi_n^{\#}(\mathcal{Z}_2(M;{\bf M}),\{0\})$.  The proof follows by simple adaptation
of the arguments in \cite{pitts}.

\subsection{Theorem}\label{pitts.min.max} {\em  {Assume $\Phi$ satisfies $(A_0)-(A_1)$}.

Let $\Pi \in  \pi_n^{\#}(\mathcal{Z}_2(M;{\bf M}),\Phi_{|I_0^n})$ with
$$\max\{{\bf M}(\Phi(x)):x\in I^n_0\}<{\bf L}(\Pi)<\infty.$$
There exists   a stationary integral varifold $\Sigma$, whose  support is a smooth embedded minimal surface, such that
$$||\Sigma||(M)={\bf L}(\Pi).$$ 
Moreover, if $S^*$ is a critical  sequence   then we can choose $\Sigma \in {\bf C}(S^*)$.
}
\begin{proof}
Consider $S=\{\varphi_i\}_{i\in \N}\in \Pi$ given by Proposition \ref{pulltight}, and let
$$
0< \varepsilon = {\bf L}(S) - \max\{{\bf M}(\Phi(x)):x\in I^n_0\}.
$$ Because every $\Sigma\in {\bf C}(S)$ satisfies $$||\Sigma||(M)={\bf L}(\Pi)>\max\{{\bf M}(\Phi(x)):x\in I^n_0\},$$ we obtain that  every $\Sigma$ in ${\bf C}(S)$ must be stationary.  Since the construction of 
\cite[Theorem 4.10]{pitts} can be made to not affect those $\varphi_i(x)$ with $${\bf M}(\varphi_i(x))\leq {\bf L}(S)-\varepsilon/2,$$ and since $${\bf M} (\varphi_i(x))\leq \max\{{\bf M}(\Phi(x)):x\in I^n_0\} + \varepsilon/2$$ for every $x \in {\rm dmn}(\varphi_i) \cap I_0^n$ and sufficiently large $i$, we can see that the competitor $\{\varphi_i^*\}_{i\in \N}$ constructed by Pitts belongs to $\Pi$. Therefore, as in \cite{pitts}, we can find
an almost-minimizing (in annular regions)  $\Sigma\in {\bf C}(S)$.  The regularity theory developed in \cite[Section 7]{pitts} implies  that $\Sigma$ is an integral varifold whose support is a smooth embedded minimal surface.
\end{proof}

%%%%%%%%%%%%%%%%%%%%%%%%%%%%%%%%%%%%%%%%%%%%%%%%%%%%%%%%
%%%%%%%%%%%%%%%%%%%%%%%%%%%%%%%%%%%%%%%%%%%%%%%%%%%%%%%%%%%

\section{Lower bound on width}\label{bound.width.section}

Let ${\mathcal T}\subset {\mathcal V}_2(S^3)$  be the set of all varifolds that correspond to a  great sphere in $S^3$
with multiplicity one. Note that  $\mathcal{T}$ is naturally homeomorphic to $\RP^3$.  

Let
 $$\Phi:I^{5}\rightarrow   \mathcal{Z}_2(S^3)$$
 be  a continuous map in the flat topology
 satisfying $(A_0)-(A_4)$ (thus  Theorem \ref{discrete.sweepout} can be applied) and the following hypotheses:
\begin{itemize}
%\item[($C_0$)] $\Phi$ satisfies  {hypotheses $(A_0)$--$(A_4)$. Thus  Theorem \ref{discrete.sweepout} %can be applied};
\item[($A_5$)] $\max\{{\bf M}(\Phi(x)):x\in I^5_0\} =4\pi,$ and  $$x\in I_0^5\quad {\mbox{and}\quad{\bf M}(\Phi(x))=4\pi \Rightarrow \Phi(x)\in \mathcal{T}};$$
%\item[($C_2$)] $|\Phi|(\partial I^4 \times [1/2]) \subset \mathcal{T}$;
\item[($A_6$)] for every $\delta>0$ there exists $\varepsilon>0$ such that 
\begin{equation}\label{Jset.width}
x \in I_0^5 {\quad\mbox{and}\quad} {\bf F}(|\Phi(x)|,{\mathcal T})\leq \varepsilon \Rightarrow x\in J_\delta=\partial I^4\times \left[\frac12-\delta,\frac12+\delta\right];
\end{equation}
\item[($A_7$)]  {$|\Phi|(\partial I^4 \times [1/2]) \subset \mathcal{T}$ and}
$$|\Phi|_{*}([\partial I^4\times\{1/2\}])\neq 0\mbox{ in }H_3(\RP^3,\mathbb{Z}).$$ 
\end{itemize}

 We define $\hat \Phi:\partial I^4 \times I \rightarrow \mathcal{T}$ by 
$$
\hat \Phi(z,t) = |\Phi(z,1/2)|,
$$
for $(z,t) \in \partial I^4 \times I$. In particular, $\hat\Phi(x)=|\Phi(x)|$ for any $x \in \partial I^4 \times \{1/2\}$.

By applying Theorem \ref{discrete.sweepout} to $\Phi$, we obtain a  $(5,{\bf M})$-homotopy sequence of mappings into  $(\mathcal{Z}_2(S^3;{\bf M}),\Phi_{|I_0^n})$
$$C=\{\tilde \phi_i\}_{i\in \N}\quad {\mbox{such that}\quad{\bf L}(C)\leq \sup\{{\bf M}(\Phi(x)):x\in I^5\}.}$$ We denote  by $\Pi$ the corresponding $(5,{\bf M})$-homotopy class.  

\subsection{Theorem}\label{big.width.discrete}{\em  {Assume $\Phi$ satisfies hypotheses $(A_0)-(A_7)$.

Then
$${\bf L}(\Pi) > 4\pi.$$}
}

This theorem has the following important corollary.

\subsection{Corollary}\label{width.8pi.discrete} {\em  {Assume $\Phi$ satisfies hypotheses $(A_0)-(A_7)$.}
If $$\sup\{{\bf M}(\Phi(x)):x\in I^5\}<8\pi,$$
then there exists a smooth embedded minimal surface $\Sigma \subset S^3$  {with genus $g\geq 1$} such that
$$\mathrm{area}(\Sigma)={\bf L}(\Pi)>4\pi.$$
}
\begin{proof}
Using Theorem \ref{big.width.discrete}  {and $(A_5)$} we obtain that
$$ {4\pi=\sup\{{\bf M}(\Phi(x)):x\in I^5_0\}<{\bf L}(\Pi)}.$$ Hence we can apply Theorem \ref{pitts.min.max} to conclude the existence of  a stationary integral varifold $\Sigma$, whose support is a  smooth embedded minimal surface, such that
$$ {4\pi}<||\Sigma||(S^3)= {\bf L}(\Pi)\leq L(C)<8\pi.$$
Every minimal surface in $S^3$ has area bounded below by $4\pi$  {and so } the inequality above  implies that $\Sigma$ has multiplicity one. Since by Almgren \cite{almgren66} the great spheres are the only minimal surfaces in $S^3$ that are topological spheres, it follows that $\Sigma$ has  genus $g\geq 1$.  This implies the desired result.
\end{proof}

\subsection{Proof of Theorem \ref{big.width.discrete}}

We argue by contradiction. Assume that ${\bf L}(\Pi)=4\pi$, and consider  the critical sequence $S=\{\phi_i\}_{i\in\N}\in \Pi$ given by Proposition \ref{pulltight}. Suppose $\phi_i$ has domain $I(5,k_i)$, and ${\bf f}(\phi_i)=\delta_i$. Note that every varifold in ${\bf C}(S)$ is  stationary, since  any varifold in $|\Phi|(I^5_0)$ with  area $4\pi$ belongs to $\mathcal{T}$.

We will use cubical singular homology groups with integer coefficients  (see Massey \cite{massey}). 
If $X$ is a topological space, we denote by $C_n(X)$ the group
of cubical singular $n$-chains in $X$ with integer coefficients. If $f:X \rightarrow Y$ is a continuous map, we denote by $f_{\#}:C_n(X)\rightarrow C_n(Y)$ and $f_*:H_n(X,\mathbb{Z}) \rightarrow H_n(Y,\mathbb{Z})$ the homomorphisms induced by $f$. 

Note that we can  identify $\alpha \in I(5,k_i)_p$  with a $p$-singular cube $\alpha:I^p \rightarrow I^5$ in $I^5$ in a natural way (through an affine map).  If $R=\sum_{\alpha \in I(5,k_i)_p} n_\alpha \alpha \in C_p(I^5)$, $n_\alpha \in \mathbb{Z}$, we denote by $R_q$ the set of all $q$-cells of $I(5,k_i)$ that are faces  of some $\alpha$
with $n_\alpha \neq 0$. In this case we say that $R$ is subordinated to $I(5,k_i)$. The support of $R$ is the union of the supports of all $\alpha$
with $n_\alpha \neq 0$.

The proof is divided in three steps.

%\subsubsection{First step:}

\subsection{First step:} We construct a 4-chain $R(i) \in C_4(I^5)$, subordinated to $I(5,k_i)$, with $${\rm support}(\partial R(i)) \subset \partial I^5.$$ The chain $R(i)$ is constructed so that   $|\phi_i(x)|$ is sufficiently close to $\mathcal{T}$ for any $x \in R(i)_0$.

\medskip
  
Let $\varepsilon_0>0$ be small, to be chosen later. Then we choose $\delta>0$  such that
\begin{equation}\label{Jdelta.defi}
x \in J_\delta =\partial I^4 \times [1/2-\delta,1/2+\delta]  \Rightarrow {\bf F}(|\Phi(x)|, \hat{\Phi}(x)) \leq \varepsilon_0.
\end{equation}
It follows from condition $(A_6)$ that there exists $0<\varepsilon\leq \varepsilon_0/2$ such that
\begin{equation}\label{Jdelta}
x\in \partial I^5,\, {\bf F}(|\Phi(x)|,\mathcal{T})< 2\varepsilon \Rightarrow x \in J_\delta.
\end{equation}

 Consider $$\bar a(i)=\left\{\alpha\in I(5,k_i)_5: {\bf F}(|\phi_i(x)|,{\mathcal T})\geq\frac{\varepsilon}{2}\ \mbox{for  all}\, x\in \alpha_0 \right\}.$$ Let $a(i)$ be
the set of 5-cells $\alpha \in \bar a(i)$ for which we can find a sequence $\{\alpha_j\}_{j=1}^l \subset \bar a(i)$ with
$\alpha_1=\alpha$, $\alpha_l=\beta \otimes [0,3^{-k_i}]$ for some $\beta \in I(4,k_i)_4$, and such that $\alpha_j$ and $\alpha_{j+1}$ share a common 4-face for each $j=1,\dots,l-1$. Because $\phi_i$ vanishes on $(I(4,k_i)\otimes \langle[0]\rangle)_0$,  {if $\varepsilon_0$ is sufficiently small} we have that $\beta \otimes [0,3^{-k_i}] \in a(i)$ for every $\beta \in I(4,k_i)_4$.  {Loosely speaking, $\cup_{\alpha\in a(i)}\alpha$ is the connected component of $\cup_{\alpha\in \bar a(i)} \alpha$ that contains $I^4\times\{0\}$.}

Let $b(i)$ denote the set of 4-cells in $I(5,k_i)$ that are faces of exactly one 5-cell
in $a(i)$. Consider the following 5-chain
$$A(i)=\sum_{\alpha \in a(i)}\alpha \in C_5(I^5).$$
We have
$$
\partial A(i)=\sum_{\alpha\in b(i)} \mathrm{sgn}(\alpha)\alpha,
$$
where $\mathrm{sgn}(\alpha)$ is equal to 1 or $-1$. Note that $\beta \otimes [0] \in b(i)$ for every $\beta \in I(4,k_i)_4$.
From the definition of the boundary homomorphism, we have that $\mathrm{sgn}(\beta \otimes [0])=-1$ for 
every $\beta \in I(4,k_i)_4$.

Let $c(i)$ be the set of 4-cells of $b(i)$ that belong to 
the subcomplex $T(5,k_i) \cup S(5,k_i)$. Then we have  the disjoint decomposition  below:
\begin{equation}\label{width.c(i)}
b(i) \cap I_0(5,k_i)_4=c(i)\cup \{\alpha=\beta \otimes [0]:\beta \in I(4,k_i)_4\}.
\end{equation}

We define  the 4-chain:
\begin{multline}\label{width.R(i)2}
R(i)=\partial A(i)-\sum_{\alpha\in b(i)\cap I_0(5,k_i)_4} \mathrm{sgn}(\alpha)\alpha\\
=\partial A(i)+\sum_{\beta \in I(4,k_i)_4} \beta \otimes [0]-\sum_{\alpha\in c(i)} \mathrm{sgn}(\alpha)\alpha.
\end{multline}
Note that ${\rm support}(\partial R(i)) \subset \partial I^5$.

\subsection{Lemma}\label{proximity.rp3.width}\textit{We have
$$\sup\{{\bf F}(|\phi_i(x)|,\mathcal T):x\in R(i)_0\}\leq \varepsilon$$
for every  sufficiently large $i$ such that $R(i)\neq 0$.
}

\begin{proof}
Let $i$ be sufficiently large such that $5\delta_i \leq \varepsilon/2$, and let $x\in R(i)_0$. From the definition of $R(i)$ we see that we can find a $4$-cell $$\alpha\in b(i) \cap ( I(5,k_i)_4\setminus I_0(5,k_i)_4) \quad\mbox{ with } x\in \alpha_0.$$ Thus $\alpha$ is the common 4-face of two distinct cells  $\beta,\gamma \in I(5,k_i)_5$. Since $\alpha \in b(i)$, we can suppose, after a possible relabeling, that $\beta \in a(i)$ and $\gamma \notin a(i)$. It follows from the definition of $a(i)$ that  $\gamma \notin \bar a(i)$. This means that there exists 
$y\in \gamma_0$ with ${\bf F}(|\phi_i(y)|,{\mathcal T})< \varepsilon/2$. Note that ${\bf d}(x,y)\leq 5$, hence 
$$ {\bf F}(|\phi_i(x)|,{\mathcal T})\leq{\bf F}(|\phi_i(y)|,{\mathcal T})+{\bf F}(|\phi_i(y)|, |\phi_i(x)|)< \frac{\varepsilon}{2}+5\delta_i<\varepsilon.$$
\end{proof}

\subsection{Second step:} We prove that the  support of $R(i)$ separates $I^4 \times \{0\}$ from $I^4 \times \{1\}$. This uses the assumption that  ${\bf L}(\Pi)=4\pi$ in a fundamental way. Then we prove that $\partial R(i)$ is homologous to $\partial I^4\times\{1/2\}$
in $H_3(\partial I^4 \times I, \Z)$.

\subsection{Lemma}\label{whites.disconnected} \textit{If $i$ is sufficiently large, then no 5-cell of the type $\beta \otimes [1-3^{-k_i},1]$, $\beta \in I(4,k_i)_4$, belongs to $a(i)$. 
}
\begin{proof}
Suppose, by contradiction, that there exists $\alpha=\beta \otimes [1-3^{-k_i},1]$, $\beta \in I(4,k_i)_4$, with  $\alpha \in a(i)$. 
Then we can find a sequence of maps
$$\gamma_i:I(1,n_i)_0\rightarrow I(5,k_i)_0,$$
with
\begin{itemize}
\item $n_i\geq k_i$ and  ${\bf d}(\gamma_i(x),\gamma_i(y))\leq 1$ if ${\bf d}(x,y)\leq 1$;
\item %  $${\bf d}(x,y)\leq 1\implies {\bf d}(\gamma_i(x),\gamma_i(y))\leq 1,$$
 $\gamma_i([0])\in (I(4,k_i)\otimes \langle[0]\rangle)_0$ and $\gamma_i([1])\in (I(4,k_i)\otimes \langle[1]\rangle)_0;$
 \item$\gamma_i(I(1,n_i)_0)\subset \cup_{\alpha \in a(i)}\alpha.$
 \end{itemize}
In particular, putting $\sigma_i=\phi_i\circ \gamma_i$, we have
\begin{equation}\label{width.inside.ai}
{\bf F}(|\sigma_i(x)|,{\mathcal T})\geq \frac{\varepsilon}{2}
\end{equation}
for all $x \in I(1,n_i)_0$.

We now show that $\gamma_i$ is {\em homotopic to a vertical path}, meaning we can find  a map
$$\psi_i:I(1,s_i)_0\times I(1,s_i)_0 \rightarrow I(5,k_i)_0$$
such that 
\begin{itemize}
\item [(a)] $\psi_i([0], \cdot)=\gamma_i\circ {\bf n}(s_i,n_i)$ and  $\psi_i([1], y)= c+ {\bf n}(s_i,k_i)(y) { e_5}$, where
$$c=\frac{1}{3}(1,1,1,1,0)\quad\mbox{and}\quad e_5=(0,0,0,0,1);$$
\item [(b)] $\psi_i(\cdot, [0])\in (I(4,k_i)\otimes\langle [0]\rangle)_0\quad\mbox{and}\quad\psi_i(\cdot, [1])\in (I(4,k_i)\otimes\langle[1]\rangle)_0$;
\item [(c)] if $x,y\in I(2,s_i)_0$,  $${\bf d}(x,y)\leq 1\Rightarrow {\bf d}(\psi_i(x),\psi_i(y))\leq 5.$$
\end{itemize}

In order to show this we associate to each $\gamma_i$ a  piecewise linear continuous curve  $\tilde \gamma_i:I\rightarrow I^5$ given by
\begin{equation}\label{continuous.curve}
\tilde \gamma_i(t) = (j+1-3^{n_i}t)\gamma_i\left(\frac{j}{3^{n_i}}\right)+(3^{n_i}t-j)\gamma_i\left(\frac{j+1}{3^{n_i}}\right)
\end{equation}
 for every $j3^{-n_i}\leq t \leq(j+1)3^{-n_i}$, $j=0,\ldots,3^{n_i}-1$.  Note that $\tilde\gamma_i(t) = \gamma_i([t])$
 if $[t] \in I(1,n_i)_0$.
 
 %We abuse notation and denote this curve by $\gamma_i$ as well. 
 
 %Because $\pi_1(I^5,I^4\times\{0\}\cup I^4\times\{1\})$ is trivial, we can find a continuous map
 Let $\psi:I^2\rightarrow I^5$ be given by  $\psi(u,t)=(1-u) \tilde\gamma_i(t) + u(c+te_5)$. Then
 \begin{equation}\label{homotopy.cond1}
 \psi(0,t)=\tilde \gamma_i(t),\quad  \psi(1,t)=c+te_5,\quad\mbox{ for all }t\in I,
 \end{equation}
 and 
 \begin{equation}\label{homotopy.cond2}
 \psi(I\times\{0\})\subset I^4\times\{0\}, \quad \psi(I\times\{1\})\subset I^4\times\{1\}.
 \end{equation}
Choose $s_i\geq n_i$ sufficiently large so that
\begin{equation}\label{closest.disconnected}
|\psi(x)-\psi(y)|\leq \frac{1}{3^{k_i+2}}\quad\mbox{for all }x,y\in I(2,s_i)_0\mbox{ with }{\bf d}(x,y)\leq 1.
\end{equation}
For $x \in I(2,s_i)_0$, we choose $\psi_i(x) \in I(5,k_i)_0$ to satisfy $$d(\psi_i(x),\psi(x))=d(\psi(x),I(5,k_i)_0).$$ Note that
such choice might not be unique.  If $\psi(x)\in I^4\times\{j\}$, $j=0$ or 1, then it follows from the definition that $\psi_i(x) \in (I(4,k_i) \otimes \langle[j]\rangle)_0$. This proves property (b) for $\psi_i$.  From \eqref{continuous.curve} and \eqref{homotopy.cond1} we obtain property (a) for $\psi_i$.  Finally, from \eqref{closest.disconnected} we have that $\psi_i(x)$ and $\psi_i(y)$ are vertices of a common 5-cell in $I(5,k_i)$ if $x,y\in I(2,s_i)_0$ satisfy ${\bf d}(x,y) \leq 1$. This establishes property (c). 

Consider the sequence $D=\{\sigma_i\}_{i\in \N}$, where  $\sigma_i=\phi_i\circ\gamma_i$.  From the fact that $\gamma_i$ is homotopic to a vertical path, we obtain that  $D$ is a $(1,{\bf M})$-homotopy sequence of mappings into $(\mathcal{Z}_2(S^3;{\bf M}),\{0\})$ that is homotopic with $\{v_i\}_{i\in \N}$, where $$v_i:I(1,k_i)_0\rightarrow  \mathcal{Z}_2 {(S^3)},\quad v_i(x)=\phi_i(c+xe_{5}).$$
Hence $D$ and $\{v_i\}_{i\in \N}$ belong to  the same element $\Omega$ in $\pi_1^{\#}(\mathcal{Z}_2(S^3;{\bf M}),\{0\}).$ 

Since $S$ is homotopic with $C$ ($S,C\in \Pi$), we obtain from Theorem \ref{discrete.sweepout}  {(iii)} that $\Omega$ is non-trivial in  $\pi_1^{\#}(\mathcal{Z}_2(S^3;{\bf M}),\{0\})$. Hence it follows from Pitts (\cite{pitts}, Theorem 4.6, Corollary 4.7)   that ${\bf L}(\Omega)>0$. From Theorem \ref{pitts.min.max} (applied to $\Omega \in \pi_1^{\#}(\mathcal{Z}_2(S^3;{\bf M}),\{0\})$), we get the existence of a stationary integral varifold $\Sigma$ whose support is a smooth embedded minimal surface  in $S^3$ and such that
\begin{equation}\label{inequality.width}
4\pi\leq||\Sigma||(S^3)={\bf L}(\Omega)\leq {\bf L}(D)\leq {\bf L}(S)={\bf L}(\Pi)=4\pi.
\end{equation}
The first inequality follows because the area of any  minimal surface in $S^3$ is at least  $4\pi$. The second inequality follows  because $D \in \Omega$, and the third inequality follows 
because the definition of $D$ implies ${\bf K}(D) \subset {\bf K}(S)$. We note that  this string of inequalities implies that $\Sigma$ must be a great sphere.

From \eqref{inequality.width} we also get that $D$ is a critical sequence (since ${\bf L}(\Omega)={\bf L}(D)$), and that  ${\bf C}(D)\subset {\bf C}(S)$ (since ${\bf L}(D)={\bf L}(S)$). In particular, every element of ${\bf C}(D)$ is a stationary varifold  {because every varifold in {\bf C}(S) is stationary}.  We know from  Theorem \ref{pitts.min.max}  that  the surface $\Sigma$ in \eqref{inequality.width} can be chosen to belong to ${\bf C}(D)$, hence
 $${\bf F}({\bf C}(D),{\mathcal T})=0.$$
 
 On the other hand, according to \eqref{width.inside.ai}, we have ${\bf F}({\bf C}(D),{\mathcal T})\geq \varepsilon/2$. This gives us a contradiction.

\end{proof}

\subsection{Lemma}\label{boundary.R(i).width} \textit{ For sufficiently large $i$, ${\rm support}(\partial R(i)) \subset J_\delta$ and
$$[\partial R(i)]=[\partial I^4\times\{1/2\}]\quad\mbox{ in }H_3(J_\delta,\Z).$$
In particular, $R(i)\neq 0$.
}
\begin{proof}
 We obtain from Lemma \ref{whites.disconnected}  that no 4-cell in $b(i)$ 
belongs to  the subcomplex $I(4,k_i)\otimes \langle[1]\rangle$. Therefore $c(i)\subset S(5,k_i)_4$. If
$$C(i)=\sum_{\alpha\in c(i)}\mathrm{sgn}(\alpha)\alpha,$$
we get that   $C(i)$ is  a 4-chain in $\partial I^4 \times I$. Since, from \eqref{width.R(i)2},
$$\partial R(i) =\partial \left(\sum_{\beta \in I(4,k_i)_4} \beta \otimes [0]\right)-\partial C(i),$$
we conclude that $\partial R(i)$ is a 3-cycle in $\partial I^4 \times I$ and
$$[\partial R(i)]=[\partial I^4\times\{0\}]=[\partial I^4\times\{1/2\}]\mbox{ in }H_3(\partial I^4\times I,\Z).$$ 
Since ${\rm support}(\partial R(i))\subset \partial I^4 \times I$, we know from Lemma \ref{homotopy.sequence.boundary} that 
$$\lim_{i\to\infty}\sup\{{\bf F}(\phi_i(x),\Phi(x)):x\in \partial R(i)_0\}=0.$$
Combining this with 
Lemma \ref{proximity.rp3.width} and \eqref{Jdelta}, we obtain  that ${\rm support}(\partial R(i))\subset J_\delta$ if $i$ is sufficiently large. Now we use  a deformation retraction of $\partial I^4 \times I$ onto $J_\delta$ to get 
$$[\partial R(i)]=[\partial I^4\times\{1/2\}]\quad\mbox{ in }H_3(J_\delta,\Z).$$
\end{proof}

\subsection{Third step:} We construct a continuous map $f_i:{\rm support}(R(i))\rightarrow \mathcal{T}$
that extends $\hat\Phi_{|{\rm support}(\partial R(i))}$. From that we derive a contradiction, using that
$|\Phi|_*([\partial I^4\times\{1/2\}])\neq 0 $ in $H_3(\RP^3,\Z)$.

\subsection{Lemma.}\label{continuous.extension.width} \textit{For all sufficiently large $i$, there exists  a continuous function
$$f_i:{\rm support}(R(i))\rightarrow \mathcal{T}$$
such that ${f_i}_{|{\rm support}(\partial R(i))}=\hat\Phi_{|{\rm support}(\partial R(i))}$.
}

\begin{proof}
Throughout the proof of this lemma, $D_r(p)$ denotes  a ball centered at $p$ of radius $r$ in $\RP^3$ with respect to the standard metric. Unless otherwise stated,  geometric quantities in $\RP^3$ such as convexity, diameter, or distances, are  computed with respect to the standard metric.

Let $\eta>0$ be chosen so that every ball of radius $11\eta$ in $\RP^3$ is geodesically convex.  The topology  induced by the ${\bf F}$-metric on ${\mathcal T}\approx \RP^3$ coincides with the topology induced by  the geodesic distance of $\RP^3$.  Therefore, by compactness, we can find $c_0>0$ so that
%\begin{equation}\label{distance.geod.width}
%{\bf B}_r^{\bf F}(p)\cap \RP^3\subset {B}_{c_0r}(p)\quad\mbox{for all }p\in \RP^3.
%\end{equation}
\begin{equation}\label{distance.geod.width}
p,q \in \mathcal{T} {\rm \ satisfy \ }{\bf F}(p,q) < \frac{\eta}{2c_0} \Rightarrow {\rm dist}(p,q) < \frac{\eta}{2}.
\end{equation}
At this point we can choose $\varepsilon_0=\frac{\eta}{20c_0}.$

%We also choose $\delta>0$ such that 
%\begin{equation}\label{J.proximity}
 %{x\in J_{\delta}}\implies{\bf F}(|\Phi(x)|, \hat \Phi(x)) \leq \varepsilon_0.
%\end{equation}

Let $i$ be sufficiently large such that Lemmas \ref{proximity.rp3.width}, \ref{whites.disconnected}, and \ref{boundary.R(i).width} apply, and we have:
\begin{itemize}
\item[(a)] ${\bf f}(\phi_i)\leq \varepsilon_0$;
\item[(b)] for every $x\in S(5,k_i)_0$ we have ${\bf F}(\phi_i(x),\Phi(x))< \varepsilon_0$ (using Lemma \ref{homotopy.sequence.boundary});
\item[(c)] for every $\alpha\in I_0(5,k_i)_4$,
$$\sup\{{\bf F}(\hat\Phi(x),\hat\Phi(y)): x,y\in \alpha\cap J_\delta\}< \varepsilon_0.$$
This combined with  \eqref{distance.geod.width} gives
$$\sup\{\mathrm{dist}(\hat\Phi(x),\hat\Phi(y)): x,y\in \alpha\cap J_\delta\}< \frac{\eta}{2}.$$
\end{itemize}

Define $f^0_i:R(i)_0\rightarrow \mathcal{T}$ as follows:   if $x\in \partial R(i)_0$ we make $f^0_i(x)=\hat\Phi(x)$, otherwise we choose
$f^0_i(x)\in \mathcal{T}$ such that 
$${\bf F}(f^0_i(x),|\phi_i(x)|)={\bf F}(|\phi_i(x)|, {\mathcal T}).$$

We now prove that 
\begin{equation}\label{width.contained.rp3}
\mathrm{diam}(f_i^0(\alpha_0))< \frac{\eta}{2}
\end{equation}
for every $\alpha\in R(i)_4$.
From  \eqref{distance.geod.width}, it suffices to show that
\begin{equation*}
{\bf F}( f_i^0(x), f_i^0(y)) < \frac{\eta}{2c_0}
\end{equation*}
for every $\alpha\in R(i)_4$ and $x,y\in\alpha_0$.
To that end, consider $\alpha\in R(i)_4$ and $x,y\in\alpha_0$. In particular we have ${\bf d}(x,y) \leq 4$. If both $x,y\in \partial R(i)_0,$ then the inequality above follows from property (c). If only one of the vertices, say $x$, belongs to $\partial R(i)_0$ then, using the definition of $\hat \Phi$,  (\ref{Jdelta.defi}), Lemma \ref{proximity.rp3.width}, Lemma \ref{boundary.R(i).width}, properties (a) and (b), we have that
 \begin{eqnarray*}
{\bf F}(f_i^0(x),f_i^0(y))&\leq& {\bf F}(\hat\Phi(x),|\Phi(x)|)+{\bf F}(|\Phi(x)|, f_i^0(y))\\
&\leq& \varepsilon_0+{\bf F}(|\Phi(x)|, f_i^0(y))\\
&\leq&\varepsilon_0+{\bf F}(|\Phi(x)|, |\phi_i(x)|)+{\bf F}(|\phi_i(x)|,f_i^0(y))\\
&\leq&2\varepsilon_0+{\bf F}(|\phi_i(x)|,|\phi_i(y)|)+{\bf F}(|\phi_i(y)|,f_i^0(y))\\
&\leq&  6\varepsilon_0+{\bf F}(|\phi_i(y)|,\mathcal{T})\\
&\leq&7\varepsilon_0<\frac{\eta}{2c_0}.
\end{eqnarray*}
Finally, if $x,y\notin \partial R(i)_0$ we have from Lemma \ref{proximity.rp3.width} and property (a) that
\begin{eqnarray*}
  &&{\bf F}(f_i^0(x),f_i^0(y))\\
&&  \hspace{1cm} \leq  {\bf F}(f_i^0(x),|\phi_i(x)|)+{\bf F}(|\phi_i(x)|,|\phi_i(y)|)+{\bf F}(|\phi_i(y)|,f_i^0(y))\\
  &&\hspace{1cm}={\bf F}(|\phi_i(x)|,\mathcal{T})+{\bf F}(|\phi_i(x)|,|\phi_i(y)|)+{\bf F}(|\phi_i(y)|,\mathcal{T})\\
  && \hspace{1cm}\leq 6\varepsilon_0<\frac{\eta}{2c_0}.
\end{eqnarray*}

We now proceed to the iterative construction of $f_i$.
We cover $\RP^3$ with a finite union of  balls $\{D_{\eta/2}(p_k)\}_{k=1}^N$, where each $D_{11\eta}(p_k)$ is geodesically convex. We denote by $R(i)^{(j)}$ ($\partial R(i)^{(j)}$) the union of the supports of all $q$-cells $\alpha \in R(i)_q$ ($\alpha \in \partial R(i)_q$) with $q\leq j$.
 The map $$f^j_i:R(i)^{(j)} \rightarrow \mathcal{T}$$ is called a {\em continuous $j$-extension of} $f_i^0$ if
\begin{itemize}
\item[(1)] $f^j_i=f^0_i$ on $R(i)_0$, and $f^j_i=\hat\Phi$ on  $\partial R(i)^{(j)}$;
\item[(2)] for every  $\alpha \in R(i)_j$, with $j\geq 1$, we have $$\mathrm{diam}(f^j_i(\alpha))\leq (2^j-2+2^{j-2})\eta.$$
\end{itemize}

Assuming the existence of    a continuous $j$-extension $f^j_i$ of $f_i^0$, $j\leq 3$, we will construct a  continuous $(j+1)$-extension
 $f^{j+1}_i$ of $f_i^0$. Let $\alpha \in R(i)_{j+1}$. If $\alpha \in \partial R(i)_{j+1}$, we set $f^{j+1}_{i}=\hat \Phi$ on $\alpha$. 
 In this case  it follows from property (c) that property (2) holds for $\alpha$.  {We note that, }since $f_i^j=\hat \Phi$ on ${\rm support}(\partial \alpha)$, we  {have}  
 $f^{j+1}_i=f^j_i$ on ${\rm support}(\partial \alpha)$. If  $\alpha \notin \partial R(i)_{j+1}$, we know from \eqref{width.contained.rp3} and property (1)  that
$$f^j_i(\alpha_0)\subset B_{\eta}(p_k)\quad\mbox{for some }k=1,\ldots,N.$$
By applying property (2) to the  $j$-faces of $\alpha$, we obtain from the inclusion above that
$$f^j_i({\rm support}(\partial \alpha))\subset  B_{(2^j-1+2^{j-2})\eta}(p_k).$$
We can now use the convexity of  $B_{11\eta}(p_k)$  to construct  a continuous map
$$f^{j+1}_i:{\rm support}(\alpha)\rightarrow B_{(2^j-1+2^{j-2})\eta}(p_k) $$
such that $f^{j+1}_i=f^j_i$ on ${\rm support}(\partial \alpha)$.
Furthermore, we have
$$\mathrm{diam}(f^{j+1}_i(\alpha))\leq 2(2^j-1+2^{j-2})\eta=(2^{j+1}-2+2^{j-1})\eta.$$
 It follows that $f^{j+1}_i$ is a continuous and well-defined $(j+1)$-extension of $f_i^0$.

Arguing inductively, we construct a 4-extension $f_i^4$ of $f_i^0$. The map $f_i=f^4_i:{\rm support}(R(i)) \rightarrow \mathcal{T}$ is continuous and satisfies $f_i=\hat \Phi$ on ${\rm support}(\partial R(i))$.
\end{proof}

We now finish the argument. The map $f_i:{\rm support}(R(i))\rightarrow \mathcal{T} \approx \RP^3$ constructed in Lemma \ref{continuous.extension.width} induces a homomorphism in homology
$${f_i}_{*}:H_{*} ({\rm support}(R(i)),\Z)\rightarrow H_{*}(\RP^3,\Z).$$
Since  $f_i=\hat \Phi$ on ${\rm support}(\partial R(i))$,  we have
$$ \hat \Phi_{*}[\partial R(i)]={f_i}_{*}[\partial R(i)]=[{f_i}_{\#}\partial (R(i))]=[\partial {f_i}_{\#}(R(i))]=0.$$
But Lemma \ref{boundary.R(i).width} implies that $$\hat \Phi_{*}[\partial R(i)]=\hat \Phi_{*}[\partial I^4\times\{1/2\}]=|\Phi|_{*}([\partial I^4\times\{1/2\}])\in H_3(\RP^3,\Z).$$ This is a  contradiction since we have assumed from the beginning that $|\Phi|_*([\partial I^4\times\{1/2\}])\neq 0 $ in $H_3(\RP^3,\Z)$.

%%%%%%%%%%%%%%%%%%%%%%%%%%%%%%%%%%%%%%%%%%%%%%%%%%%%%%%%
%%%%%%%%%%%%%%%%%%%%%%%%%%%%%%%%%%%%%%%%%%%%%%%%%%%%%%%%%%%

%%%%%%%%%%%%%%%%%%%%%%%%%%%%%%%%%%%%%%%%%%%%%%%%%%%%%
%%%%%%%%%%%%%%%%%%%%%%%%%%%%%%%%%%%%%%%%%%%%%%%%%%%%%%%%%%%

\section{Proof of Theorem B}\label{thmb.section}
Let 
\begin{multline*}
\mathcal{F}_1=\{ S \subset S^3: S \mbox{ is an embedded closed minimal surface of}  \\
\mbox{genus } g(S)\geq 1\}.
\end{multline*}
The Jacobi operator of $\Sigma$ is given by  $L=\Delta + |A|^2 + 2$, where
$A$ denotes the second fundamental form of $\Sigma.$  The index of $\Sigma$, denoted by ${\rm index}(\Sigma)$, is defined to be the number of negative eigenvalues of $L$.

Theorem B follows from the next theorem.
\subsection{Theorem}\textit{We have 
$$2\pi^2=\inf_{S \in \mathcal{F}_1} {\rm area}(S)$$
 and, for every $\Sigma\in \mathcal{F}_1$, 
${\rm area}(\Sigma)=2\pi^2$ if and only if $\Sigma$ is the Clifford torus up to isometries of $S^3$.}
\begin{proof}
From  Theorem \ref{existence.minimizer}, choose $\Sigma\in \mathcal{F}_1$ such that
$${\rm area}(\Sigma)= \inf_{S \in \mathcal{F}_1} {\rm area}(S)\leq 2\pi^2.$$
Consider the min-max family $\Phi$ (see Definition \ref{Fi.family})  and the homotopy class $\Pi$ (see  Definition \ref{homotopy.class.sigma}) associated with $\Sigma$.    {Theorem \ref{modified.family} (vi), (vii), and (viii) implies that hypotheses $(A_5)$, $(A_6)$, and $(A_7)$ are satisfied}.  {Thus we can apply} Corollary \ref{width.8pi.discrete}  {and conclude} the existence of $S\in  \mathcal{F}_1$ so that, from Theorem \ref{modified.family}(iii), we have
$$ {\rm area}(S)={\bf L}(\Pi)\leq \sup\{{\bf M}(\Phi(x)):x\in I^5\}\leq\mathcal{W}(\Sigma)={\rm area}(\Sigma).$$
Thus ${\bf L}(\Pi)={\rm area}(\Sigma)$.

We want to show that ${\rm index}(\Sigma)\leq 5$ because, by a theorem of Urbano \cite{urbano},  that implies $\Sigma$ must be  the Clifford torus up to  isometries of $S^3$. Before we do so, we need to establish a non-degeneracy lemma for  the Jacobi operator on $\Sigma$.

Let $\{e_1,e_2,e_3,e_4\}$ be the standard orthonormal basis of $\mathbb{R}^4$. For $x \in \Sigma$, define $\psi_i(x)=\langle N(x),e_i\rangle$  for each $1\leq i\leq 4$, and $\psi_5(x)=1$. Denote by $E$ the subspace of $C^\infty(\Sigma)$  spanned by $\{\psi_j\}_{1\leq j\leq 5}$. 

Notice that $L\psi_i=2\psi_i$ for $1\leq i\leq 4$ (see \cite{urbano}).

Recall the definition of $F_v$ in Section \ref{associated.t} and $N_v$ in Remark \ref{canonical.remark} (1). Choose $\delta>0$ such that the map
$$
P: B^4_\delta(0) \times (-\delta,\delta) \times \Sigma \rightarrow S^3, \quad P_{(v,t)}(x)=(\cos t)F_v(x)+(\sin t)N_{v}(x) 
$$
has $\Sigma_{(v,t)}=P_{(v,t)}(\Sigma)$, where $\{\Sigma_{(v,t)}\}$ is the canonical family defined in  Definition \ref{associated.family}, and such that $P_{(v,t)}$ is an embedding of $\Sigma$ into $S^3$.

If $1\leq i\leq 4$, $x\in\Sigma$, we have
$$\left\langle\frac{d}{ds}_{|s=0}P_{(se_i,0)}(x), N(x)\right\rangle= -2\langle e_i, N(x) \rangle= -2\psi_i(x)$$
and so
\begin{equation}\label{phi.thmb}
\frac{d^2}{(ds)^2}_{|s=0}{\rm area}\left(P_{(se_i,0)}(\Sigma)\right)=-4 \int_{\Sigma}\psi_iL\psi_i\,d\Sigma.
\end{equation}
Similarly,
\begin{equation}\label{phi5.thmb}
\frac{d^2}{(ds)^2}_{|s=0}{\rm area}\left(P_{(0,s)}(\Sigma)\right)=-\int_{\Sigma}\psi_5L\psi_5\,d\Sigma.
\end{equation}

\subsection{Lemma}\label{index.five}\textit{
$$
-\int_{\Sigma} \psi L\psi \, d\Sigma < 0 \quad\mbox{for every }\psi \in E \setminus \{0\}.
$$
}
\begin{proof}
Let
$$f(v,t)={\rm area}(\Sigma_{(v,t)})={\rm area}(P_{(v,t)}(\Sigma)), \quad (v,t)\in B^4_\delta(0) \times (-\delta,\delta).$$
 Since $\Sigma$ is minimal, we have  $f(0,0)=\mathcal{W}(\Sigma)$ and $Df(0,0) = 0$.
We also know, from Theorem \ref{heintze.karcher}, that $f(v,t) \leq f(0,0)$ for every $(v,t) \in B^4_\delta(0) \times (-\delta,\delta)$. Hence $D^2f(0,0) \leq 0$
and this means that
$$
-\int_{\Sigma} \psi L\psi \, d\Sigma \leq 0\quad\mbox{for every }\psi \in E.
$$

Suppose the lemma  is not true. We can find $\phi \in E\setminus \{0\}$ such that
$$
-\int_{\Sigma} \phi L\psi \,d\Sigma= -\int_{\Sigma} \psi L\phi \,d\Sigma=0 \quad\mbox{for every }\psi \in E.
$$
Hence 
\begin{equation}\label{lemma.five.orthogonal}
\int_{\Sigma} \phi \, \psi_i \,d\Sigma =0\quad \mbox{ for every }1\leq i\leq 4, \quad\mbox{and }\int_{\Sigma}L\phi\, d\Sigma=0.
\end{equation}
This implies, since $\psi_5=1\in E$, the existence of  $c\in \mathbb{R}$ such that  $$1= c \phi + \psi,\quad\mbox{where}\quad\psi = \sum_{i=1}^4 a_i\psi_i.$$
Hence, because $\psi$ is an eigenfunction of $L$, we have
$$
\int_{\Sigma} (|A|^2 +2)\, d\Sigma = \int_{\Sigma} L(1) \,d\Sigma
=\int_{\Sigma} (cL\phi + L\psi) \, d\Sigma
=2 \int_{\Sigma} \psi \, d\Sigma.
$$

On the other hand, we also have $1=c^2 \phi^2 + 2c\phi \psi+\psi^2$.
If we integrate over $\Sigma$, we obtain from \eqref{lemma.five.orthogonal}
\begin{multline*}
{\rm area}(\Sigma) = \int_{\Sigma} (c^2 \phi^2 + 2c\phi \psi+\psi^2)\,d\Sigma
\geq \int_{\Sigma}\psi^2\, d\Sigma
=  \int_{\Sigma}\psi (1-c\phi)\, d\Sigma\\
=\int_{\Sigma} \psi\, d\Sigma.
\end{multline*}
%Putting these two inequalities together we have

Hence
$$
2\, {\rm area}(\Sigma) \leq \int_{\Sigma} (|A|^2 +2)\, d\Sigma
= 2\int_{\Sigma} \psi \, d\Sigma
\leq 2\, {\rm area}(\Sigma).
$$
This implies $A=0$ and so $\Sigma$ is a great sphere. This contradicts our assumption that $\Sigma\in \mathcal{F}_1$.
\end{proof}

 Suppose, by contradiction, that ${\rm index}(\Sigma) \geq 6$. {The idea is to construct a comparison map
 $$C':\overline{B}^4 \times [-\pi,\pi]\rightarrow \mathcal{Z}_2(S^3)$$
 which coincides with $C$, the map given by Theorem \ref{canonical.family.continuous}, outside a neighborhood of the origin. Using this map we will conclude that
 $${\rm area}(\Sigma_w)= {\rm area}(\Sigma)\quad\mbox{for some}\quad w\in B^4\setminus\{0\}.$$
 Finally, we show that this identity implies $\Sigma$ is totally geodesic, which gives us the desired contradiction.
 }

Because ${\rm index}(\Sigma) \geq 6$, there exists $\varphi \in C^\infty(\Sigma)$ such that 
\begin{itemize}
\item $-\int \varphi L\varphi \,d\Sigma < 0$,
\item $-\int \varphi L\psi_i d\Sigma=0$ for $1\leq i\leq 5$.
\end{itemize} 
Let $X$ be any vector field such that $X=\varphi N$ along $\Sigma$, and let $\{\Gamma_s\}_{s\geq 0}$ be the one parameter group of diffeomorphisms generated by $X$.

Define $f: B^4_\delta(0) \times (-\delta,\delta)\times (-\delta,\delta) \rightarrow \mathbb{R}$ by
$$
f(v,t,s) = {\rm area}(\Gamma_s\circ P_{(v,t)}(\Sigma)).
$$
We have $f(0,0,0)={\rm area}(\Sigma)$, and $Df(0,0,0)=0$ since $\Sigma$ is minimal.
It follows from the choice of $\varphi$, \eqref{phi.thmb}, \eqref{phi5.thmb}, and Lemma \ref{index.five}, that $D^2f(0,0,0) <0$.
This means that there exists $0<\delta_1\leq \delta$ such that
\begin{equation}\label{thmb.area.menor}
{\rm area}(\Gamma_s\circ P_{(v,t)}(\Sigma)) < f(0,0,0)={\rm area}(\Sigma)
\end{equation}
for every $(v,t,s) \in  (B^4_{\delta_1}(0) \times (-\delta_1,\delta_1)\times (-\delta_1,\delta_1))\setminus \{(0,0,0)\}$.

Let $\beta:\mathbb{R}^5 \rightarrow \mathbb{R}$ be a smooth function such that $0\leq \beta(y) \leq \delta_1/2$ for $y\in \mathbb{R}^5$, $\beta(y)=0$ if $|y| \geq \delta_1/2$ and $\beta(y) =\delta_1/2$ if $|y| \leq \delta_1/4$.
We then define
$$
C'(v,t)= [|\Gamma_{\beta(v,t)}\circ P_{(v,t)}(\Sigma)|]\in \mathcal{Z}_2(S^3)\quad\mbox{for }|(v,t)|<\delta_1.
$$
We have that $C'(v,t)=C(v,t)$ if $\delta_1/2<|(v,t)|<\delta_1$, where $C$ is the map given by Theorem \ref{canonical.family.continuous},  and this means we can extend $C'$ to a continuous  map in the flat topology $$C':\overline{B}^4 \times [-\pi,\pi]\rightarrow \mathcal{Z}_2(S^3)$$ by defining $C'(v,t)=C(v,t)$ if $|(v,t)|\geq \delta_1$.  Note that from \eqref{thmb.area.menor} we have
\begin{equation}\label{area.menor.C}
\sup\{{\bf M}(C'(v,t)):|(v,t)|\leq \delta_1\}<{\rm area}(\Sigma)
\end{equation}
and so
$$\sup\{{\bf M}(C'(v,t)):(v,t)\in \overline B^4\times[-\pi,\pi]\}\leq {\rm area}(\Sigma).$$
{We use the map $C'$ to show:}
\subsection{Lemma}\label{conforme.igual}{\em {There is $w\in B^4\setminus\{0\}$ so that ${\rm area}(\Sigma_w)={\rm area}(\Sigma)$.}}
\begin{proof}

If we replace $C$  {by} $C'$ in Definition \ref{Fi.family}, we get  a continuous map in the flat topology $\Phi':I^5 \rightarrow \mathcal{Z}_2(S^3)$ that,  {according to Theorem \ref{modified.family}}, satisfies  hypotheses $(A_0)$--$(A_7)$ and thus  Theorem \ref{big.width.discrete} can be applied.

Consider the $(5,{\bf M})$-homotopy sequence $S=\{\phi_i\}_{i\in\N}$ of mappings into $(\mathcal{Z}_2(S^3;{\bf M}), {\Phi'_{|I^5_0}})$ given by  Theorem \ref{discrete.sweepout}, and denote by $\Pi'$ the corresponding $(5,{\bf M})$-homotopy class. From Corollary  \ref{width.8pi.discrete} we get the existence of  a smooth embedded minimal surface $\Sigma'$  {with genus $g\geq 1$} such that
$$4\pi < \mathrm{area}(\Sigma')={\bf L}(\Pi').$$
Thus
$$\mathrm{area}(\Sigma)\leq \mathrm{area}(\Sigma')={\bf L}(\Pi')\leq {\bf L}(S)\leq \sup\{ {\bf M}(\Phi'(x)):x\in I^5\}\leq{\rm area}(\Sigma).$$
This implies that  $S$ is a critical sequence  {and hence}, according to Theorem \ref{pitts.min.max}, we can choose $\Sigma' \in {\bf C}(S)$.

 {After passing to a subsequence, pick}  $x_i\in{\mathrm{dmn}}(\phi_i)$ so that $|\phi_i(x_i)|$ converges to $\Sigma'$ in the sense of varifolds. It follows
 {from Theorem \ref{discrete.sweepout} (i)} that,  {for some sequence $\{l_i\}_{i\in\N}$ tending to infinity, we have }  
\begin{multline*}
\mathrm{area}(\Sigma')=\lim_{i\to\infty}{\bf M}(\phi_i(x_i))\leq \lim_{i\to\infty}\sup\{{\bf M}(\Phi'(y)):\alpha\in I(5,l_i)_5, x_i,y\in \alpha\}\\
\leq \mathrm{area}(\Sigma).
\end{multline*}
Thus  {we  obtain from {Theorem \ref{discrete.sweepout} (ii)}  the existence of  a sequence $\{y_i\}_{i\in\N}$ in $I^5$ such that 
\begin{equation}\label{xieyi}
\lim_{i\to\infty}\mathcal{F}(\Phi'(y_i),\phi_i(x_i))=0\quad\mbox{and}\quad\lim_{i\to\infty}{\bf M}(\Phi'(y_i))={\mathcal W}(\Sigma) ={\rm area}(\Sigma).
\end{equation}}
From the definition of $\Phi'$ we have $\Phi'(y_i)=C'(v_i,t_i)$ for some sequence $(v_i,t_i)\in \overline B^4\times[-\pi,\pi]$ and we can extract a  subsequence $\{(v_i,t_i)\}_{i\in\N}$ converging to $(v,t)\in \overline B^4\times[-\pi,\pi]$.

Moreover,  \eqref{area.menor.C} implies that $C(v_i,t_i)=C'(v_i,t_i)$  and $|(v_i,t_i)|\geq \delta_1/2$ for all $i$ sufficiently large. 

\subsection{Lemma}\label{v.no.bordo} $w=T(v )\in B^4.$
\begin{proof}
%{The fact that $|(v_i,t_i)|\geq \delta_1/2$ for all $i$ sufficiently large implies  $v\neq 0$ and so $w=T(v)\neq 0$.}

Suppose $T(v)\in S^3$, i.e.,  $v\in S^3\cup \overline\Omega_{\varepsilon}$.  Theorem \ref{canonical.family.continuous} implies the existence of a geodesic sphere $S$  such that, after passing to a further subsequence, we have 
$$  \lim_{i\to\infty}{\mathcal F}(\Phi'(y_i),S)=\lim_{i\to\infty}{\mathcal F}(C'(v_i,t_i),S)=\lim_{i\to\infty}{\mathcal F}(C(v_i,t_i),S)=0.$$
If $\mathcal{F}(S)={\rm area}(S)=0$, we obtain from   Proposition \ref{convergence.sets} the existence of  $q\in S^3$ such that  for every $r$ we have
$$\Sigma_{(T(v_i),t_i)}\subset B_{r}(q)\quad\mbox{for all }  i\mbox{ sufficiently large}.$$
Thus, Theorem  \ref{no.concentration.mass} gives us that ${\bf M}(C(v_i,t_i))$ tends to zero. This is a contradiction and hence $\mathcal{F}(S)>0$.

Combining with \eqref{xieyi}, we obtained two subsequences $\{x_i\}_{i\in\N},$ $\{y_i\}_{i\in\N}$ in $I^5$ and a geodesic sphere $S$ with $\mathcal{F}(S)>0$ such that
$$ \lim_{i\to\infty}{\mathcal F}(\Phi'(y_i),S)=0,\quad\lim_{i\to\infty}{\mathcal F}(\Phi'(y_i),\phi_i(x_i))=0,$$
and
$$\lim_{i\to\infty}{\bf F}(|\phi_i(x_i)|,\Sigma')=0.$$
Lower semicontinuity of mass implies that  $S\llcorner(S^3\setminus \Sigma') =0$ and so {$S\subset \Sigma'$.} This is a contradiction because $S$ is a geodesic sphere and  {$\Sigma'$ has genus $g\geq 1$}.
\end{proof}

From the lemma above we have  $w=T(v)\in B^4$. Recall that ${\bf M}(C(v_i,t_i))={\bf M}(C'(v_i,t_i))$ tends to ${\rm area}(\Sigma)$,  and so we obtain from Theorem \ref{heintze.karcher} that either $t=0$ or $|t|=\pi$, because otherwise $\Sigma$ would be totally geodesic.

We argue $|t|=\pi$ does not occur. Choose $p\in \Sigma_{T(v)}$.
Theorem \ref{no.concentration.mass} tells us that there exists $r>0$ and $\delta'>0$ such that
\begin{equation}\label{equality.concentration}
{\rm area}(\Sigma_{(u,s)} \cap B_r(-p)) \leq \delta'<{\mathcal W}(\Sigma)
\end{equation}
for every $(u,s) \in B^4 \times [-\pi,\pi]$.

For all $i$ sufficiently large we have  $$\pi-r/2< |t_i|\leq \pi\quad\mbox{and}\quad d_H(\Sigma_{T(v_i)},\Sigma_{T(v)}) \leq r/2,$$ where $d_H$ denotes the Hausdorff distance. Hence
\begin{multline*}
d(\Sigma_{(T(v_i),t_i)},p)\geq d(\Sigma_{(T(v_i),t_i)},\Sigma_{T(v)})\\
\geq
 d(\Sigma_{(T(v_i),t_i)},\Sigma_{T(v_i)})-d_H(\Sigma_{T(v_i)},\Sigma_{T(v)})\\
=|t_i|-d_H(\Sigma_{T(v_i)},\Sigma_{T(v)})\geq \pi-r.
\end{multline*}
Thus $\Sigma_{(T(v_i),t_i)} \subset B_r(-p)$ and  \eqref{equality.concentration} contradicts the fact that ${\bf M}(C(v_i,t_i))$ tends to $\mathcal{W}(\Sigma).$

Thus $t=0$ and so, recalling that $|(v_i,t_i)|\geq \delta_1/2$ for all $i$ sufficiently large, we have $v\neq 0$ which means that $${\rm area}(\Sigma_w) ={\rm area}(\Sigma),\quad w=T(v){\in B^4\setminus\{0\}}.$$

\end{proof}

{Using Lemma \ref{conforme.igual} we now claim that $\Sigma$ must be totally geodesic.}

From formula (1.12) of \cite{montiel-ros}, by substituting $g=\frac{-2w}{(1+|w|^2)}$, we have that 
$$
{\rm area}(\Sigma_w) ={\rm area}(\Sigma)  -4\int_\Sigma \frac{\langle w, N(x)\rangle ^2}{|x-w|^4} d\Sigma.
$$
Thus  Lemma \ref{conforme.igual} implies that $\langle w, N(x)\rangle =0$ for every $x\in\Sigma$. 

On the other hand, let $h:S^3 \rightarrow \mathbb{R}$ be given by $h(x)=\langle x,w\rangle$. Because $\langle w, N(x)\rangle =0$ for every $x\in\Sigma$, 
the conformal vector field $V(x)=\nabla h(x)$ of $S^3$ satisfies $V(x) \in T_x\Sigma$ for all $x\in \Sigma$. This means $\Sigma$ is
invariant by the flow generated by $V$, but this is only possible if $\Sigma$ is totally geodesic. 

This is impossible because $\Sigma\in{\mathcal F}_1$ and thus ${\rm index}(\Sigma)\leq 5$. Hence we obtain from \cite{urbano} that $\Sigma$ is the Clifford torus up to ambient isometries.
\end{proof}

%is the Clifford torus up to isometries of $S^3$.

%%%%%%%%%%%%%%%%%%%%%%%%%%%%%%%%%%%%%%%%%%%%%%%%%%%%%
%%%%%%%%%%%%%%%%%%%%%%%%%%%%%%%%%%%%%%%%%%%%%%%%%%%%%%%%%%%

\section{Proof of Theorem A}\label{thma.section}

Let $\Sigma \subset S^3$ be an embedded closed surface of genus $g\geq 1$. We can assume $\mathcal{W}(\Sigma)<8\pi$.

Consider the min-max family $\Phi$ (see Definition \ref{Fi.family})  and the homotopy class $\Pi$ (see  Definition \ref{homotopy.class.sigma}) associated with $\Sigma$.  We have from Theorem \ref{modified.family} that all conditions required in Section \ref{bound.width.section} are met  and so we can apply Corollary \ref{width.8pi.discrete} to conclude the existence of a minimal surface $\Sigma'$ with genus $g\geq 1$ so that, from Theorem \ref{modified.family} (iii), we have
$${\rm area}(\Sigma')={\bf L}(\Pi)\leq \sup\{{\bf M}(\Phi(x)):x\in I^5\}\leq\mathcal{W}( \Sigma).$$
From Theorem B we have ${\rm area}( \Sigma')\geq 2\pi^2$ and so we have proved that $\mathcal{W}(\Sigma)\geq 2\pi^2$.

 Suppose  now $\mathcal{W}(\Sigma) =2\pi^2$. 
 
 \subsection{Lemma}\label{willmore.equality}{\em There is $w\in B^4$ so that
 ${\rm area}(\Sigma_w)=\mathcal{W}(\Sigma)=2\pi^2.$
 }
 \begin{proof}
 Consider the map $C$  given by Theorem \ref{canonical.family.continuous}
and the  $(5,{\bf M})$-homotopy sequence $S=\{\phi_i\}_{i\in\N}\in\Pi$ of mappings into $(\mathcal{Z}_2(S^3;{\bf M}), {\Phi_{|I^5_0}})$ given by  Theorem \ref{discrete.sweepout}.

Thus, from Theorem B,
\begin{multline*}
2\pi^2\leq \mathrm{area}(\Sigma')={\bf L}(\Pi)\leq {\bf L}(S)\leq \sup\{ {\bf M}(\Phi(x)):x\in I^5\}\\
\leq\mathcal{W}(\Sigma)=2\pi^2.
\end{multline*}
This implies that  $S$ is a critical sequence  {and hence}, according to Theorem \ref{pitts.min.max}, we can choose $\Sigma' \in {\bf C}(S)$.

 {After passing to a subsequence, pick}  $x_i\in{\mathrm{dmn}}(\phi_i)$ so that $|\phi_i(x_i)|$ converges to $\Sigma'$ in the sense of varifolds. It follows
 {from Theorem \ref{discrete.sweepout} (i)} that,  {for some sequence $\{l_i\}_{i\in\N}$ tending to infinity, we have }  
\begin{multline*}
\mathrm{area}(\Sigma')=\lim_{i\to\infty}{\bf M}(\phi_i(x_i))\\
\leq \lim_{i\to\infty}\sup\{{\bf M}(\Phi(y)):\alpha\in I(5,l_i)_5, x_i,y\in \alpha\}
\leq \mathcal{W}(\Sigma).
\end{multline*}
Thus  {we  obtain from {Theorem \ref{discrete.sweepout} (ii)}  the existence of  a sequence $\{y_i\}_{i\in\N}$ in $I^5$ such that 
$$\lim_{i\to\infty}\mathcal{F}(\Phi(y_i),\phi_i(x_i))=0\quad\mbox{and}\quad\lim_{i\to\infty}{\bf M}(\Phi(y_i))={\mathcal W}(\Sigma).$$}
From the definition of $\Phi$ we have $\Phi(y_i)=C(v_i,t_i)$ for some sequence $(v_i,t_i)\in \overline B^4\times[-\pi,\pi]$ and we can extract a  subsequence $\{(v_i,t_i)\}_{i\in\N}$ converging to $(v,t)\in \overline B^4\times[-\pi,\pi]$.

\subsection{Lemma} $w=T(v)\in B^4$.
 \begin{proof}
 If $v\in S^3\cup \overline\Omega_{\varepsilon}$ we argue like in  Lemma \ref{v.no.bordo}, and obtain two subsequences $\{x_i\}_{i\in\N},$ $\{y_i\}_{i\in\N}$ in $I^5$ and a geodesic sphere $S$ with $\mathcal{F}(S)>0$ such that
$$ \lim_{i\to\infty}{\mathcal F}(\Phi(y_i),S)=0,\quad\lim_{i\to\infty}{\mathcal F}(\Phi(y_i),\phi_i(x_i))=0,$$
and
$$\lim_{i\to\infty}{\bf F}(|\phi_i(x_i)|,\Sigma')=0.$$
Lower semicontinuity of mass implies that  $S\llcorner(S^3\setminus \Sigma') =0$ and so {$S\subset \Sigma'$.} This is a contradiction because $S$ is a geodesic sphere and  {$\Sigma'$ has genus $g\geq 1$}.
 \end{proof}
 Because ${\bf M}(C(v_i,t_i))$ tends to ${\mathcal W}(\Sigma)$, we combine the above lemma with Theorem \ref{heintze.karcher} to conclude that either $t=0$ or $|t|=\pi$. The same arguments as in Lemma \ref{conforme.igual} show that  $|t|=\pi$ does not occur. 
 
 Thus $t=0$, which means that $${\rm area}(\Sigma_w) ={\mathcal W}(\Sigma)=2\pi^2.$$
  \end{proof}

Lemma \ref{willmore.equality} implies at once Theorem A because in that case $\Sigma_w$ must be a minimal surface with genus $g\geq 1$ and area $2\pi^2$ and thus, by Theorem B, the Clifford torus up to ambient isometries. As a result, $\Sigma$ is the Clifford torus up to conformal transformations.

%%%%%%%%%%%%%%%%%%%%%%%%%%%%%%%%%%%%%%%%%%%%%%%%%%%%%
%%%%%%%%%%%%%%%%%%%%%%%%%%%%%%%%%%%%%%%%%%%%%%%%%%%%%%%%%%%

\part*{Part II. Technical work}\label{technical.work}
\medskip
\section{No area concentration}
\label{concentration.section}

 The goal of this section is to prove Theorem \ref{no.concentration.mass}:

\subsection*{Theorem \ref{no.concentration.mass}}
%\label{no.concentration.mass}
\textit{
For every  $\delta>0$, there exists $r>0$ such that
$$
{\rm area}(\Sigma_{(v,t)} \cap B_r(q)) \leq \delta\quad\mbox{for every }q\in S^3\mbox{ and } (v,t) \in B^4 \times [-\pi,\pi].
$$
}
The strategy for the proof is the following. From Remark \ref{canonical.remark} we know that $\Sigma_{(v,t)}$ is  contained in  the immersed surface
$$P_{(v,t)}=\psi_{(v,t)} \circ F_v:\Sigma \rightarrow S^3,$$
where
\begin{align}\label{mapP}
P_{(v,t)}(x) = &(\cos t)\, F_v(x)+ (\sin t)\, \frac{{DF_v}_{|x}(N)}{|{DF_v}_{|x}(N)|}\\ \notag
= & (\cos t)\, \left((1-|v|^2)\frac{x-v}{|x-v|^2}-v\right)\\ \notag
 &+ (\sin t)\, \left(N(x)+2\langle N(x),v\rangle \frac{x-v}{|x-v|^2}\right). \notag
\end{align}
It suffices to show that $P_{(v,t)}(\Sigma)$ has no area concentration, meaning that ${\rm area}(P_{(v,t)}(\Sigma)\cap B_r(q))$ is small if $r$ is small. The Jacobian of $P_{(v,t)}$ is uniformly bounded outside a tubular neighborhood of $\Sigma$ and so we need to analyze what happens when $v$ approaches $p\in \Sigma$.  We will do that by dividing $\Sigma$ in three regions: a tiny disc $D$ around $p$, where  $P_{(v,t)}(D)$ tends to a geodesic sphere and so there is no area concentration, a small annular region $N$, where $P_{(v,t)}(N)$ is forming a neck with area smaller than $\delta$ and so there is no area concentration, and the remaining region $\Sigma\setminus (D\cup N)$, where  the Jacobian of $P_{(v,t)}$ is uniformly bounded and so there is no area concentration.

Theorem \ref{no.concentration.mass} is proven at the end of this section.

\subsection{Preliminary results}

We derive three auxiliary results. Recall the definition of $\Lambda$ in Section \ref{associated.t}.

\subsection{Lemma}\label{estimates.P}\textit{
There exists a constant $C>0$ such that if   $v=\Lambda(p,s) \in B^4$ with  $|s|<C^{-1}$, then
\begin{eqnarray*}
|{DP_{(v,t)}}|(x) &\leq& C \left( 1+ \frac{|s|}{|s|^2 + |x-p|^2}\right),\\
|D^2 P_{(v,t)}|(x) &\leq& C \left( 1+ \frac{1}{|s|^2 + |x-p|^2}\right),
\end{eqnarray*}
for all $p, x\in\Sigma$.
}
\begin{proof}
For $v\in B^4$, consider $$h_v:\Sigma \rightarrow \mathbb{R}^4, \quad h_v(x) = \frac{x-v}{|x-v|^2}.$$ 
We claim the existence of $C_1>0$ such that if   $v=\Lambda(p,s) \in B^4$ with  $|s|<C_1^{-1}$, then
\begin{equation}\label{estimate.hv}
|D^kh_v|(x) \leq \frac{C_1}{(|s|^2+ |x-p|^2)^\frac{k+1}{2}}, \quad\mbox{for all }p,x\in \Sigma,\quad  k=0,1,2.
\end{equation}
There is $C_2>0$ so that,  for all $x,p \in\Sigma$,
\begin{equation}\label{estimate.hv1}
1-\langle x, p\rangle=\frac{|x-p|^2}{2}\quad\mbox{and}\quad|\langle x, N(p)\rangle|\leq C_2 |x-p|^2.
\end{equation}
Therefore, recalling
$$\Lambda(p,s)=(1-s_1)(\cos(s_2)p+\sin (s_2)N(p)),$$
we obtain
\begin{eqnarray*}
|x-v|^2 &=& 1-2\langle x,v\rangle +|v|^2\\
&=& 1-2\langle x, (1-s_1)(\cos(s_2) \, p + \sin(s_2)\, N(p)) \rangle +(1-s_1)^2\\
&=& 1-2\cos s_2 \langle x,p\rangle -2\sin s_2 \langle x,N(p)\rangle +2s_1\cos s_2\langle x,p\rangle \\
&&+ \, 2s_1\sin s_2  \langle x, N(p) \rangle +1-2s_1+s_1^2\\
&=& (1-s_1) (2-2\langle x,p\rangle) + s_1^2+s_2^2+O(s_1s_2^2+s_2^4+|s_2||x-p|^2).
\end{eqnarray*}
Thus, from \eqref{estimate.hv1} we see that we can find $C_3>0$ such that  
\begin{equation}\label{estimate.hv2}
|x-v|^2 \geq \frac12(|x-p|^2 +|s|^2)\quad\mbox{if}\quad |s|\leq C_3^{-1}.
\end{equation}
Direct computation shows that
$$|D^kh_v|(x)=O\left(|x-v|^{-(k+1)}\right)\quad\mbox{for }k=0,1,2$$
and thus the claim follows from \eqref{estimate.hv2}.

From \eqref{estimate.hv1} we have for $|s|<C_2^{-1}$
\begin{multline}\label{estimate.hv5}
|\langle N(x), v \rangle| = (1-s_1)|\cos s_2\, \langle N(x),p\rangle + \sin s_2\, \langle N(x),N(p)\rangle|\\
\leq C_2(s_2+|x-p|^2).
\end{multline}
Using the fact that $\langle DN_{|x}(Z),x\rangle=0$ for all $Z\in T_x\Sigma$, we have 
$$|\langle DN_{|x}(Z),v\rangle |=|\langle DN_{|x}(Z),v-x\rangle |= O(|Z||x-v|)$$
for all $x\in \Sigma$ and $Z\in T_x\Sigma$.
Finally we have  
\begin{equation}\label{estimate.hv6}
1-|v|^2 = 2s_1-s_1^2=O(s_1).
\end{equation}
Since 
\begin{eqnarray*}
P_{(v,t)}(x) &=& (\cos t)\, \left((1-|v|^2)h_v(x)-v\right)\\
&& \hspace{1cm}+ (\sin t)\, \left(N(x)+2\langle N(x),v\rangle h_v(x)\right),
\end{eqnarray*}
we use \eqref{estimate.hv}, \eqref{estimate.hv5}, and \eqref{estimate.hv6},  to conclude the existence of  $C>0$ such that if $|s|\leq 1/C$ then
$$
|DP_{(v,t)}|(x) \leq C \left( 1+ \frac{|s|}{|s|^2 + |x-p|^2}\right)
$$
and
$$
|D^2 P_{(v,t)}|(x) \leq C \left( 1+ \frac{1}{|s|^2 + |x-p|^2}\right).
$$

\end{proof}

Let $E_p:T_p\Sigma \rightarrow \Sigma \subset S^3$ be the exponential map of $\Sigma$ at $p$. 
We denote by $D_r(0)\subset T_p\Sigma$ the disk of radius $r$, centered
at the origin, and by $D_r(p) \subset \Sigma$ the geodesic disk of radius $r$, centered at $p$, with respect to the
induced metric.

\subsection{Lemma}\label{intermediate.mass}\textit{
For every $\delta>0$, there exist $L>0$ and $\alpha>0$ such
that the following holds: if $v=\Lambda(p,(s,ks))$ and  $(1+k^2)s^2\leq \alpha$, then
$$
\int_{D_\alpha(0)\setminus D_{L\sqrt{1+k^2}s}(0)} |{\rm Jac }(P_{(v,t)} \circ E_p)| \,dw \leq \delta.
$$
}
 
 \begin{proof}
It follows from Lemma \ref{estimates.P} that 
 \begin{multline*}
 \int_{D_\alpha(0)\setminus D_{L\sqrt{1+k^2}s}(0)} |{\rm Jac }(P_{(v,t)} \circ E_p)|dw\\
 \leq C_1 \int_{D_\alpha(0)\setminus D_{L\sqrt{1+k^2}s}(0)}  \left( 1+ \frac{|(s,ks)|}{|(s,ks)|^2 + |E_p(w)-p|^2}\right)^2 dw\\
  \leq C_2\alpha^2 + C_2 \int_{\mathbb{R}^2 \setminus D_{L\sqrt{1+k^2}s}(0)}  \left( \frac{|(s,ks)|}{|(s,ks)|^2 + |w|^2}\right)^2 dw,
 \end{multline*}
 for some constants $C_1,C_2>0$ depending only on $\Sigma$.
 
 After the change of variables $\tilde{w}=\frac{w}{|(1,k)|s}$, we obtain
 \begin{multline*}
  \int_{\mathbb{R}^2 \setminus D_{L\sqrt{1+k^2}s}(0)}  \left( \frac{|(s,ks)|}{|(s,ks)|^2 + |w|^2}\right)^2 dw\\
   = \int_{\mathbb{R}^2 \setminus D_{L}(0)}  \left( \frac{1}{1+ |\tilde{w}|^2}\right)^2 d\tilde{w}
   \leq \frac{\pi}{L^{2}}.
 \end{multline*}
Hence, if $\alpha>0$ is sufficiently small and $L>0$ is sufficiently large, we have 
 $$
\int_{D_\alpha(0)\setminus D_{L\sqrt{1+k^2}s}(0)} |{\rm Jac }(P_{(v,t)} \circ E_p)| \leq \delta.
$$
 \end{proof}

For every $x\in S^3$ denote by $\pi_x:S^3\setminus \{x\} \rightarrow \{x\}^\perp$ the stereographic projection centered at $x$:
 $$
 \pi_x(p) = x+ \frac{1}{1-\langle p,x\rangle}(p-x).
 $$
 The inverse of  $\pi_x$  is given by
$$
\pi_x^{-1}(w) = \frac{2}{1+|w|^2}(w-x)+x,\quad w\in \{x\}^\perp.
$$

 \subsection{Lemma}\label{renormalization}\textit{Let $(v_n,t_n) \in B^4\times [-\pi,\pi]$ with $v_n$ tending to $v=p \in \Sigma$.  After passing to a subsequence, write
 $$v_n=\Lambda(p_n,({s_n},k_n{s_n}))\quad\mbox{with}\quad \lim_{n\to\infty} k_n=k \in [-\infty, +\infty],$$
 and set
$$f_n(w)=P_{(v_n,t_n)} \circ E_{p_n} (\sqrt{1+k_n^2}s_nw).$$
 Then $f_n$ converges uniformly in $C^1_{{\rm loc}}$  to
$$
 f(w)=(\cos t+k\sin t) \left(\frac{\pi_x^{-1}(w)-x}{\sqrt{1+k^2}}\right)-(\cos t \, p- \sin t \, N(p)),
$$
 where $x=-\frac{1}{\sqrt{1+k^2}}\,p+\frac{k}{\sqrt{1+k^2}}N(p)\in S^3$.}

\subsection{Remark} 
\begin{enumerate}
\item With $x=-\frac{1}{\sqrt{1+k^2}}\,p+\frac{k}{\sqrt{1+k^2}}N(p)$, we have
$$  \frac{\pi_x^{-1}(w)-x}{\sqrt{1+k^2}}-p=\pi^{-1}_{-p}\left(\sqrt{1+k^2}w-kN(p)\right)\quad\mbox{for all }w\in T_p\Sigma.$$
Thus, as expected when $t=0$, 
$$f(T_p\Sigma)= \frac{\pi_x^{-1}\left(T_p\Sigma\right)-x}{\sqrt{1+k^2}}-p=\partial B_{\overline r_k}(\overline Q_{p,k}).$$
\item For the definition of $f_n$ to make sense we choose sequences of orthonormal sets $\{e_n^1,e_n^2\} \subset T_{p_n}\Sigma$ such that $e_n^i \rightarrow e^i\in T_p\Sigma$, $i=1,2$. Then we identify
$w =(w_1,w_2)\in \mathbb{R}^2$ with $w_1e_n^1+w_2e_n^2 \in T_{p_n}\Sigma$ for each $n$.
\end{enumerate}
\begin{proof} Note that both $s_n$ and $k_ns_n$  must tend to zero.
We have
\begin{align}\label{vn.renormalization}
 v_n&=(1-s_n) (\cos\, (k_ns_n) p_n +\sin\, (k_ns_n) N(p_n))\\ \notag
 &=p_n -s_np_n + k_ns_nN(p_n)+O((1+k_n^2)s_n^2),\notag
\end{align}
and
$$
E_{p_n}(\sqrt{1+k_n^2}s_nw) = p_n + \sqrt{1+k_n^2}s_nw + O((1+k_n^2)s_n^2|w|^2).
$$
Hence
\begin{multline*}
E_{p_n}(\sqrt{1+k_n^2}s_nw) -v_n=  \sqrt{1+k_n^2}s_nw +  s_np_n - k_ns_nN(p_n)\\
+O((1+k_n^2)s_n^2(1+|w|^2)),
\end{multline*}
and, using the fact that $\{w,p_n,N(p_n)\}$ is a orthogonal set of vectors, 
\begin{eqnarray*}
|E_{p_n}(\sqrt{1+k_n^2}s_nw) -v_n|^2 &=&  (1+k_n^2)s_n^2(1+|w|^2)\big(1+ O(\sqrt{1+k_n^2}s_n)\big).
\end{eqnarray*}
Therefore
\begin{multline}\label{quotient.renormalization}
\frac{E_{p_n}(\sqrt{1+k_n^2}s_nw) -v_n}{|E_{p_n}(\sqrt{1+k_n^2}s_nw) -v_n|^2} \\
= 
\frac{ \sqrt{1+k_n^2}w +  p_n - k_nN(p_n)}{(1+k_n^2)s_n(1+|w|^2)\big(1+ O(\sqrt{1+k_n^2}s_n)\big)}+O(1).
\end{multline}
Combining $1- |v_n|^2= 2s_n-s_n^2$ with \eqref{quotient.renormalization} we obtain
\begin{multline}\label{conformal.renormalization}
\lim_{n \rightarrow \infty} F_{v_n}\circ E_{p_n}(\sqrt{1+k_n^2}s_nw)\\
=\frac{2}{(1+|w|^2)}\left(\frac{w}{\sqrt{1+k^2}}+\frac{p}{1+k^2}-\frac{kN(p)}{1+k^2}\right) -p\\
=\frac{\pi_x^{-1}(w)-x}{\sqrt{1+k^2}}-p,
\end{multline}
where   $x=-\frac{1}{\sqrt{1+k^2}}\,p+\frac{k}{\sqrt{1+k^2}}N(p)$.

From the fact that $\langle N(x)-N(p_n),N(p_n)\rangle=O(|x-p_n|^2)$, we obtain from \eqref{estimate.hv1} and \eqref{vn.renormalization} that 
\begin{multline*}
\langle N(x),v_n\rangle=k_ns_n\langle N(x),N(p_n)\rangle+O(|x-p_n|^2+(1+k_n^2)s_n^2)\\
=k_ns_n+k_ns_n\langle N(x)-N(p_n),N(p_n)\rangle+O(|x-p_n|^2+(1+k_n^2)s_n^2)\\
=k_ns_n+O(|x-p_n|^2+(1+k_n^2)s_n^2).\\
\end{multline*}
Thus
$$ \langle N\circ E_{p_n}(\sqrt{1+k_n^2}s_nw),v_n\rangle=k_ns_n+O((1+k_n^2)s_n^2(1+|w|^2)),$$
which when combined with \eqref{quotient.renormalization} implies
\begin{multline}\label{normal.renormalization}
\lim_{n\rightarrow \infty} 2\langle N\circ E_{p_n}(\sqrt{1+k_n^2}s_nw),v_n\rangle\frac{E_{p_n}(\sqrt{1+k_n^2}s_nw) -v_n}{|E_{p_n}(\sqrt{1+k_n^2}s_nw) -v_n|^2}\\
=\lim_{n\rightarrow \infty} \left(2k_ns_n\frac{ \sqrt{1+k_n^2}w +  p_n - k_nN(p_n)}{(1+k_n^2)s_n(1+|w|^2)}\right)\\
=\frac{2}{(1+|w|^2)}\left(\frac{kw}{\sqrt{1+k^2}}+\frac{kp}{1+k^2}-\frac{k^2N(p)}{1+k^2}\right)\\
=\frac{k(\pi_x^{-1}(w)-x)}{\sqrt{1+k^2}},
\end{multline}
where   $x=-\frac{1}{\sqrt{1+k^2}}\,p+\frac{k}{\sqrt{1+k^2}}N(p)$.

From \eqref{mapP}, \eqref{conformal.renormalization}, and \eqref{normal.renormalization} we obtain that $f_n$ converges to $f$ pointwise

Fix $K>0$. It follows from Lemma \ref{estimates.P} that for every $w\in D_K(0)$,
\begin{multline*}
|Df_n(w)|\\
\leq   C\sqrt{1+k_n^2}s_n\left(1+ \frac{|(s_n,k_ns_n)|}{|(s_n,k_ns_n)|^2 + |E_{p_n}(\sqrt{1+k_n^2}s_nw)-p_n|^2}\right)\\
\leq C
\end{multline*}
and 
\begin{multline*}
|D^2f_n(w)|\\
\leq   C(1+k_n^2)s_n^2\left(1+ \frac{1}{|(s_n,k_ns_n)|^2 + |E_{p_n}(\sqrt{1+k_n^2}s_nw)-p_n|^2}\right)\\
\leq C.
\end{multline*}
Since we already know that $f_n$ converges to $f$ pointwise, the estimates above give $C^1$ convergence on compact subsets. 
\end{proof}

\subsection{Proof of Theorem \ref{no.concentration.mass}}
It suffices to show that for every $\delta>0$ and $q\in S^3$, we can find $r=r(q,\delta)$ so that 
$${\rm area}(\Sigma_{(v,t)} \cap B_r(q))\leq \delta$$
because, via a standard finite covering argument, we can then find $r$ independent of $q$.

Suppose this statement is false. There exist $q\in S^3$, $\delta>0$, and a sequence $(v_n,t_n) \in B^4 \times [-\pi,\pi]$ such that 
$$ {\rm area}(\Sigma_{(v_n,t_n)} \cap B_{1/n}(q))\geq \delta$$
for every $n \in \mathbb{N}$. By passing to a subsequence, we can assume $(v_n,t_n)$ converges to $(v,t) \in \overline{B}^4 \times [-\pi,\pi]$.

In what follows, we use repeatedly the fact that, from the area formula, 
$$
{\rm area}(\Sigma_{(v,t)} \cap B_r(q)) \leq \int_{P_{(v,t)}^{-1}(B_r(q))} |{\rm Jac\,}P_{(v,t)}| \, d\Sigma\quad\mbox{for all }r>0.
$$

If $v \in B^4$, then $P_{(v_n,t_n)}$ converges uniformly to $P_{(v,t)}$ in the $C^\infty$ topology and so we can find  $r>0$ such that, for all $n$ sufficiently large,
$$
\int_{P_{(v_n,t_n)}^{-1}(B_r(q))} |{\rm Jac\,}P_{(v_n,t_n)}| \, d\Sigma \leq \frac{\delta}{2}.
$$
This gives us a contradiction.

If $v\in S^3 \setminus \Sigma$, we see from \eqref{mapP} that  again $P_{(v_n,t_n)}$ converges uniformly, in the $C^\infty$ topology,
to some  $P_1:\Sigma \rightarrow S^3$. The proof proceeds as in the case $v\in B^4$.

Finally we have to consider the case $v=p\in \Sigma$. After passing to a subsequence, we can write
 $$v_n=\Lambda(p_n,({s_n},k_n{s_n}))\quad\mbox{with}\quad \lim_{n\to\infty} k_n=k \in [-\infty, +\infty].$$
 According to Lemma \ref{intermediate.mass}, we can choose $L>0$ and $\alpha >0$ so that
\begin{equation}\label{no.concentration.neck}
\int_{D_\alpha(0)\setminus D_{L\sqrt{1+k_n^2}s_n}(0)} |{\rm Jac }(P_{(v_n,t_n)} \circ E_{p_n})| \,dw \leq \frac{\delta}{6}
\end{equation}
if $n$ is sufficiently large.

Using Lemma \ref{estimates.P}, we extract a subsequence  $P_{(v_n,t_n)}$ that converges, $C^1$ uniformly, on $\Sigma \setminus D_{\alpha/4}(p)$ to some $C^1$
map $P_2:\Sigma \setminus D_{\alpha/4}(p) \rightarrow S^3.$ There exists $r_1>0$ such that
$$
\int_{P_2^{-1}(B_{2r_1}(q))} |{\rm Jac\, }P_2| \, d\Sigma \leq \frac{\delta}{12}
$$
and so, if $n$ is sufficiently large, we have 
\begin{equation}\label{no.concentration.outside}
\int_{P_{(v_n,t_n)}^{-1}(B_{r_1}(q))\setminus D_{\alpha/2}(p)} |{\rm Jac \,}P_{(v_n,t_n)}| \, d\Sigma \leq \frac{\delta}{6}.
\end{equation}
Consider $f_n: \overline{D}_{2L}(0) \rightarrow S^3$ given by $f_n(w)=P_{(v_n,t_n)} \circ E_{p_n}(\sqrt{1+k_n^2}s_nw)$. The sequence $f_n$ converges in the $C^1$ topology to $f$ given by Lemma \ref{renormalization} and hence we can find $r_2>0$ such that
$$
\int_{f^{-1}(B_{2r_2}(q))\cap \overline{D}_{2L}(0)} |{\rm Jac \,}f| \, dw \leq \frac{\delta}{12}.
$$
Therefore, if $n$ is sufficiently large, we have
\begin{equation}\label{no.concentration.bubble}
\int_{f_n^{-1}(B_{r_2}(q))\cap\overline{D}_{L}(0)} |{\rm Jac \,}f_n| \, dw \leq \frac{\delta}{6}.
\end{equation}

If $r=\min\{r_1,r_2\}$, we have the decomposition
\begin{eqnarray*}
&&\int_{P_{(v_n,t_n)}^{-1}(B_r(q))} |{\rm Jac}\,P_{(v_n,t_n)}| \, d\Sigma \\
&&= \int_{P_{(v_n,t_n)}^{-1}(B_r(q))\cap D_{L\sqrt{1+k_n^2}s_n}(p_n)} |{\rm Jac}\,P_{(v_n,t_n)}| \, d\Sigma\\
&&\hspace{0.5cm}+\int_{P_{(v_n,t_n)}^{-1}(B_r(q))\cap\big(D_\alpha(p_n) \setminus D_{L\sqrt{1+k_n^2}s_n}(p_n)\big)} |{\rm Jac}\,P_{(v_n,t_n)}| \, d\Sigma\\
&&\hspace{0.5cm}+\int_{P_{(v_n,t_n)}^{-1}(B_r(q)) \setminus D_\alpha(p_n)} |{\rm Jac}\,P_{(v_n,t_n)}| \, d\Sigma\\
&&\leq \int_{f_n^{-1}(B_{r_2}(q))\cap\overline{D}_{L}(0)} |{\rm Jac \,}f_n| \, dw\\
&&\hspace{.5cm}+\int_{D_\alpha(0)\setminus D_{L\sqrt{1+k_n^2}s_n}(0)} |{\rm Jac }\,(P_{(v_n,t_n)} \circ E_{p_n})| \,dw\\
&&\hspace{.5cm}+\int_{P_{(v_n,t_n)}^{-1}(B_{r_1}(q)) \setminus D_{\alpha/2}(p)} |{\rm Jac }\,P_{(v_n,t_n)}| \, d\Sigma.
\end{eqnarray*}

Using \eqref{no.concentration.neck}, \eqref{no.concentration.outside}, and \eqref{no.concentration.bubble} in the identity above we obtain
\begin{eqnarray*}
\int_{P_{(v_n,t_n)}^{-1}(B_r(q))} |{\rm Jac}\,P_{(v_n,t_n)}| \, d\Sigma\leq \frac{\delta}{2}
\end{eqnarray*}
for all $n$ sufficiently large. This is a contradiction.

%%%%%%%%%%%%%%%%%%%%%%%%%%%%%%%%%%%%%%%%%%%%%%%%%%%%%%%%
%%%%%%%%%%%%%%%%%%%%%%%%%%%%%%%%%%%%%%%%%%%%%%%%%%%%%%%%%%%

\section{Interpolation results: Continuous to discrete}\label{continuous.discrete}

{In this section we prove an interpolation theorem and use it to show Theorem \ref{discrete.sweepout}.}

{Assume} that we have  a continuous map in the flat topology
$$\Phi:I^{n}\rightarrow   \mathcal{Z}_2(M)$$
with  the following properties:
\begin{itemize}
\item $\Phi_{|I_0^n}$ is continuous in the ${\bf F}$-metric;
\item ${\bf L}(\Phi)=\sup\{{\bf M}(\Phi(x)):x\in I^{n}\}<+\infty;$
\item $\limsup_{r\to 0}{\bf m}(\Phi,r)=0.$
\end{itemize}

\subsection{Theorem}\label{flattomass}
 {\em There exist sequences of mappings 
$$\phi_i:I(n,k_i)_0\rightarrow  \mathcal{Z}_2(M),$$
$$\psi_i:I(1,k_i)_0\times I(n,k_i)_0\rightarrow  \mathcal{Z}_2(M)$$
with  $k_i < k_{i+1}$, $\psi_i([0],\cdot)=\phi_i,$  $\psi_i([1],\cdot)=(\phi_{i+1})_{|I(n,k_i)_0},$
and   sequences $\{\delta_i\}_{i\in\N}$ tending to zero and   $\{l_i\}_{i\in\N}$ tending to  infinity, such that
\begin{itemize}
\item[(i)] {for every $y\in I(n,k_i)_0$
$${\bf M}(\phi_i(y))\leq\sup \{{\bf M}(\Phi(x)):\alpha\in I(n,l_i)_n, x,y\in \alpha\}+\delta_i.$$
 In particular
 $$\max\{{\bf M}(\phi_i(x)):x\in I(n,k_i)_0\}\leq {\bf L}(\Phi)+\delta_i;$$
}
\item[(ii)] ${\bf f}(\psi_i)<\delta_i$;
\item[(iii)]
$$\sup\{{\mathcal{F}}(\psi_i(y,x)-\Phi(x))\,|\, y\in I(1,k_i)_0, x\in I(n,k_i)_0\}\leq \delta_i;$$
\item[(iv)] if  $x\in I_0(n,k_i)_0$ and $y\in I(1,k_i)_0$ we have
$${\bf M}(\psi_i(y,x))\leq {\bf M}(\Phi(x))+\delta_i.$$
\end{itemize}
Moreover, if $\Phi_{|\{0\}\times I^{n-1}}$ is continuous in the mass topology then we can
choose $\phi_i$ so that
$$\phi_i(x)=\Phi(x)\quad\mbox{for all } x\in {B(n,k_i)_0}.$$
}

\medskip

{For the reader's convenience we recall Theorem \ref{discrete.sweepout}.}
%Before proving Theorem \ref{flattomass}, we will derive an important consequence.
Let $$c=\frac{1}{3}(1,\ldots,1,0)\in I^{n-1}\times\{0\},$$ and $e_{n}$ be the coordinate vector corresponding to the  $x_{n}$-axis.

 We recall the following hypotheses for the continuous map in the flat topology $\Phi:I^{n}\rightarrow   \mathcal{Z}_2(M)$. Set $c=\frac{1}{3}(1,\ldots,1,0)\in I^{n-1}\times\{0\}.$
\begin{enumerate}
 \item[($A_0$)] $\Phi_{|I_0^n} \mbox{ is continuous in the ${\bf F}$-metric.}$
  \item[($A_1$)] $\Phi(I^{n-1}\times\{0\})=\Phi(I^{n-1}\times\{1\})=0.$
\item[($A_2$)] ${\bf L}(\Phi)=\sup\{{\bf M}(\Phi(x)):x\in I^{n}\}<+\infty.$
\item[($A_3$)] $\lim_{r\to 0}{\bf m}(\Phi,r)=0.$
\item[($A_4$)] The map $t\mapsto \Phi(c+tx_n)$, $0\leq t\leq 1$, defines a non-trivial class in $\pi_1(\mathcal{Z}_2(M;\mathcal{F}),\{0\})$.
 %there exists a family $\{U(t)\}_{0\leq t\leq 1}$ of open sets of $M$ of finite perimeter with
%\begin{itemize}
%item $U(0)=\emptyset$  and $U(1)=M$;
%\item $ \Phi(c+tx_n)=\partial [|U(t)|]\quad\mbox{for all }0\leq t\leq 1;$
%where $c=\frac{1}{3}(1,\ldots,1,0)\in I^{n-1}\times\{0\};$
%\item the map $t\to [|U(t)|]$ is continuous in the mass norm.
%\end{itemize}
\end{enumerate}

Then:
\subsection*{Theorem \ref{discrete.sweepout} }{\em  {Assume $\Phi$ satisfies hypotheses $(A_0)$--$(A_4)$}.

There exists an $(n,{\bf M})$-homotopy sequence of mappings into  $(\mathcal{Z}_2(M;{\bf M}),\Phi_{|I_0^n})$ $$\tilde{\phi}_i:I(n,k_i)\rightarrow  \mathcal{Z}_2(M),$$
 with the following properties:
\begin{itemize}
\item[(i)]  {There is a sequence $\{l_i\}_{i\in \N}$ tending to infinity such that for every sequence $x_i\in I(n,k_i)_0$ we have
$$\limsup_{i\to\infty}{\bf M}(\tilde\phi_i(x_i))\leq \limsup_{i\to\infty}\{{\bf M}(\Phi(x)):\alpha\in I(n,l_i)_n, x,x_i\in \alpha\}.$$
 In particular
  $${\bf L}(\{\tilde{\phi}_i\}_{i\in\N})\leq \sup\{{\bf M}(\Phi(x)):x\in I^{n}\};$$

 }
 \item[(ii)] {$$\lim_{i\to\infty}\sup\{{\mathcal{F}}(\tilde\phi_i(x)-\Phi(x))\,|\, x\in I(n,k_i)_0\}=0;$$}
\item[(iii)] The sequence of mappings 
$$v_i:I(1,k_i)_0\rightarrow  \mathcal{Z}_2(M;{\bf M}),\quad v_i(x)=\tilde{\phi}_i(c+xe_{n}),$$
 is a $(1,{\bf M})$-homotopy sequence of mappings into $(\mathcal{Z}_2(M;{\bf M}),\{0\})$  that belongs
 to a non-trivial element of $\pi_1^{\#}(\mathcal{Z}_2(M;{\bf M}),\{0\})$.
\end{itemize}
}

\begin{proof}
Let $\phi_i$, $\psi_i$, $\delta_i$ be given by Theorem \ref{flattomass}. It follows from property (iv) of   Theorem \ref{flattomass} and $(A_1)$ that 
\begin{equation}\label{mass.bound.zero.interpolation}
{\bf M}(\psi_i(y,x))\leq \delta_i
\end{equation}
for all $y\in I(1,k_i)_0$ and $x\in T(n,k_i)_0\cup B(n,k_i)_0$.

Define $\tilde{\psi}_i:I(1,k_i)_0\times I(n,k_i)_0\rightarrow  \mathcal{Z}_2(M)$ by $\tilde{\psi}_i(y,x)=0$ if $x\in T(n,k_i)_0\cup B(n,k_i)_0$ and $\tilde{\psi}_i(y,x)=\psi_i(y,x)$ otherwise. Define also $\tilde{\phi}_i(x)=\tilde{\psi}_i([0],x)$ for $x \in I(n,k_i)_0$. Note that ${\bf f}(\tilde{\psi_i}) < 2\delta_i$, by \eqref{mass.bound.zero.interpolation} and Theorem \ref{flattomass} part (ii). It follows that $\{\tilde{\phi}_i\}_{i \in \mathbb{N}}$ is a $(n,{\bf M})$-homotopy sequence of mappings into  $(\mathcal{Z}_2(M;{\bf M}),\Phi_{|I_0^n})$. Theorem \ref{flattomass} part (i)  {and part (iii) imply  Theorem \ref{discrete.sweepout} (i) and (ii), respectively.}

It remains to prove property  {(iii)} of Theorem \ref{discrete.sweepout}. 
Consider the auxiliary sequence 
$$\gamma_i:I(1,k_i)_0\rightarrow   \mathcal{Z}_2(M), \quad\gamma_i(x)=\Phi(c+xe_n)$$
and the continuous map in the flat topology
$$\gamma:[0,1] \rightarrow   \mathcal{Z}_2(M), \quad\gamma(x)=\Phi(c+xe_n).$$

Because $\Phi$ is continuous in the flat topology we have that
\begin{equation}\label{one.dimensional.zero.interpolation}
\lim_{i\to\infty}\sup_{\alpha\in I(1,k_i)_1}\{\mathcal{F}(\gamma_i(x)-\gamma_i(y)):x,y\in \alpha_0\}=0.
\end{equation}
 From that we get that $\tilde\gamma=\{\gamma_i\}_{i\in\N}$ is a $(1,\mathcal{F})$-homotopy sequence of mappings into $(\mathcal{Z}_2(M;\mathcal{F}),\{0\})$.  Furthermore, it follows from Theorem \ref{flattomass} (ii)  that  
$$\sup\{{\mathcal{F}}(\tilde{\phi}_i(x)-\Phi(x))\,:\, x\in I(n,k_i)_0\}\leq 2\delta_i.$$
This implies that  $v=\{v_i\}_{i\in\N}$ and $\tilde\gamma=\{\gamma_i\}_{i\in\N}$ are in the same $(1,\mathcal{F})$-homotopy class of mappings into $(\mathcal{Z}_2(M;\mathcal{F}),\{0\})$:
$$[v]=[\gamma]\in \pi_1^{\#}(\mathcal{Z}_2(M;\mathcal{F}),\{0\}).$$ 
 Since $\pi_1^{\#}(\mathcal{Z}_2(M;{\bf M}),\{0\})$, $\pi_1^{\#}(\mathcal{Z}_2(M;\mathcal{F}),\{0\})$, and  $\pi_1(\mathcal{Z}_2(M;\mathcal{F}),\{0\})$,  are all naturally isomorphic by \cite[Theorem 4.6]{pitts},  we get that  $[v]$ is non-trivial in $\pi_1^{\#}(\mathcal{Z}_2(M;{\bf M}),\{0\})$ if and only if $[\tilde\gamma ]$ is non-trivial in $\pi_1^{\#}(\mathcal{Z}_2(M;\mathcal{F}),\{0\})$, which occurs if and only if  $[\gamma]$ is non-trivial in $\pi_1^{}(\mathcal{Z}_2(M;\mathcal{F}),\{0\})$. The latter condition is assured by hypothesis $(A_4)$.

\end{proof}

The remainder of this section is devoted to the proof of Theorem \ref{flattomass}.

\subsection{Technical Results} We prove  two technical results which will be used in the proof  of Theorem \ref{flattomass}. 

The first proposition is an extension result. It  states  that if $T\in \mathcal{Z}_2(M)$  {and} $l,m\in\N$ are fixed, then we can find $k\in \N$, $k\geq l$, such that any map $\phi$ that sends  $I_0(m,l)_0$ into a small neighborhood of $T$ (with respect to the flat metric) can be extended to $I(m,k)_0$ in a way that the fineness of the extension $\tilde\phi$ and the maximum value of ${\bf M}(\tilde \phi)$ are not much bigger than  the fineness of $\phi$ and the maximum value of ${\bf M}(\phi)$, respectively. The issue of controlling the fineness of $\tilde \phi$ is nontrivial because {a priori} we only know that $\phi(I_0(m,l)_0)$ is close to $T$ in the flat metric, which is  weaker than the mass norm. A similar problem was addressed by Pitts in \cite[Lemma 3.7]{pitts}.  The fact, proven in Section \ref{concentration.section}, that there is no mass concentration will be used in the proof (although we think it might not be necessary). 

%It is possible this condition could be removed. 

Let $a(n)=2^{-4(n+2)^2-2}$ where $n\in\N$ is fixed.

\subsection{Proposition}\label{main.lemma.singleT}\textit{
Let $l,m\in \N$, with $m\leq n+1$, and let   $\delta, r, L$ be positive real numbers. Fix
$$T \in\mathcal{Z}_2(M)\cap\{S:{\bf M}(S)\leq 2L\}.$$ 
There exist $0<\varepsilon=\varepsilon(l,m,T, \delta,r, L)<\delta$ and $k=k(l,m,T,\delta, r, L)\in \N$ for which the following holds:
given $0<s<\varepsilon$ and 
$$\phi:I_0(m,l)_0\rightarrow  {\bf B}_{s}^{\mathcal{F}}(T)\cap\{S:{\bf M}(S)\leq 2L\}$$
with ${\bf m}(\phi,r)\leq \delta/4$, 
there exists
$$\tilde\phi:I(m,k)_0\rightarrow  {\bf B}_{s}^{\mathcal{F}}(T)$$
with
\begin{itemize}
\item[(i)] ${\bf f}(\tilde\phi) \leq \delta$ if $m=1$, and ${\bf f}(\tilde \phi)\leq m({\bf f}(\phi)+\delta) $ if $m \neq1$;
\item [(ii)] $\tilde\phi=\phi\circ {\bf n}(k,l)$ on $I_0(m,k)_0;$
\item[(iii)]$$\sup_{x\in I(m,k)_0}\{{\bf M}(\tilde\phi(x))\}\leq \sup_{x\in I_0(m,l)_0}\{{\bf M}(\phi(x))\}+\frac{\delta}{n+1};$$
\item[(iv)]${\bf m}(\tilde\phi,r)\leq 2({\bf m}(\phi,r)+a(n)\delta).$
\end{itemize}
}

\begin{proof}

We assume $m>1$ (the case $m=1$ is easier) and argue by contradiction. In this  case we can find 
$$\phi_k:I_0(m,l)_0\rightarrow  {\bf B}_{\varepsilon_k}^{\mathcal{F}}(T)\cap\{S:{\bf M}(S)\leq 2L\}$$
for each $k> \max\{ l,\delta^{-1}\}$, with $\varepsilon_k< 1/k$ and ${\bf m}(\phi_k,r)\leq \delta/4$, such that there is no extension $\tilde{\phi}_k$ of $\phi_k$ to $I(m,k)_0$  
satisfying (i)--(iv). 

The next lemma is a straightforward adaptation of \cite[Lemma 3.7]{pitts}.

\subsection{Lemma}\label{lemma.pitts.singleT} \textit{There exists $N\in \N$, $N \geq l$,  such that for  a subsequence $\{\phi_j\}$ we can find
$$\psi_j:I(1,N)_0\times I_0(m,l)_0\rightarrow  {\bf B}_{\varepsilon_j}^{\mathcal{F}}(T)$$
satisfying
\begin{itemize}
\item[(i)] ${\bf f}(\psi_j) \leq \delta$ if $m=1$ and ${\bf f}(\psi_j)\leq {\bf f}(\phi_j)+\delta$ if $m \neq1$;
\item[(ii)] $\psi_j([0],x)=\phi_j(x)$ and $\psi_j([1],x)=T$ for all $x\in I_0(m,l)_0$;
\item[(iii)] 
\begin{multline*}
\sup\{{\bf M}(\psi_j(y,x)):(y,x)\in I(1,N)_0\times I_0(m,l)_0\}\\
\leq \sup_{x\in  I_0(m,l)_0}\{{\bf M}(\phi_j(x))\}+\frac{\delta}{n+1};
\end{multline*}
\item[(iv)]${\bf m}(\psi_j,r)\leq  2({\bf m}(\phi_j,r)+a(n)\delta).$
\end{itemize}
}

\begin{proof}
Since the set of varifolds in $\mathcal{V}_2(M)$ with mass bounded above by $2L$ is compact in the weak topology, we can find  a subsequence $\{\phi_j\}$ of $\{\phi_k\}_{k\in \mathbb{N}}$  and  a map
$$V: I_0(m,l)_0\rightarrow \mathcal{V}_2(M)$$
so that 
$$ \lim_{j\to\infty}|\phi_j(x)|=V(x)\mbox{ as varifolds,}$$
for each $x\in  I_0(m,l)_0.$
Note that 
$$ \lim_{j\to\infty}\phi_j(x)=T\mbox{ as currents.}$$

Since the mass is  lower semicontinuous in the flat topology, and  since ${\bf m}(\phi_j,r)\leq \delta/4$, we have
\begin{equation}\label{3.7concentration}
||T||(B_r(p))\leq ||V(x)||(B_r(p))\leq {\bf m}(\phi_j,r)+a(n){\delta}<\frac{\delta}{3}
\end{equation}
for all $j$ sufficiently large,  $ p\in M,$ and $x\in I_0(m,l)_0.$ 

We can choose points $\{p_i\}_{i=1}^v$, and positive real numbers $\{r_i\}_{i=1}^v$, $r_i< r$, so that
$$B_{r_{i_1}}(p_{i_1})\cap B_{r_{i_2}}(p_{i_2})=\emptyset \quad\mbox{if }i_1\neq i_2,$$
and such that
\begin{equation}\label{3.7area.no.concentration}
||T||(B_{r_i}(p_i))\leq ||V(x)||(B_{r_i}(p_i))<\frac{\delta}{3}, 
\end{equation} 
\begin{equation}\label{3.7area.no.boundary}
||T||(\partial B_{r_i}(p_i))=||V(x)||(\partial B_{r_i}(p_i))=0, 
\end{equation} 
and 
\begin{equation}\label{3.7area.no.boundary2}
 ||V(x)||(M\setminus \cup_{i=1}^vB_{r_i}(p_i))<\frac{\delta}{3},
\end{equation}
for all $x\in I_0(m,l)_0$ and $i=1,\ldots,v$.  We can assume $v=3^N-1$ for some $N \in \mathbb{N}$ satisfying $N \geq l$.

From  \cite[Corollary 1.14]{almgren}, we get that there exists $Q_j(x)\in {\bf I}_3(M)$, for all $j$ sufficiently large and  $x\in I_0(m,l)_0$, such that
$$\partial Q_j(x)=\phi_j(x)-T, \quad {\bf M}(Q_j(x))=\mathcal{F}(\phi_j(x)-T).$$
In particular we have ${\bf M}(Q_j(x)) < \varepsilon_j < 1/j$.
%\begin{equation}\label{3.7mass.zero}
%\lim_{j\to\infty}{\bf M}(Q_j(x))=0.
%\end{equation}

For each $i=1,\ldots,v$, consider the distance function $d_i(x)=d(p_i,x)$. Using \cite[Lemma 28.5]{simon}, we  find a decreasing subsequence $\{r^j_i\}$ converging to $r_i$ with $r_i^j<r$ and such that the slices $\langle Q_j(x),d_i,r_i^j\rangle $ are in ${\bf I}_2(M)$ and satisfy
\begin{equation}\label{3.7slice.formula}
\langle Q_j(x),d_i,r_i^j\rangle=\partial (Q_j(x)\llcorner B_{r_i^j}(p_i))-(\phi_j(x)-T)\llcorner B_{r_i^j}(p_i)
\end{equation}
for every $x\in  I_0(m,l)_0$. 
Note that since $\lim_{j\to\infty}{\bf M}(Q_j(x))=0$, by the coarea formula we can choose $\{r^j_i\}$ such that 
\begin{equation}\label{3.7.area.slice}
\sum_{x\in I_0(m,l)_0}\sum_{i=1}^v {\bf M}(\langle Q_j(x),d_i,r_i^j\rangle)\leq a(n){\delta}<\frac{\delta}{2(n+1)}
\end{equation}
for every sufficiently large  $j$.
Furthermore, using \eqref{3.7area.no.concentration}, \eqref{3.7area.no.boundary}, \eqref{3.7area.no.boundary2}, and the lower semicontinuity of the mass functional, we get that
\begin{equation}\label{3.7area.inequality}
||\phi_j(x)||(B_{r^j_i}(p_i))<\frac{\delta}{3},  \quad ||T||(B_{r^j_i}(p_i))<\frac{\delta}{3},
\end{equation}
\begin{equation}\label{3.7area.inequality.complement}
||\phi_j(x)||{(}M\setminus \cup_{i=1}^vB_{r_i}(p_i))<\frac{\delta}{3}, \quad ||T||(M\setminus \cup_{i=1}^vB_{r_i}(p_i))<\frac{\delta}{3},
\end{equation}
and
\begin{equation}\label{3.7.lower.inequality}
(||T||-||\phi_j(x)||)(B_{r_i^j}(p_i))\leq \frac{\delta}{2(n+1)v}
\end{equation}
for every sufficiently large $j$, $i=1,\ldots,v,$  and  $x\in I_0(m,l)_0$.

We consider the map given by
\begin{eqnarray*}
\psi_j\left(\left[\frac{i}{3^N}\right],x\right)&=&\phi_j(x)-\sum_{a=1}^i\partial (Q_j(x)\llcorner B_{r^j_a}(p_a))\quad\mbox{ if }0\leq i\leq 3^N-1,\\
 \psi_{j}([1],x)&=&T,
 \end{eqnarray*}
defined on $I(1,N)_0\times I_0(m,l)_0$. 

%For simplicity, we denote the above map  by $\psi_k(i,x)$.

Note that
$$\psi_j\left(\left[\frac{i}{3^N}\right],x\right)-T=\partial (Q_j(x)\llcorner (M\setminus \cup_{a=1}^iB_{r^j_a}(p_a)),$$
from which it follows that  $\psi_j\left(\left[\frac{i}{3^N}\right],x\right)\in {\bf B}_{\varepsilon_j}^{\mathcal{F}}(T)$. From \eqref{3.7slice.formula},  we also have
\begin{multline}\label{3.7.formula2}
\psi_j\left(\left[\frac{i}{3^N}\right],x\right)=\phi_j(x)\llcorner (M\setminus \cup_{a=1}^iB_{r^j_a}(p_a))+\sum_{a=1}^iT\llcorner B_{r_a^j}(p_a)\\
-\sum_{a=1}^i\langle Q_j(x),d_a,r_a^j\rangle\llcorner B_{r_a^j}(p_a).
\end{multline}

  It follows from  \eqref{3.7.area.slice}, \eqref{3.7area.inequality}, \eqref{3.7area.inequality.complement}, and \eqref{3.7.formula2} that
\begin{eqnarray*}
&&{\bf M}\left(\psi_j\left(\left[\frac{i}{3^N}\right],x\right)-\psi_j\left(\left[\frac{i-1}{3^N}\right],x\right)\right) \\
&&\hspace{2cm}\leq \frac{\delta}{3}+ {\bf M}(\phi_j(x)\llcorner B_{r_i^j}(p_i))+ {\bf M}(T\llcorner B_{r_i^j}(p_i))<\delta
\end{eqnarray*}
for $1\leq i\leq v=3^N-1$, and
\begin{multline*}
{\bf M}\left(\psi_j\left(\left[1-\frac{1}{3^N}\right],x\right)-T\right)
\leq {\bf M}(\phi_j(x)\llcorner (M\setminus \cup_{a=1}^vB_{r^j_a}(p_a)))\\+ {\bf M}(T \llcorner(M\setminus \cup_{a=1}^vB_{r^j_a}(p_a)))+\frac{\delta}{3}<\delta.
\end{multline*}
If ${\bf d}(x,y)=1$, we also have
\begin{eqnarray*}
&&{\bf M}\left(\psi_j\left(\left[\frac{i}{3^N}\right],x\right)-\psi_j\left(\left[\frac{i}{3^N}\right],y\right)\right)\\
&&\hspace{1cm}\leq {\bf M}\left((\phi_j(x)-\phi_j(y))\llcorner(M\setminus \cup_{a=1}^iB_{r^j_a}(p_a))\right)+\frac{\delta}{2}\\
&&\hspace{1cm}\leq {\bf f}(\phi_j)+{\delta}.
\end{eqnarray*}
Hence ${\bf f}(\psi_j)\leq {\bf f}(\phi_j)+\delta$.

To prove Lemma \ref{lemma.pitts.singleT}(iii)  we use \eqref{3.7.area.slice}, \eqref{3.7.lower.inequality}, and \eqref{3.7.formula2}, to conclude
\begin{eqnarray*}
{\bf M}\left(\psi_j\left(\left[\frac{i}{3^N}\right],x\right)\right)&\leq& ||\phi_j(x)||(M\setminus \cup_{a=1}^iB_{r^j_a}(p_a))\\
&&\hspace{.5cm}+\sum_{a=1}^i ||T||(B_{r^j_a}(p_a))+\frac{\delta}{2(n+1)}\\ 
&\leq&  ||\phi_j(x)||(M)\\
&&\hspace{.5cm}+\sum_{a=1}^i (||T||-||\phi_j(x)||)(B_{r^j_a}(p_a))+\frac{\delta}{2(n+1)}\\
&\leq& ||\phi_j(x)||(M)+\frac{\delta}{n+1}.
\end{eqnarray*}
Finally, Lemma \ref{lemma.pitts.singleT}(iv) follows  from \eqref{3.7concentration}, \eqref{3.7.area.slice}, and \eqref{3.7.formula2}:
\begin{multline*}
\left|\left|\psi_j\left(\left[\frac{i}{3^N}\right],x\right)\right|\right|(B_r(p))\leq ||\phi_j(x)||(B_r(p))+||T||(B_r(p))+a(n)\delta\\
\leq   2{\bf m}(\phi_k,r)+2a(n){\delta}.
\end{multline*}
\end{proof}

In order to finish the proof of Proposition \ref{main.lemma.singleT}, we  will use Lemma \ref{lemma.pitts.singleT} to construct an extension $\tilde\phi_j$ for every sufficiently large $j$. 
This will  imply  a contradiction.

Define 
$$\hat\phi_j:I(1,N)_0\times I_0(m,N)_0\rightarrow  {\bf B}_{\varepsilon_j}^{\mathcal{F}}(T)$$
by
 $$\hat\phi_j (y,x)=\psi_j(y,{\bf n}(N,l)(x)).$$
 Recall that $S(m+1,N)_0=I(1,N)_0\times I_0(m,N)_0$. We extend $\hat\phi_j$ to  $$S(m+1,N)_0\cup T(m+1,N)_0$$ by setting it equal to $T$ on $T(m+1,N)_0$.   The extension $\tilde{\phi}_j:I(m,j)_0\rightarrow  {\bf B}_{\varepsilon_j}^{\mathcal{F}}(T)$ is defined by 
$$ \tilde\phi_j=\hat\phi_j\circ {\bf r}_m(N)\circ {\bf n}(j,N+q),$$
where ${\bf r}_m(N)$
 and   $q$ are as in Appendix \ref{appendix.map}.

%Replacing $\phi,\tilde \phi$ with $\phi_k,\tilde \phi_k$  in the statement of Proposition \ref{main.lemma.singleT}, we see that   Proposition \ref{main.lemma.singleT}(i) follows from %Lemma \ref{lemma.pitts.singleT} (i) and \eqref{prop2restriction},  Proposition \ref{main.lemma.singleT}(ii) follows from   Lemma \ref{lemma.pitts.singleT} (ii)  and 
%\eqref{prop1restriction},   Proposition \ref{main.lemma.singleT}(iii) follows   from Lemma \ref{lemma.pitts.singleT} (iii), and  Proposition \ref{main.lemma.singleT}(iv) follows from %Lemma \ref{lemma.pitts.singleT} (iv). This gives us a contradiction.

\end{proof}

The next result removes the dependence of $\varepsilon$ and $k$ on the parameters $l$ and $m$,  in Proposition \ref{main.lemma.singleT}. Roughly speaking, it says that with $T\in {\mathcal Z}_2(M)$ fixed we can find $k\in \N$ such that every map $\phi$ from $I_0(m,j)_0$ into a small neighborhood of $T$ (with respect to the flat metric) can be extended to a map $\tilde \phi$ from $I(m,k+j)_0$ into the same neighborhood of $T$ and having the same properties as the map constructed in Proposition \ref{main.lemma.singleT}.

The constant $b(n)$ mentioned below is universal.

\subsection{Proposition}\label{singleT}  \textit{Let  $\delta,r,L$ be positive real numbers, and  let
 $$T\in  \mathcal{Z}_2(M)\cap\{S:{\bf M}(S)\leq 2L\}.$$  There exist $0<\varepsilon=\varepsilon(T,\delta,r,L)<\delta$ and $k=k(T,\delta,r,L)\in \N$ for which  the following holds:
given $0<s<\varepsilon$,   $j,m\in \N$ with  $m\leq n+1,$ and 
$$\phi:I_0(m,j)_0\rightarrow  {\bf B}_s^{\mathcal{F}}(T)\cap\{S:{\bf M}(S)\leq 2L-\delta\}$$
 with $$2^{n+2}({\bf m}(\phi,r)+a(n)\delta)\leq \delta/4,$$
there exists
$$\tilde\phi:I(m,j+k)_0\rightarrow  {\bf B}_{s}^{\mathcal{F}}(T)$$
with 
\begin{itemize}
\item[(i)]  ${\bf f}(\tilde\phi) \leq \delta$ if $m=1$ and ${\bf f}(\tilde \phi)\leq b(n)({\bf f}(\phi)+\delta)$ if $m\neq 1$; 
\item[(ii)] $\tilde\phi=\phi\circ {\bf n}(k+j,j)$ on $I_0(m,k+j)_0$;
\item[(iii)] $$\sup_{x\in I(m,k+j)_0}\{{\bf M}(\tilde\phi(x))\}\leq \sup_{x\in I_0(m,j)_0}\{{\bf M}(\phi(x))\}+\delta;$$
\item[(iv)]${\bf m}(\tilde\phi,r)\leq 2^{n+2}({\bf m}(\phi,r)+a(n)\delta).$
\end{itemize}
}

\begin{proof}
Assume $m>1$ (the case $m=1$ is easier).
 Using  the notation of Proposition \ref{main.lemma.singleT}, set 
$$k_0=0, \quad k_1=k(0,1,T,{\delta},r,L), \quad k_i=k(k_{i-1},i,T,\delta,r, L),$$
where $i=1,\ldots,n+1$, 
and
$$\varepsilon=\min\{\varepsilon(k_{i-1},i,T,\delta,r,L):i=1,\ldots,n+1\}.$$
In what follows, we will apply Proposition \ref{main.lemma.singleT} to maps defined on  {vertices} of a $p$-cell $\alpha \in I(m,j)_p$, after
identifying $\alpha$ with $I^p$ through an affine map.

Let $V_p$ be the set of vertices of $I(m,j+k_p)$ that 
belong to the $p$-skeleton of $I(m,j)$, i.e.,  $V_p = \cup_{\alpha\in I(m,j)_{p}} \alpha(k_p)_0$. We say a map 
$$\phi_p:V_p \rightarrow {\bf B}_{s}^{\mathcal{F}}(T)\cap\{S:{\bf M}(S)\leq 2L\}$$  
 is a $p$-extension of $\phi$ if the following conditions are met:
\begin{enumerate}
%\item $\phi_p$ is defined on $I(m,j)_0$ and on $\alpha(k_p)_0$ for every $\alpha\in I(m,j)_{p}$.
\item $\phi_p(x)=\phi \circ{\bf n}(j+k_p,j)(x)$ for  $x \in I^m_0$.
\item If $p=1$,   we require ${\bf f}(\phi_1)\leq {\bf f}(\phi)+\delta.$
If $p>1$, we ask that there exists  $\phi_{p-1}$, a $(p-1)$-extension of $\phi$, so that $${\bf f}(\phi_p)\leq p({\bf f}(\phi_{p-1})+\delta).$$ 
\item $$ \sup_{ x\in\,V_p}\{{\bf M}(\phi_p(x))\} \leq \sup_{x\in I_0(m,j)_0}\{{\bf M}(\phi(x))\}+\frac{p\delta}{n+1}.$$
\item ${\bf m}(\phi_p,r)\leq 2^p{\bf m}(\phi,r)+2(2^p-1)a(n)\delta.$
\end{enumerate} 

We will now construct a $1$-extension $\phi_1$ of $\phi$. First fix $y\in I_0(m,j)_0$, and 
 define 
$$\phi_0:I(m,j)_0\rightarrow  {\bf B}_s^{\mathcal{F}}(T)\cap\{S:{\bf M}(S)\leq 2L-\delta\}$$
by $\phi_0(x)=\phi(x)$ if   $x \in I_0(m,j)_0$ and $\phi_0(x)=\phi(y)$ if $x \notin I_0(m,j)_0$. By applying
Proposition \ref{main.lemma.singleT} to $\phi_0$ in each 1-cell of $I(m,j)$, we get a map
$\tilde\phi_0:V_1\rightarrow  {\bf B}_{s}^{\mathcal{F}}(T)$. 

Let $\alpha \in I(m,j)_1$. If $\alpha$ is a $1$-cell of $I_0(m,j)$, we set 
$\phi_1=\phi \circ{\bf n}(j+k_1,j)$ on $\alpha(k_1)_0$. If $\alpha \notin I_0(m,j)$, we set
 $\phi_1=\tilde\phi_0$ on $\alpha(k_1)_0$.  The fact that $\phi_1$ is a $1$-extension of $\phi$ follows directly from the construction and Proposition  \ref{main.lemma.singleT}.

 %Given a $1$-cell $\alpha\in I(m,j)_1$, set ${\phi_1}$ on $\alpha(k_1)_0$ to be $\phi_{|\alpha_0} \circ{\bf n}(k_1,0)$ if %$\alpha$ is also a $1$-cell of $I_0(m,j)$. If not,   set ${\phi_1}$ on $\alpha(k_1)_0$ to be $\tilde \phi$ given by  %Proposition \ref{main.lemma.singleT}, which we apply   with $l=0, m=1$, and  $\phi$ replaced by $\phi_0$ restricted %to $\alpha_0$. 

\subsection{Lemma}\label{iterate1}\textit{
Given a $p$-extension  $\phi_p$  of $\phi$, we can find  a $(p+1)$-extension $\phi_{p+1}$ of $\phi$.
}
\begin{proof}
By applying
Proposition \ref{main.lemma.singleT} to $\phi_p$ in a $(p+1)$-cell $\alpha$ of $I(m,j)$, we get a map
$\tilde\phi_{p,\alpha}:\alpha(k_{p+1})_0 \rightarrow  {\bf B}_{s}^{\mathcal{F}}(T)$.  If $\alpha$ and $\overline{\alpha}$
are adjacent $(p+1)$-cells of $I(m,j)$, then property (ii) of Proposition \ref{main.lemma.singleT} guarantees that
$\tilde\phi_{p,\alpha}=\tilde\phi_{p,\overline{\alpha}}$ on $\alpha(k_{p+1})_0 \cap \overline{\alpha}(k_{p+1})_0$.
Therefore there exists $\tilde\phi_p: V_{p+1} \rightarrow {\bf B}_{s}^{\mathcal{F}}(T)$ such that 
$\tilde\phi_p=\tilde\phi_{p,\alpha}$ on $\alpha(k_{p+1})_0$ for each $\alpha$ of $I(m,j)_{p+1}$.

Note that  $\tilde\phi_{p}$ satisfies:
\begin{itemize}
%\item $\tilde \phi= \phi_{|\alpha_0}\circ{\bf n}(k_{p+1},0)$ if $\alpha\in  I_0(m,j)_{p+1}$.
%\item  $\tilde \phi= \phi_{p}\circ{\bf n}(k_{p+1},k_p)$ on $\alpha_0(k_{p+1})_0$. 
\item  ${\bf f}(\tilde \phi_p)\leq (p+1)({\bf f}(\phi_p)+\delta)$.
\item \begin{eqnarray*}
\sup_{x\in V_{p+1}}\{{\bf M}(\tilde \phi_p(x))\}&\leq& \sup_{x\in V_p}\{{\bf M}(\phi_p(x))\}+\frac{\delta}{n+1}\\
&\leq& \sup_{x\in I_0(m,j)_0}\{{\bf M}(\phi(x))\}+\frac{(p+1)\delta}{n+1}.
\end{eqnarray*}
\item  $${\bf m}(\tilde \phi_p,r)\leq 2({\bf m}(\phi_p,r)+a(n)\delta)\leq 2^{p+1}{\bf m}(\phi,r)+2(2^{p+1}-1)a(n)\delta.$$
\end{itemize}

Let $\alpha \in I(m,j)_{p+1}$. If $\alpha$ is a $(p+1)$-cell of $I_0(m,j)$, we set 
$\phi_{p+1}=\phi \circ{\bf n}(j+k_{p+1},j)$ on $\alpha(k_{p+1})_0$. If $\alpha \notin I_0(m,j)_{p+1}$, we set
 $\phi_{p+1}=\tilde\phi_p$ on $\alpha(k_{p+1})_0$. The fact that $\phi_{p+1}$ is a $(p+1)$-extension follows from the
 construction and the properties of $\tilde\phi_p$ listed above. 

\end{proof}

It follows by induction that there exists an $m$-extension $\phi_m:V_m \rightarrow {\bf B}_{s}^{\mathcal{F}}(T)\cap\{S:{\bf M}(S)\leq 2L\}$ of $\phi$. Note that $V_m=I(m,j+k_m)_0$.  To finish the proof of Proposition \ref{singleT}, we make
$k=k_m$ and $\tilde\phi=\phi_m.$
\end{proof}

\subsection{Proof of Theorem \ref{flattomass}}

The idea of the proof is the following. First, we cover $\{T:{\bf M}(T)\leq 2{\bf L}(\Phi)\}$ with a finite union of balls $\{B_i\}_{i=1}^N$ such that Proposition \ref{singleT} can be applied in each ball. Then, we choose $j$ large enough so that, for every $\alpha\in I(n,j)_n$,  $\Phi(\alpha_0)$ belongs to some ball $B_i$. Finally, we use Proposition \ref{singleT} to first construct $\phi$ along $3^k$ subdivisions of $1$-cells in $I(n,j)$, then along $3^{2k}$ subdivisions of $2$-cells of $I(n,j)$, and argue inductively until we have constructed $\phi$ defined on $I(n,j+nk)_0$.
Some care is in order to make sure that at every step of the inductive construction  the hypotheses of Proposition \ref{singleT} are still satisfied. The procedure is straightforward but slightly long and tedious.

Choose $\delta$, $r$ small so that 
\begin{equation}\label{flattomass.concentra0}
L= {\bf L}(\Phi)<2L-2(n+1)\delta \quad\mbox{and}\quad {\bf m}(\Phi,r)<a(n){\delta}.
\end{equation}
Compactness of $\mathcal{Z}_2(M)\cap\{T:{\bf M}(T)\leq 2L\}$ in the flat topology implies  we can cover this set with a finite number of balls ${\bf B}^{\mathcal{F}}_{\varepsilon_i}(T_i)$, $i=1,\ldots,N$, where $$T_i\in \mathcal{Z}_2(M)\cap\{T:{\bf M}(T)\leq 2L\} \quad\mbox{and}\quad \varepsilon_i=\frac{{\varepsilon(T_i,\delta,r,L)}}{9n+4}.$$ Here we use the notation of Proposition \ref{singleT}. We can assume that $\varepsilon_1<\ldots<\varepsilon_N$. Note that $(9n+4)\varepsilon_N<\delta.$ Let  $k_i=k(T_i,\delta,r,L)$ denote the constant  given by Proposition \ref{singleT}, and let  $k=\max\{k_i\}_{i\in \{1,\ldots,N\}}$.

Choose $j$ sufficiently large so that for all $\alpha\in I(n,j)_{n}$ and $\beta\in I_0(n,j)_{n-1},$  we have 
\begin{equation}\label{dist5}
\sup_{x,y\in \alpha}\{\mathcal{F}(\Phi(x)-\Phi(y))\}<\varepsilon_1,
\end{equation}
and
\begin{equation}\label{dist6}
\sup_{x,y\in \beta}|{\bf M}(\Phi(x))-{\bf M}(\Phi(y))|<\delta.
\end{equation}
Additionally, if $\Phi_{| \{0\}\times I^{n-1}}$ is continuous in the mass norm, we also require that for all $\gamma\in [0]\otimes I(n-1,j)$ we have 
\begin{equation}\label{dist7}
\sup_{x,y\in \gamma}\{{\bf M}(\Phi(x)-\Phi(y))\}<\delta.
\end{equation}

Consider the function
$$ {\bf c}:I(n,j)\rightarrow \{1,\ldots,N\}$$
given by 
$${\bf c}(x)=\max\{i:\Phi(x)\in {\bf B}^{\mathcal{F}}_{\varepsilon_i}(T_i)\}\quad\mbox{if }x\in I(n,j)_0,$$
and 
$$ {\bf c}(\alpha)=\max\{{\bf c}(x):x\in \alpha_0\}\quad\mbox{if }\alpha\in I(n,j)_p.$$
The key property of ${\bf c}$ is described below.

\subsection{Lemma}\label{cfunction}\textit{
Let  $\alpha\in I(n,j)$. Then $\Phi(x)\in  {\bf B}^{\mathcal{F}}_{2\varepsilon_{{\bf c}(\alpha)}}(T_{{\bf c}(\alpha)})$ for every $x\in \alpha$. If $x\in \alpha_0$, we also have
 $$T_{{\bf c}(x)}\in  {\bf B}^{\mathcal{F}}_{3\varepsilon_{{\bf c}(\alpha)}}(T_{{\bf c}(\alpha)}).$$
}

\begin{proof}
There exists a vertex $y\in\alpha_0$ such that ${\bf c}(y)={\bf c}(\alpha)$. Hence, by \eqref{dist5}, we know that
\begin{equation*}\label{cfunction1}
\mathcal{F}(\Phi(x)-\Phi(y))<\varepsilon_1\leq \varepsilon_{{\bf c}(\alpha)}.
\end{equation*}
Furthermore, from the definition of ${\bf c}$, we get
\begin{equation*}\label{cfunction2}
\Phi(y)\in {\bf B}^{\mathcal{F}}_{\varepsilon_{{\bf c}(\alpha)}}(T_{{\bf c}(\alpha)}).
\end{equation*}
Hence $\Phi(x)\in  {\bf B}^{\mathcal{F}}_{2\varepsilon_{{\bf c}(\alpha)}}(T_{{\bf c}(\alpha)})$ for every $x\in \alpha$.

If $x\in \alpha_0$, we also have $\Phi(x)\in {\bf B}^{\mathcal{F}}_{\varepsilon_{{\bf c}(x)}}(T_{{\bf c}(x)})$.
The lemma follows from the triangle inequality and   the fact that ${\bf c}(x)\leq {\bf c}(\alpha)$.
\end{proof}

Let $V_p$ be the set of vertices of $I(n,j+pk)$ that belong to the $p$-skeleton of $I(n,j)$, i.e., $V_p=\cup_{\alpha\in I(n,j)_p} \alpha(pk)_0.$ In particular, $V_n=I(n,j+nk)_0$.
We say a map 
$$\phi_p:V_p \rightarrow \mathcal{Z}_2(M)$$
 is a $p$-extension of $\Phi$ if the following conditions are met:
\begin{enumerate}
%\item $\phi_p$ is defined on $I(n,j)_0$ and on $\alpha(pk)_0$ for every $\alpha\in I(n,j)_p$.
\item   If $p=1$,   we require that ${\bf f}(\phi_1)\leq \delta.$ If $p>1$, we ask that there exists a $(p-1)$-extension $\phi_{p-1}$ of $\Phi$ so that $${\bf f}(\phi_p)\leq b(n)({\bf f}(\phi_{p-1})+\delta).$$ 
\item For every $\alpha\in I(n,j)_p$ we have
$$
\sup_{x\in \alpha(pk)_0}\{{\bf M}(\phi_p(x))\} \leq\sup_{x\in \alpha_0}\{{\bf M}(\Phi(x))\}+p\delta<2L-\delta.
$$
\item  For every $\alpha\in I(n,j)_q$ with $q\leq p$ we have
$$\phi_p(\alpha(pk)_0)\in {\bf B}^{\mathcal{F}}_{3p\varepsilon_{{\bf c}(\alpha)}}(T_{{\bf c}(\alpha)}).
$$
\item $
{\bf m}(\phi_p,r)\leq 2^{p(n+2)}(p+1)a(n)\delta.
$
\end{enumerate}

We start by constructing a $1$-extension of $\Phi$.  In what follows, we will apply Proposition \ref{singleT} to maps defined on  {vertices} of a $p$-cell $\alpha \in I(n,j)_p$, after
identifying $\alpha$ with $I^p$ through an affine map.

Let $\phi_0: I(n,j)_0 \rightarrow \mathcal{Z}_2(M)$ be the restriction of $\Phi$ to $I(n,j)_0$. Given a $1$-cell  $\alpha\in I(n,j)$,  we  have from Lemma \ref{cfunction} that
$$
\phi_0(\alpha_0) \subset {\bf B}^{\mathcal{F}}_{3\varepsilon_{{\bf c}(\alpha)}}(T_{{\bf c}(\alpha)}).
$$
By applying Proposition \ref{singleT} to $\phi_0$ on $\alpha$, with $T=T_{c(\alpha)}$, we get a map $\tilde\phi_{0,\alpha}:\alpha(k_{c(\alpha)}) \rightarrow {\bf B}^{\mathcal{F}}_{3\varepsilon_{{\bf c}(\alpha)}}(T_{{\bf c}(\alpha)}).$ Since 
$\tilde\phi_{0,\alpha}(x)=\phi_0(x)$ for $x \in \alpha_0$, the map
 ${\phi_1}:V_1 \rightarrow \mathcal{Z}_2(M)$ given by $\phi_1=\tilde\phi_{0,\alpha} \circ {\bf n}(j+k,j+k_{c(\alpha)})$ on  $\alpha(k)_0$, $\alpha \in I(n,j)_1$, is well-defined. It follows directly from Proposition \ref{singleT} that $\phi_1$ is a $1$-extension of $\Phi$.

\subsection{Lemma}\label{iterate2} \textit{Assume $1\leq p\leq n-1$. Given  a  $p$-extension $\phi_p$ of $\Phi$, we can find a  $(p+1)$-extension $\phi_{p+1}$ of $\Phi$.
}

\begin{proof}
Given   $\alpha\in I(n,j)_{p+1}$, we have from condition (3) and Lemma \ref{cfunction} that
$$\phi_p(\alpha_0(pk)_0)\in {\bf B}^{\mathcal{F}}_{3(p+1)\varepsilon_{{\bf c}(\alpha)}}(T_{{\bf c}(\alpha)}).
$$
Because of  conditions (2) and (4) we can apply
Proposition \ref{singleT} to $\phi_p$ in $\alpha$, with $j=pk$, $m=p+1$,  $T=T_{{\bf c}(\alpha)}$, and get a map
$$\tilde\phi_{p,\alpha}:\alpha(pk+k_{c(\alpha)})_0 \rightarrow   {\bf B}^{\mathcal{F}}_{3(p+1)\varepsilon_{{\bf c}(\alpha)}}(T_{{\bf c}(\alpha)}).$$ By property (ii) of Proposition \ref{singleT} we get that $$\tilde\phi_{p,\alpha}=\phi_p \circ {\bf n}(j+pk+k_{c(\alpha)},j+pk)$$
on the $p$-faces of $\alpha$. Hence the map $\phi_{p+1}: V_{p+1} \rightarrow \mathcal{Z}_2(M)$ given by
$$\phi_{p+1}=\tilde\phi_{p,\alpha}\circ {\bf n}(j+(p+1)k,j+pk+k_{c(\alpha)})$$ on $\alpha((p+1)k)_0$, $\alpha\in I(n,j)_{p+1}$, is well-defined. 

Note that $\phi_{p+1}$ satisfies:
\begin{itemize}
\item if $\alpha\in I(n,j)_{p+1}$, then 
$$\phi_{p+1}=\phi_p\circ{\bf n}(j+(p+1)k,j+pk)\quad {\rm on \ } \alpha_0((p+1)k)_0;$$
\item${\bf f}(\phi_{p+1})\leq b(n)({\bf f}(\phi_{p})+\delta)$;
\item \begin{multline}\label{inclusion3.thm.flattomass}
\sup_{x\in \alpha((p+1)k)_0}\{{\bf M}( \phi_{p+1}(x))\} \leq \sup_{x\in \alpha_0(pk)_0}\{{\bf M}(\phi_p(x))\}+\delta\\
\leq \sup_{x\in \alpha_0}\{{\bf M}(\Phi(x))\}+(p+1)\delta;
\end{multline}
\item if $\alpha\in I(n,j)_{p+1}$, then 
$$\phi_{p+1}(\alpha((p+1)k)_0)\in {\bf B}^{\mathcal{F}}_{3(p+1)\varepsilon_{{\bf c}(\alpha)}}(T_{{\bf c}(\alpha)});
$$
\item \begin{equation*}\label{flattomass.concentra2}
{\bf m}(\phi_{p+1},r)\leq 2^{n+2}(2^{p(n+2)}(p+1)a(n)\delta+a(n)\delta)\leq 2^{(p+1)(n+2)}(p+2)a(n)\delta.
\end{equation*}
 \end{itemize}

Furthermore, if  $\beta\in I(n,j)_q$ with $q\leq p$, we can find  $\alpha\in I(n,j)_{p+1}$ such that $\beta$ is a face of $\alpha$. Hence, by the first property of $\phi_{p+1}$ listed above, 
$$\phi_{p+1}(\beta((p+1)k)_0)=\phi_p(\beta(pk)_0)\subset {\bf B}^{\mathcal{F}}_{3p\varepsilon_{{\bf c}(\beta)}}(T_{{\bf c}(\beta)}).$$

We conclude that $\phi_{p+1}$ is a $(p+1)$-extension of $\Phi$.
 \end{proof}
 
 By applying Lemma \ref{iterate2} inductively,  we obtain  the existence of an $n$-extension $\phi_\delta=\phi_n$ of $\Phi$:
$$\phi_{\delta}:I(n,j+nk)_0\rightarrow \mathcal{Z}_2(M).$$
The map $\phi_\delta$ has   the following properties: 
\begin{itemize}
\item[a)] ${\bf f}(\phi_{\delta})\leq c(n)\delta$ for some universal constant $c(n)$;
\item[b)]  {for every $x\in I(n,j+nk)_0$ 
\begin{equation}\label{mass.thm.pointwise} 
{\bf M}( \phi_{\delta}(x))\leq \sup\{{\bf M}(\Phi(y)):\alpha\in I(n,j)_0, x,y\in \alpha\}+n\delta.
\end{equation}
In particular} 
\begin{equation}\label{mass.thm.flattomass}
\sup_{x\in I(n,j+nk)_0}\{{\bf M}(\phi_{\delta}(x))\}\leq{\bf L}(\Phi)+n\delta;
\end{equation}
\item[c)]
\begin{equation*}
{\bf M}(\phi_{\delta}(x))\leq{\bf M}(\Phi(x))+(n+1)\delta\quad\mbox{for all}\quad x\in I_0(n,j+nk)_0;
\end{equation*}
\item[d)] \begin{equation*}\label{inclusion4.thm.flattomass}
{\bf m}(\phi_{\delta},r)\leq 2^{n(n+2)}(n+1)a(n){\delta};
\end{equation*}
\item[e)] For every $\alpha\in I(n,j)_p$ with $p\leq n$
\begin{equation*}\label{inclusion2.thm.flattomass}
\phi_{\delta}(\alpha(nk)_0)\in {\bf B}^{\mathcal{F}}_{3n\varepsilon_{{\bf c}(\alpha)}}(T_{{\bf c}(\alpha)}).
\end{equation*}
\end{itemize}
We note that property c) follows from  \eqref{dist6} and \eqref{inclusion3.thm.flattomass}. Furthermore, Lemma \ref{cfunction}, \eqref{dist5}, and property e)  imply that
\begin{equation}\label{inclusion5.thm.flattomass}
\sup\{\mathcal{F}(\phi_{\delta}(x)-\Phi(x)): x\in I(n,j+nk)_0\}\leq 3(n+1)\varepsilon_N<\delta.
\end{equation}

Before proceeding with the construction, we need one more definition. A map $$\bar \phi:I(n,\bar k)_0\rightarrow  \mathcal{Z}_2(M)\cap\{S:{\bf M}(S)\leq 2L\}$$
is called a $(n,\delta,\bar k)$-extension of $\Phi$ if it satisfies 
\begin{itemize}
\item[a')] ${\bf f}(\bar \phi)\leq c(n)\delta$;
\item[b')] \begin{equation*}\label{mass.thm.flattomass.2}
\sup_{x\in I(n,\bar k)_0}\{{\bf M}(\bar \phi(x))\}\leq{\bf L}(\Phi)+n\delta;
\end{equation*}
\item[c')]
\begin{equation*}
{\bf M}(\bar \phi(x))\leq{\bf M}(\Phi(x))+(n+1)\delta\quad\mbox{for all}\quad x\in I_0(n,\bar k)_0;
\end{equation*}
\item[d')] \begin{equation*}
{\bf m}(\bar \phi,r)\leq 2^{n(n+2)}(n+1)a(n){\delta};
\end{equation*}
\item[e')] \begin{equation*}\label{hyp.homotopic.sequence}
\sup\{\mathcal{F}(\Phi(x)-\bar \phi(x)):x\in I(n,\bar k)_0\}<\varepsilon_1.
\end{equation*}
\end{itemize}

The constant $d(n)$ mentioned below is universal.

\subsection{Proposition}\label{homotopic.sequence}  \textit{Let $\bar \phi$ be an $(n,\delta,\bar k)$-extension of $\Phi$, with $\bar k\geq j+nk$. Then there exists
$$\psi:I(1,\hat k)_0\times I(n,\hat k)_0\rightarrow  \mathcal{Z}_2(M),$$
with $\hat k=(n+1)k+\bar k$, such that
  $$\psi([0],\cdot)=\phi_{\delta}\circ{\bf n}(\hat k, j+nk), \quad  \psi([1],\cdot)=\bar \phi\circ{\bf n}(\hat k, \bar k),$$
and
\begin{itemize}
\item[(i)] ${\bf f}(\psi)<d(n)\delta$;
\item[(ii)]
$$\sup\{{\mathcal{F}}(\psi(y,x)-\Phi(x))\,:\, y\in I(1,\hat k)_0, x\in I(n,\hat k)_0\}\leq \delta;$$
\item[(iii)]
$${\bf M}(\psi(y,x))\leq {\bf M}(\Phi(x))+2(n+2)\delta\quad \mbox{for all } (y,x)\in I(1,\hat k)_0\times I_0(n,\hat k)_0.$$
\end{itemize}
}

\begin{proof}
Let $\tilde \phi_{\delta}=\phi_{\delta}\circ{\bf n}(\bar k,nk+j)$ on $I(n,\bar k)_0$. We also define 
$${\bf \bar c}: I(n,\bar k)\rightarrow\{1,\ldots,N\}$$  
by 
$${\bf \bar c}(\alpha)=\sup\{{\bf c}(\beta):\beta \in I(n,j)\mbox{ and }\alpha\cap \beta\neq \emptyset \}.$$
Note that ${\bf \bar c}(\alpha) \leq {\bf \bar c}(\alpha')$ if $\alpha \subset \alpha'$.
The next lemma is similar to Lemma \ref{cfunction}.
\subsection{Lemma}\label{lemm.incl}\textit{
Let  $\alpha\in I(n,\bar k)$. We have
$$\tilde \phi_\delta(\alpha_0),\,\bar\phi(\alpha_0)\subset  {\bf B}_{(3n+4)\varepsilon_{{\bf \bar c}(\alpha)}}^{\mathcal{F}}(T_{{\bf \bar c}(\alpha)}),$$
and
$$\Phi(x)\in  {\bf B}^{\mathcal{F}}_{3\varepsilon_{{\bf \bar c}(\alpha)}}(T_{{\bf \bar c}(\alpha)})$$
for every $x\in \alpha_0$. In particular, if $\alpha,\alpha' \in I(n,\bar k)$ satisfy $\alpha \subset \alpha'$ then
$$
T_{{\bf \bar c}(\alpha)}\in  {\bf B}^{\mathcal{F}}_{6\varepsilon_{{\bf \bar c}(\alpha')}}(T_{{\bf \bar c}(\alpha')}).
$$
%$$\Phi(x),\, T_{{\bf \bar c}(x)}\in  {\bf B}^{\mathcal{F}}_{3\varepsilon_{{\bf \bar c}(\alpha)}}(T_{{\bf \bar c}(\alpha)}).$$
}

\begin{proof}
Let $\eta \in I(n,j)$ with $\alpha \subset \eta$. From the definition of ${\bf \bar c}$, there exists $\beta \in  I(n,j)$ with  $\alpha \cap \beta \neq \emptyset$  such that ${\bf \bar c}(\alpha)={\bf c}(\beta)$. In particular, $c(\eta) \leq c(\beta)$ and $\beta \cap \eta \neq \emptyset$. It follows from  Lemma \ref{cfunction} and property e) that
\begin{equation*}
\tilde \phi_{\delta}(\alpha_0)\subset \phi_{\delta}(\eta(nk)_0)  \subset {\bf B}^{\mathcal{F}}_{3n\varepsilon_{c(\eta)}}(T_{c(\eta)}),
\end{equation*}
and
$$\Phi(y)\in {\bf B}^{\mathcal{F}}_{2\varepsilon_{c(\beta)}}(T_{c(\beta)})
\cap  {\bf B}^{\mathcal{F}}_{2\varepsilon_{c(\eta)}}(T_{c(\eta)})\quad\mbox{for all }y\in \beta \cap \eta.$$
Hence $T_{c(\eta)} \in  {\bf B}^{\mathcal{F}}_{4\varepsilon_{c(\beta)}}(T_{c(\beta)})$. It follows that
$\tilde \phi_{\delta}(\alpha_0)\subset {\bf B}^{\mathcal{F}}_{(3n+4)\varepsilon_{c(\beta)}}(T_{c(\beta)})$.

Let $y\in \alpha_0 \cap \beta \subset \eta$. From property e') and (\ref{dist5}), we get $\bar \phi(x) \in {\bf B}^{\mathcal{F}}_{2\varepsilon_1}(\Phi(y))$ and $\Phi(x) \in {\bf B}^{\mathcal{F}}_{\varepsilon_1}(\Phi(y))$ for each $x\in \alpha_0$. Therefore
\begin{equation*}\label{incl.prop.homotopic.sequence}
\bar\phi(\alpha_0) \subset {\bf B}^{\mathcal{F}}_{4\varepsilon_{c(\beta)}}(T_{c(\beta)}) \quad {\rm and} \quad \Phi(x)\in  {\bf B}^{\mathcal{F}}_{3\varepsilon_{c(\beta)}}(T_{c(\beta)}).
\end{equation*}

\end{proof}

We say a $p$-cell  $\alpha$ of $I(n+1,\bar k) = I(1,\bar k)\otimes I(n,\bar k)$ is  {\em horizontal} if $\alpha=[y]\otimes\beta$ for some $[y] \in I(1,\bar k)_0$ and  $\beta  \in I(n,\bar k)_p$,  and  {\em vertical} if $\alpha=\gamma \otimes\beta$ for some $\gamma \in I(1,\bar k)_{1}$ and $\beta  \in I(n,\bar k)_{p-1}$.

Let $W_p$ be the set of vertices of $I(n+1,\bar k +pk)$ that belong to the $p$-skeleton of $I(n+1,\bar k)$, i.e., $W_p=\cup_{\alpha\in I(n+1,\bar k)_p} \alpha(pk)_0.$ In particular, $W_{n+1}=I(n+1,\hat k)_0$.

Consider
$$\psi_0:I(1,\bar k)_0\times I(n,\bar k)_0 \rightarrow   \mathcal{Z}_2(M)$$
given by 
$$\psi_0([0],x)=\tilde\phi_{\delta}(x),\quad \psi_0([i\cdot 3^{-\bar k}],x)=\bar \phi(x),$$
where $0<i\leq 3^{\bar k}$. 
We say that a map 
$$\psi_p:W_p \rightarrow \mathcal{Z}_2(M)$$
is a $p$-homotopy if the following conditions hold:
\begin{enumerate}
%\item $\psi_p$ is defined on $I(1,0)_0\times I(n,\bar k)_0$ and on $\alpha(kp)_0$ for every $\alpha$ a $p$-cell of 
%$I(1,0)\times I(n,\bar k)$.
\item $$\psi_p([0],\cdot)=\tilde\phi_{\delta}\circ{\bf n}(\bar k+pk, \bar k), \quad  \psi_p([1],\cdot)=\bar \phi\circ{\bf n}(\bar k+pk, \bar k).$$
\item If $p=1$,   we require that ${\bf f}(\psi_1)\leq c(n)\delta.$ If $p>1$, we ask that  there exists  a $(p-1)$-homotopy
$\psi_{p-1}$ so that $${\bf f}(\psi_p)\leq b(n)({\bf f}(\psi_{p-1})+\delta).$$ 
\item  If $\alpha=\gamma \otimes \beta$ is a  $p$-cell of $I(n+1,\bar k)$,  then
$$\sup_{(y,x)\in \alpha(pk)_0}\{{\bf M}(\psi_p(y,x))\}\leq \sup_{x\in \beta_0}\{{\bf M}(\tilde \phi_{\delta}(x)), {\bf M}(\bar\phi(x))\}+p\delta.
$$
\item If $\alpha=\gamma \otimes \beta$ is a $p$-cell of $I(n+1,\bar k)$, then
 $$\psi_p(\alpha(pk)_0)\subset  {\bf B}_{(3n+6p-2)\varepsilon_{{\bf \bar c}(\beta)}}^{\mathcal{F}}(T_{{\bf \bar c}(\beta)}).$$
\item ${\bf m}(\psi_p,r)<2^{(n+p)(n+2)}(n+p+1)a(n){\delta}.$
\item If $\alpha$ is a horizontal  $p$-cell of $I(n+1,\bar k)$, then 
$$\psi_p =\psi_0 \circ  {\bf n}(\bar k+pk,\bar k)$$
on $\alpha(pk)_0$.
\end{enumerate}

We start by defining a $1$-homotopy $\psi_1$. 
Let $\alpha=\gamma  \otimes \beta$ be a vertical 1-cell of $I(n+1, \bar k)$. By applying Proposition \ref{singleT} to $\psi_0$ on $\alpha$, with $T=T_{{\bf \bar c}(\beta)}$, we get a map $\tilde\psi_{0,\alpha}:\alpha(k_{\bar {\bf c}(\beta)}) \rightarrow {\bf B}^{\mathcal{F}}_{(3n+4)\varepsilon_{\bar {\bf c}(\beta)}}(T_{\bar {\bf c}(\beta)}).$ Note that we can
apply Proposition \ref{singleT} here because of  Lemma \ref{lemm.incl} and  properties b), d), b'), d') above. Since 
$\tilde\psi_{0,\alpha}(x)=\psi_0(x)$ for $x \in \alpha_0$, the map
 ${\psi_1}:W_1 \rightarrow \mathcal{Z}_2(M)$ given by $\psi_1=\tilde\psi_{0,\alpha} \circ {\bf n}(\bar k+k,\bar k+k_{\bar{\bf c}(\beta)})$ on  $\alpha(k)_0$ if $\alpha$ is a vertical 1-cell, and by 
 $\psi_1=\psi_{0} \circ {\bf n}(\bar k+k,\bar k)$ on  $\alpha(k)_0$ if $\alpha$ is a horizontal 1-cell, is well-defined. It follows directly from Proposition \ref{singleT} that $\psi_1$ is a $1$-homotopy.

\subsection{Lemma}\label{iterate3}\textit{Assume $p\leq n$. Given  a $p$-homotopy $\psi_p$, we can find  a $(p+1)$-homotopy 
$\psi_{p+1}$. 
}

\begin{proof}
Let $\alpha=\gamma \otimes \beta $ be a vertical $(p+1)$-cell of $I(n+1,\bar k)$. Hence $\beta \in I(n,\bar k)_p$.  From condition (4) of  the definition of a $p$-homotopy and Lemma \ref{lemm.incl} we have
$$\psi_p(\alpha_0(pk)_0)\subset  {\bf B}_{(3n+6(p+1)-2)\varepsilon_{{\bf \bar c}(\beta)}}^{\mathcal{F}}(T_{{\bf \bar c}(\beta)}).$$
From condition  (3) of the definition of a  $p$-homotopy, and properties b), b'), we also have
$$\sup_{x \alpha_0(pk)_0}\{{\bf M}(\psi_p(y,x))\}\leq L+(n+p)\delta<2L-\delta.$$
Now because of condition (5) we  can apply
Proposition \ref{singleT} to $\psi_p$ in $\alpha$, with  $T=T_{{\bf \bar c}(\beta)}$, $m=p+1$, $j=pk$, to get a map
$$\tilde\psi_{p,\alpha}:\alpha(pk+k_{\bar {\bf c}(\beta)})_0 \rightarrow  {\bf B}_{(3n+6(p+1)-2)\varepsilon_{{\bf \bar c}(\beta)}}^{\mathcal{F}}(T_{{\bf \bar c}(\beta)}).$$ By property (ii) of Proposition \ref{singleT} we get that $$\tilde\psi_{p,\alpha}=\psi_p \circ {\bf n}(\bar k+pk+k_{\bar {\bf c}(\beta)},\bar k+pk)$$
on the $p$-faces of $\alpha$. 

If $\alpha=\gamma \otimes \beta $ is a horizontal $(p+1)$-cell of $I(n+1,\bar k)$, we define
$$\tilde\psi_{p,\alpha}:\alpha(pk+k_{\bar {\bf c}(\beta)})_0 \rightarrow  {\bf B}_{(3n+4)\varepsilon_{{\bf \bar c}(\beta)}}^{\mathcal{F}}(T_{{\bf \bar c}(\beta)})$$
 by $\tilde\psi_{p,\alpha}={\psi_0}\circ{\bf n}(\bar k+ pk+k_{\bar {\bf c}(\beta)},\bar k)$. Since the $p$-faces
 of $\alpha$ are again horizontal cells, we get from condition (6) of the definition of a $p$-homotopy that
 $$
 \tilde\psi_{p,\alpha} = \psi_p \circ {\bf n}(\bar k+pk+k_{\bar {\bf c}(\beta)},\bar k+pk)
 $$
 on the $p$-faces of $\alpha$.

Hence the map $\psi_{p+1}: W_{p+1} \rightarrow \mathcal{Z}_2(M)$  given by
$$\psi_{p+1}=\tilde\psi_{p,\alpha}\circ {\bf n}(\bar k+(p+1)k,\bar k+ pk+k_{\bar {\bf c}(\beta)})$$ on $\alpha((p+1)k)_0$, $\alpha=\gamma \otimes \beta\in I(n,\bar k)_{p+1}$, is well-defined. 

Arguing as in the proofs of Lemmas \ref{iterate1} and Lemma \ref{iterate2}, we can check  that $\psi_{p+1}$ is a $(p+1)$-homotopy.

\end{proof}

Proceeding inductively, we construct  an $(n+1)$-homotopy 
$$\psi=\psi_{n+1}:I(n+1,\hat k)_0 \rightarrow \mathcal{Z}_2(M).$$
 From condition (2) of the definition of a $p$-homotopy it follows that there exists   a universal constant $d(n)$ so that ${\bf f}(\psi)\leq d(n)\delta$. From condition (4) of the definition of a $p$-homotopy we have that
  $$\psi(\alpha((n+1)k)_0)\subset  {\bf B}_{(9n+4)\varepsilon_{{\bf \bar c}(\beta)}}^{\mathcal{F}}(T_{{\bf \bar c}(\beta)})$$
 Thus, we obtain from Lemma \ref{lemm.incl} that
$$\sup\{\mathcal{F}(\psi(y,x)-\Phi(x)): y\in I(1,\hat k)_0, x\in I(n,\hat k)_0\}\leq (9n+7)\varepsilon_N<\delta.$$
Finally, from \eqref{dist6}, and property c), we have that for every $\beta\in I_0(n,\bar k)_{n-1}$ and $z\in \beta$
$$\sup_{x\in \beta_0}\{{\bf M}(\tilde\phi_{\delta}(x)), {\bf M}(\bar\phi(x))\}\leq {\bf M}(\Phi(z))+(n+2)\delta.$$
Therefore, condition (3) of the definition of a $p$-homotopy implies that 
$${\bf M}(\psi(y,x))\leq {\bf M}(\Phi(x))+2(n+2)\delta\quad \mbox{for all } (y,x)\in I(1,\hat k)_0\times I_0(n,\hat k)_0.$$
\end{proof}

We now finish the proof of  Theorem \ref{flattomass}. Let  $$e(n)=\max\{d(n),c(n),2(n+2)\},$$ and let $\{\delta_i\}_{i\in\N}$ be a decreasing sequence of positive numbers converging to zero. Consider 
$$\varphi_i=\phi_{\delta_i/e(n)}:I(n,k_i)_0 \rightarrow \mathcal{Z}_2(M),$$ 
$k_i \rightarrow \infty$, defined as before.  From \eqref{mass.thm.pointwise} we see that
 {for every $y\in I(n,k_i)_0$
$${\bf M}(\phi_i(y))\leq\sup \{{\bf M}(\Phi(x)):\alpha\in I(n,l_i)_n, x,y\in \alpha\}+\delta_i$$}
and this proves Theorem \ref{flattomass}(i).

We can extract a subsequence $\{\phi_i=\varphi_{j_i}\}$   such that $\phi_{i+1}$ is an $(n,\delta_{j_i},k_{j_{i+1}})$-extension of $\Phi$. Proposition \ref{homotopic.sequence}  applied to $\phi_i$ and $\phi_{i+1}$ (replacing $\phi_\delta$ and $\bar \phi$, respectively) gives us a map $\psi_i$ that satisfies Theorem \ref{flattomass}(ii),(iii),  and(iv). 

To prove Theorem \ref{flattomass}(v), we change the construction of the $p$-extension  $\phi_p$ of $\Phi$ so that, whenever $\alpha\in [0]\otimes I(n-1,j)_{p}$, we have $\phi_p=\Phi\circ{\bf n}(j+pk,j)$ on $\alpha(pk)_0$. This is still a $p$-extension because of \eqref{dist7}. Then we redefine  $\phi_{\delta}$ so that,  instead of having $\phi_{\delta}=\Phi\circ{\bf n}(j+nk,j)$ on $\alpha\in [0]\otimes I(n-1,j)_{n-1}$, we have $\phi_{\delta}=\Phi$  on $\alpha(nk)_0$.   The rest of the construction follows exactly as in the previous case.

%%%%%%%%%%%%%%%%%%%%%%%%%%%%%%%%%%%%%%%%%%%%%%%%%%%%%%%%
%%%%%%%%%%%%%%%%%%%%%%%%%%%%%%%%%%%%%%%%%%%%%%%%%%%%%%%%%%%
\section{Interpolation results: Discrete to continuous}\label{discrete.continuous}

In this section we give conditions under which a discrete map is approximated by a continuous map in the mass norm. The main result is important to prove Proposition \ref{pulltight} in Section \ref{proof.pulltight}.

We observe from Corollary 1.14 in \cite{almgren} that there exists  $\delta_0>0$, depending only on $M$, such that for every  $$\psi:I(n,0)_0\rightarrow \mathcal{Z}_2(M)$$  with ${\bf f}(\psi)<\delta_0$, and  $\alpha\in I(n,0)_1$ with  $\partial \alpha = [b]-[a]$,  we can find $Q(\alpha)\in {\bf I}_3(M)$ with  
$$\partial Q(\alpha)=\psi([b])-\psi([a])\quad \mbox{and}\quad{\bf M}(Q(\alpha))={\mathcal F}(\partial Q(\alpha)).$$

The main result of this section is:

\subsection{Theorem}\label{continuous.approximation2} \textit{There exists  $C_0>0$, depending only on $M$ and $n$, such that for every map 
$$\psi:I(n,0)_0\rightarrow \mathcal{Z}_2(M)$$  with ${\bf f}(\psi)<\delta_0$, we can find a continuous map in the mass norm
$$\Psi:I^n\rightarrow \mathcal{Z}_2(M;{\bf M})$$
such that
\begin{itemize}
\item [(i)] $\Psi(x)=\psi(x)\mbox{ for all }x\in I(n,0)_0$;
\item [(ii)] for every $\alpha\in I(n,0)_p$, $\Psi_{|\alpha}$ depends only on the values assumed by $\psi$ on the vertices of $\alpha$;
\item [(iii)] 
%\begin{multline*}
$$
\sup\{{\bf M}(\Psi(x)-\Psi(y)): x,y\in I^n\}\\
\leq C_0\sup_{ \alpha\in I(n,0)_1}\{{\bf M}(\partial Q(\alpha))\}.
%\end{multline*}
$$
\end{itemize}
}

An immediate consequence  is:

\subsection{Theorem}\label{continuous.approximation}\textit{For every map 
$$\psi:I(n,k)_0\rightarrow \mathcal{Z}_2(M)$$  with ${\bf f}(\psi)<\delta_0$, we can find a continuous map in the mass norm
$$\Psi:I^n\rightarrow \mathcal{Z}_2(M;{\bf M})$$
such that
\begin{itemize}
\item $\Psi(x)=\psi(x)\mbox{ for all }x\in I(n,k)_0$;
\item for every $\alpha\in I(n,k)_p$
$$\sup\{{\bf M}(\Psi(x)-\Psi(y)): x,y\in \alpha\}\leq C_0{\bf f}(\psi).$$
\end{itemize}
}

\begin{proof}
Let $\alpha$ be an $n$-cell of $I(n,k)$. By identifying $\alpha$ with $I^n$ and  applying Theorem \ref{continuous.approximation2} to  $\psi_{|\alpha_0}$, we get a continuous map $\Psi_\alpha:\alpha
\rightarrow  \mathcal{Z}_2(M;{\bf M})$ satisfying
$$
\sup\{{\bf M}(\Psi_\alpha(x)-\Psi_\alpha(y)): x,y\in \alpha\}\\
\leq C_0{\bf f}(\psi).
$$  It follows from Theorem \ref{continuous.approximation2} (ii) that these continuous maps  obtained from different $n$-cells coincide along common  faces,  thus giving us a well-defined map $\Psi:I^n\rightarrow \mathcal{Z}_2(M;{\bf M})$. 
\end{proof}

\subsection{Proof of Theorem \ref{continuous.approximation2}}
We note that a similar result  was proven by Almgren in Theorem 6.6 of \cite{almgren}. In our case the situation is simpler because we are dealing  with codimension one currents (2-currents in a three-manifold). The work of Almgren gives us  a map $\Psi$ that is continuous in the flat metric and satisfies  (i), (ii), and
\begin{itemize}
\item [(iii')] $$
\sup\{{\mathcal F}(\Psi(x)-\Psi(y)): x,y\in I^n\}\\
\leq C_0\sup_{ \alpha\in I(n,0)_1}\{{\bf M}(Q(\alpha))\}.
$$
\end{itemize}
In Theorem 4.6 of \cite{pitts}, Pitts explains how to adapt the methods of \cite{almgren} to  make them work in
the context of  maps that are continuous in the mass norm. This involves the construction of the continuous
map $\Psi:I^n\rightarrow \mathcal{Z}_2(M;{\bf M})$.  It follows from the proof of \cite[Theorem 6.6]{almgren}, with no modification whatsoever,  that properties (i) and (ii) of Theorem  \ref{continuous.approximation2} are satisfied. Hence the  statement of Theorem \ref{continuous.approximation2} which requires justification is the third one. We will briefly sketch  the proof of Theorem 6.6 of \cite{almgren} and show that Theorem  \ref{continuous.approximation2} (iii) indeed holds.

Let $\Delta$ be a differentiable triangulation of $M$. Hence if $s \in \Delta$ then the faces of $s$ also belong to  $\Delta$. Choose a linear order $\prec$ on $\Delta$ such that  $s'\prec s$ if  {${\rm dim(s')}<{\rm dim}(s)$. Given $s,s'\in\Delta$, we use  $s'\subset s$ if  $s'$ is a face of $s$}. Let  {$U(s)=\cup_{s\subset s'}s'$.} In what follows we will denote by $C$ varying   constants that depend only on $\Delta$ and $n$.

The first ingredient in the construction of $\Psi$ is to consider,  {for every $s\in \Delta$},  a deformation map
$${ {\mathcal D(s)}}: I\times {\bf I}_2(U(s);{\bf M})\rightarrow {\bf I}_2(U(s);{\bf M})$$
such that
\begin{itemize}
\item ${\mathcal D(s)}$ is continuous  in the mass norm;
\item ${\mathcal D}(s,0,T)=T$ and ${\mathcal D}(s,1,T)=0$ for every $s\in \Delta$, $T \in {\bf I}_2(U(s);{\bf M})$;
\item for all $s\in \Delta, t\in I$, and $T\in  {\bf I}_2(U(s))$, we have
\begin{equation}
\label{mass.dilation}{\bf M}({\mathcal D}(s,t,T))\leq C{\bf M}(T).
\end{equation}
Here ${\mathcal D}(s,t,T)={\mathcal D}(s)(t,T)$.
\end{itemize}
The construction of such maps uses the deformation maps of  \cite[Theorem 4.5]{pitts}. In the context of flat metrics, this construction was carried out  in  \cite[Section 5]{almgren}. 

The second ingredient is to consider the cutting functions, which we describe now. Let  $\Lambda\subset {\bf I}_3(M)$ be a finite set with $q$ elements.  Almgren \cite[Section 5]{almgren} associates to every $s\in\Delta$ a   {neighborhood} $L(s)$ of $s$ and constructs a function
$$C_{\Lambda}:\Delta\times \Lambda\rightarrow  {\bf I}_3(M)$$
satisfying, according to Definition 5.4, Theorem 5.8 and Lemma 5.9 of \cite{almgren},
\begin{itemize}
\item \begin{equation}\label{cut.first.property}
C_{\Lambda}(s,T)=\left(T-\sum_{s'\prec s}C_{\Lambda}(s',T)\right)\cap L(s);
\end{equation}
\item \begin{multline}\label{mass.cut.inductive}
{\bf M}\left(\partial C_{\Lambda}(s,T)-\partial\left(T-\sum_{s'\prec s}C_{\Lambda}(s',T)\right)\cap L(s)\right)\\
\leq C\cdot q \cdot  {\bf M}\left(T-\sum_{s'\prec s}C_{\Lambda}(s',T)\right);
\end{multline}
\item  {$\mathrm{support}\left(C_{\Lambda}(s,T)\right)\subset U(s)$ for all $(s,T)\in \Delta\times \Lambda$.}
\end{itemize}
From \eqref{cut.first.property} we see that 
\begin{equation*}\label{mass.cut.3chain}
{\bf M}(C_{\Lambda}(s,T))\leq C{\bf M}(T)\quad\mbox{ for every }(s,T)\in \Delta\times \Lambda.
\end{equation*}
This inequality and \eqref{mass.cut.inductive} imply that
$$
{\bf M}(\partial C_{\Lambda}(s,T))\leq  Cq{\bf M}(T)+{\bf M}(\partial T)+\sum_{s'\prec s}{\bf M}(\partial C_{\Lambda}(s',T)).
$$
Thus we conclude that
\begin{equation}\label{cut.mass}
{\bf M}(\partial C_{\Lambda}(s,T))\leq  C\cdot q \cdot ({\bf M}(T)+{\bf M}(\partial T))\quad\mbox{ for every }(s,T)\in \Delta\times \Lambda.
\end{equation}

Having defined the basic ingredients we recall Almgren's construction of the map $\Psi$. 
For every $p$-cell $\alpha$ of $I(n,0)$,  we consider the  continuous function
$$h_{\alpha}:I^p\rightarrow  \mathcal{Z}_2(M;{\bf M})$$
given by  $h_{\alpha}(0)=\psi(\alpha)$ if $p=0$, and by the following formula if $p>0$  \cite[Interpolation Formula 6.3]{almgren}: 
\begin{multline}\label{halpha}
h_{\alpha}(x_1,\ldots,x_p)=\sum_{\gamma\in \Gamma_\alpha}\mathrm{sign}\,(\gamma)\\ \cdot \sum_{s_1,\ldots,s_p\in\Delta}
{\mathcal D}(s_1,x_1)\circ\ldots \circ{\mathcal D}(s_p,x_p)\circ \partial \circ C_{\Lambda(\gamma_p)}(s_p)\circ\ldots \circ C_{\Lambda(\gamma_1)}(s_1)\left(Q(\gamma_1) \right)
\end{multline}
where
\begin{itemize}
\item $\Gamma_{\alpha}$ denotes the set of all sequences $\{\gamma_i\}_{i=1}^p$ such that $\gamma_p=\alpha$ and such that, for each $1\leq i \leq p-1,$   {$\gamma_i$ is a 
$\left({\rm dim}(\gamma_{i+1})-1\right)$-face of  $\gamma_{i+1}$. }
\item $\mathrm{sign}\,(\gamma)$ is equal to 1 or $-1$, according to   \cite[Definition 6.2]{almgren}.
\item the  finite  sets $\Lambda(\beta)$ are defined inductively in the following way: if  $\beta\in I(n,0)_1$, we have   $\Lambda(\beta)=\{Q(\beta)\}$; if  $\beta\in I(n,0)_j$ with $j>1$, we have
\begin{multline*}
\Lambda(\beta)=\{C_{\Lambda(\beta_{j-1})}(s_{j-1})\circ \ldots\circ C_{\Lambda(\beta_1)}(s_1)(Q(\beta_1))\\
:s_k\in \Delta\mbox{ and } \beta_k\mbox{ is a \ } k\mbox{-cell of }\beta\mbox{ for every }k=1,\ldots,j-1\}. 
\end{multline*}
 \end{itemize}
 
 Having fixed the triangulation $\Delta$, the deformation maps (which depend only on $\Delta$), and the cutting function $C_{\Lambda(\beta)}$ for each cell $\beta$, it is clear  that $h_{\alpha}$ is continuous in the mass norm and that it depends only on  the values assumed by  $\psi$ on the vertices of $\alpha$.   {Almgren describes in  \cite[Section 6.5]{almgren} an inductive procedure to construct $\Psi$ using  the various maps $h_{\alpha}$ described above.}

\subsection{Lemma}\label{mass.bound.almgren} \textit{For every $x\in I^p$ and $\alpha\in I(n,0)_p$, with $p\geq 1$, we have
$${\bf M}(h_{\alpha}(x))\leq C\sup\{{\bf M}(\partial Q(\beta)):\beta\in I(n,0)_1, \beta\subset \alpha\}.$$
}
\begin{proof}
The cardinality of every finite set $\Lambda(\beta)$ is bounded above by a constant depending only on $\Delta$ and $n$. Hence we obtain from \eqref{cut.mass} that 
\begin{multline}\label{cut.mass.bound}
{\bf M}\left (\partial \circ C_{\Lambda(\gamma_p)}(s_p)\circ\ldots \circ C_{\Lambda(\gamma_1)}(s_1)(Q(\gamma_1)) \right)\\\leq C\left({\bf M}(Q(\gamma_1))+{\bf M}(\partial Q(\gamma_1))\right)\leq C{\bf M}(\partial Q(\gamma_1))
\end{multline}
for every
$\{\gamma_i\}_{i=1}^p\in \Gamma_{\alpha}$, 
where the last inequality comes from the fact that $${\bf M}(Q(\gamma_1))={\mathcal F}(\partial Q(\gamma_1))\leq{\bf M}(\partial Q(\gamma_1)).$$
The number of elements of $\Gamma_{\alpha}$ is bounded above by a constant depending only on $n$, hence the desired result follows from the expression \eqref{halpha} for $h_{\alpha}$,  combined with \eqref{mass.dilation} and \eqref{cut.mass.bound}.
\end{proof}
 {Using Lemma \ref{mass.bound.almgren}, the proof of \cite[Theorem 6.6 (2) (b)]{almgren} applies with no modifications to conclude Theorem \ref{continuous.approximation2} (iii).}
%Using   \cite[Extension Criterion (4.5)]{almgren}, Almgren constructs the continuous extension $\Psi$ by defining it as a linear combination  of the functions $\{h_{\alpha}\}_{\alpha\in I(n,0)}$, where  each $h_{\alpha}$ is extended from a map defined on $I^p$ to a map defined on $I^n$ via projections. The number of terms in the linear combination depends only on $n$. Because  $h_{\alpha}$ is a constant function if $\alpha$ is a $0$-cell, we have that $\Psi(x)-\Psi(y)$ is a linear combination of functions $h_{\alpha}$ with $\dim(\alpha)\geq 1$. Thus Theorem \ref{continuous.approximation2} (iii) follows  from Lemma \ref{mass.bound.almgren}.

%%%%%%%%%%%%%%%%%%%%%%%%%%%%%%%%%%%%%%%%%%%%%%%%%%%%%%%%
%%%%%%%%%%%%%%%%%%%%%%%%%%%%%%%%%%%%%%%%%%%%%%%%%%%%%%%%%%%

\section{Pull-tight}\label{proof.pulltight}

Assume we have a continuous map in the flat topology
$$\Phi:I^{n}\rightarrow   \mathcal{Z}_2(M)$$
 {which satisfies the following hypotheses:
\begin{itemize}
\item[($B_0$)] $\Phi_{|I_0^n}$ is continuous in the ${\bf F}$-metric;
\item[($B_1$)] 
$\Phi(I^{n-1}\times\{0\})=\Phi(I^{n-1}\times\{1\})=0.$
\end{itemize}}
We denote by   $|\Phi|:I^n \rightarrow \mathcal{V}_2(M)$ the map given by $|\Phi|(x)=|\Phi(x)|$ for $x\in I^n$.

Consider $\Pi \in  \pi_n^{\#}(\mathcal{Z}_2(M;{\bf M}),\Phi_{|I_0^n})$.  

\subsection*{Proposition \ref{pulltight}}\textit{
There exists  a critical sequence $S^* \in \Pi$. For each  critical sequence $S^*$, there exists a critical sequence $S\in\Pi$ such that
\begin{itemize}
\item ${\bf C}(S)\subset {\bf C}(S^*)$;
\item every $\Sigma\in {\bf C}(S)$ is either a stationary varifold or belongs to $|\Phi|(I^n_0)$.
\end{itemize}
}
\begin{proof}

We start with a basic lemma that proves the existence of critical sequences. This is just like \cite[Subsection 4.1, Proposition 4]{pitts}.

\subsection{Lemma}\label{critical.pitts} \textit{There exists  a critical sequence $S^*\in \Pi$.
}
\begin{proof}
We choose $S^j=\{\phi^j_i\}_{i\in\N}\in \Pi$ such that ${\bf L}(S^j)\leq {\bf L}(\Pi)+1/j$, and pick an increasing sequence $\{n_j\}_{j\in\N}$ so that we have, for all $i\geq n_j,$ 
\begin{itemize}
\item$\max\{{\bf M}(\phi^j_i(x)):x\in\dmn(\phi^j_i)\}\leq {\bf L}(S^j)+1/j;$
\item $\phi^j_i\mbox{ is $n$-homotopic to }\phi^j_{i+1}\mbox{ with fineness  }1/j;$
\item {$\phi^1_i$ and $\phi^j_i\mbox{ are  $n$-homotopic to }\phi^{j+1}_{i}\mbox{ with fineness  }1/j$.}% for any $1\leq k \leq j$.
\end{itemize}
Let $\phi^*_i$ be given by  $\phi_i^*=\phi^1_i$ if $i\leq n_2-1$, and $\phi_i^*=\phi^j_i$ if $n_j\leq i\leq n_{j+1}-1$. Then $S^*=\{\phi^*_i\}\in \Pi$ and ${\bf L}(S^*)={\bf L}(\Pi)$.
\end{proof}

Given a critical sequence $S^*\in \Pi$, we apply a ``pull-tight'' procedure to $S^*$ to find another critical sequence $S\in\Pi$ such that all elements of ${\bf C}(S)$ are either stationary varifolds or belong to $|\Phi|(I^n_0)$. We essentially follow the method of \cite[Theorem 4.3]{pitts}.

Suppose $S^*=\{\phi_i^*\}_{i\in \N}$,  and set
$$c=\sup \{{\bf M}(\phi^*_i(x)):i\in\N, x\in \dmn(\phi_i^*)\}.$$
We define the following compact sets of $\mathcal{V}_2(M)$:
\begin{align*}
A&=\{V\in \mathcal{V}_2(M):||V||(M)\leq c\},\\
B&=|\Phi|(I^n_0)\subset A,\\
A_0&=B\cup\{V\in A: V\mbox{ is stationary in M}\},\\
A_1&= \{V\in A: {\bf F}(V,A_0)\geq 2^{-1}\},\\
A_i&= \{V\in A: 2^{-i}\leq {\bf F}(V,A_0)\leq 2^{-i+1}\},\quad i\in \{2,3,\ldots\}.
\end{align*}

For every $V\in A_i$, $i \geq 1$, we choose a vector field $X_V\in \mathcal{X}(M)$  with $|X_V|_{C^1}\leq 1$ and such that
$$\delta V(X_V)\leq \frac{2}{3}\inf\{\delta V(Y):Y\in \mathcal{X}(M)\mbox{ with }|Y|_{C^1}\leq 1\}<0.$$
The map $S \in \mathcal{V}_2(M) \mapsto \delta S(X_V)$ is continuous, hence we can find for every $V\in A_i$, $i\geq 1$, a radius $0<r_V<2^{-i}$ so that  we have 
$$\delta S(X_V)\leq \frac{1}{2}\inf\{\delta S(Y):Y\in \mathcal{X}(M)\mbox{ with }|Y|_{C^1}\leq 1\}<0$$
for every $S\in {\bf B}^{\bf F}_{r_V}(V)$.
The compactness of $A_i$ implies that the open cover ${\bf B}^{\bf F}_{r_V}(V)$  admits a finite subcover. Thus
we can find $q_i \in\N$ and
\begin{itemize}
 \item a set of radii $\{r_{ij}\}_{j=1}^{q_i}$, $r_{ij}< 2^{-i}$;
 \item  a set of varifolds $\{V_{ij}\}_{j=1}^{q_i}\subset A_i$;
 \item a set of vector fields $\{X_{ij}\}_{j=1}^{q_i}\subset \mathcal{X}(M)$ with $|X_{ij}|_{C^1} \leq 1$;
 \item a set of balls $U_{ij}={\bf B}^{\bf F}_{r_{ij}}(V_{ij})\cap A$, $j=1,\ldots,q_i$, with $A_i\subset \bigcup_{j=1}^{q_i}U_{ij};$
 \item a set of positive real numbers $\{\varepsilon_{ij}\}_{j=1}^{q_i}$ such that 
  $$\delta S(X_{ij})\leq -\varepsilon_{ij}<0\mbox{ for all }S\in U_{ij}, j=1,\ldots,q_i.$$
\end{itemize}

The condition  $r_{ij}<2^{-i}$ implies that $\{U_{ij}\}_{i\in \N,1\leq j\leq q_i}$ is a locally finite covering of  $A\setminus A_0$. Therefore   we can choose  a partition of unity $\{\phi_{ij}\}_{i\in \N,1\leq j\leq q_i}$ of $A\setminus A_0$ with $\mathrm{support}(\phi_{ij})\subset U_{ij}$. 

We define
$$X:A\rightarrow \mathcal{X}(M),$$
continuous in the ${\bf F}$-metric, by
\begin{align*}
X(V)&=0 \quad\mbox{if }V\in A_0,\\
X(V)&= {{\bf F}(V,A_0)} \sum_{i\in\N, 1\leq j\leq q_i}\phi_{ij}(V)X_{ij}\quad\mbox{if }V\in A\setminus A_0.
\end{align*}
It follows that  $$\delta V(X(V))=0\mbox{ if }V\in A_0\mbox{ and }\delta V(X(V))<0\mbox{ if }V\in A\setminus A_0.$$
This implies we can find a continuous function $$h:A\rightarrow [0,1]$$ such that 
\begin{itemize}
\item $h=0$ on $A_0$ and $h(V)>0$ if $V\in A\setminus A_0$, 
\item and $||{f(s,V)}_{\#}(V)||(M)< ||{f(t,V)}_{\#}(V)||(M)$ if $0\leq t<s\leq h(V)$,
\end{itemize}
where $f(t,V)$ denotes the  one-parameter group of diffeomorphisms generated by $X(V)$.

Now let 
$$H:[0,1]\times(\mathcal{Z}_2(M;{\bf F})\cap \{S:{\bf M}(S)\leq c\})\rightarrow \mathcal{Z}_2(M;{\bf F})\cap \{S:{\bf M}(S)\leq c\}$$
be given by
\begin{align*}
H(t,T)&= {f(t,|T|)}_{\#}(T)\quad  \mbox{if } 0\leq t\leq h(|T|),\\
H(t,T)&= {f(h(|T|),|T|)}_{\#}(T)\quad \mbox{if } h(|T|)\leq t\leq 1.
\end{align*} 
The key properties of $H$ are
\begin{itemize}
\item[(i)]$H$ is continuous in the product topology;
\item[(ii)] $H(t,T)=T$ for all $0\leq t\leq 1$ if $|T|\in A_0$;
\item[(iii)] $||H(1,T)||(M)<||T||(M)$ unless $T\in A_0$;
\item[(iv)] for every $\varepsilon>0$, there exists $\delta>0$ so that {for all $x\in I^n_0$ and all $0\leq t\leq 1$}
$${\bf F}(T,\Phi(x))<\delta\Rightarrow {\bf F}(H(t,T),\Phi(x))<\varepsilon.$$
This property is a consequence of the first two, since $B=|\Phi|(I_0^n) \subset A_0$ { and } $\Phi_{| I^n_0}$ is continuous in the ${\bf F}$-metric.
\end{itemize}

We now proceed to the construction of $S=\{\phi_i\}_{i\in\N}\in \Pi$ with ${\bf C}(S)\subset A_0\cap {\bf C}(S^*)$. 
We would like to put $\phi_i=H(1,\phi_i^*)$. Since the map 
$$G:\mathcal{Z}_2(M)\rightarrow \mathcal{Z}_2(M),\quad G(T)=F_{\#}(T),$$
where $F \in {\rm Diff}(M)$ is fixed, 
is continuous in the ${\bf F}$-metric but not  in the mass norm, the fineness of $\phi_i$ could be large even when ${\bf f}(\phi^*_i)$ is small. Thus we need to  interpolate $H(1,\phi_i^*)$ one more time, as in Theorem \ref{flattomass}.  When doing this, it is important to check that the values assumed by  $\phi_i$ stay close in the ${\bf F}$-metric to those assumed by $H(1,\phi_i^*)$. 

This minor issue was overlooked by Pitts \cite[page 153]{pitts}. We overcome this difficulty using the Interpolation Theorem \ref{continuous.approximation} of Section \ref{discrete.continuous}. This requires a bit of extra work that we do now.

Denote the domain of $\phi_i^*$ by $I(n,k_i)_0$, and let  $\delta_i={\bf f}(\phi_i^*)$. Apply Theorem \ref{continuous.approximation} to obtain a continuous map  in the mass norm
$$
\bar \Omega_i:I^n\rightarrow  \mathcal{Z}_2(M;{\bf M}),
$$
such that for all $x\in  I(n,k_i)_0$ and $\alpha\in I(n,k_i)_n$ we have
\begin{equation}\label{massomega} 
 \bar\Omega_i(x)=\phi_i^*(x),\quad\mbox{and}\quad\sup_{y,z\in\alpha}\{{\bf M}(\bar \Omega_i(z)-\bar \Omega_i(y))\}\leq C_0\delta_i.
\end{equation}
We claim that
\begin{equation}\label{baromega.phi}
\lim_{i\to\infty}\sup\{{\bf F}(\bar \Omega_i(x), \Phi(x)):x\in I^n_0\}=0.
\end{equation}
Indeed, from Lemma \ref{homotopy.sequence.boundary} we have that 
$$
\lim_{i\to\infty}\sup\{{\bf F}(\phi^*_i(x), \Phi(x)):x\in I_0(n,k_i)_0\}=0.
$$
The claim then follows from \eqref{massomega}.

Consider the continuous map in the ${\bf F}$-metric
$$\Omega_i:I \times I^n\rightarrow  \mathcal{Z}_2(M;{\bf F}),\quad  \Omega_i(t,x)=H(t,\bar \Omega_i(x)).$$
From property (iii) of $H$ we have
\begin{equation}\label{mass.omega}
\max \{{\bf M}(\Omega_i(t,x)):(t,x)\in I\times I^n\}\leq 
\max_{x\in I^n}\{{\bf M}(\bar \Omega_i(x))\}.
\end{equation}
 From property (iv) of $H$ and \eqref{baromega.phi} it follows that
\begin{equation}\label{omega.phi}
\lim_{i\to\infty}\sup\{{\bf F}(\Omega_i(t,x), \Phi(x)):(t,x)\in I\times I^n_0\}=0.
\end{equation}

\subsection{Lemma}\textit{ For every $i\in\N$, $\lim_{r\to 0}{\bf m}(\Omega_i,r)=0.$
}

\begin{proof}
Let $\delta>0$.  Note that ${\mathcal C}=\Omega_i(I\times I^n)$ is a compact subset of $ \mathcal{Z}_2(M;{\bf F}).$  For every $p\in M$ and $T\in {\mathcal C}$, and since $T$ is an integral current, we can choose  $r=r(p,T)>0$ so that 
$$||S||(B_r(p))<\delta\quad\mbox{for all }S\in {\bf B}^{\bf F}_{r}(T).$$
By compactness, we can select a finite covering $\{B_{r_k}(p_k)\times {\bf B}^{\bf F}_{r_k}(T_k)\}_{k=1}^N$ of $M\times {\mathcal C}$, where $r_k=r(p_k,T_k)/2$.

If $\bar r=\min\{r_k\}_{k=1}^N$, then
$$||T||(B_{\bar r}(p))<\delta\mbox{ for all }(p, T)\in M\times{\mathcal C}.$$
\end{proof}

We can now apply  Theorem \ref{flattomass} to $\Omega_i$ and obtain
$$\bar \phi_{ij}: I(1,s_{ij})_0 \times I(n,s_{ij})_0 \rightarrow   \mathcal{Z}_2(M)$$
such that 
\begin{itemize}
\item[(a)] 
$$\sup\{{\bf M}(\bar \phi_{ij}(t,x)):(t,x)\in I(n+1,s_{ij})_0\}\leq \max_{x\in I^n}\{{\bf M}(\bar \Omega_i(x))\}+\frac{1}{j};
 $$
 \item[(b)] ${\bf f}(\bar \phi_{ij})<\frac{1}{j};$
 \item[(c)]
 $$\sup\{{\mathcal{F}}(\bar \phi_{ij}(t,x)-\Omega_i(t,x))\,:\, (t,x)\in I(n+1,s_{ij})_0\}\leq \frac{1}{j};
$$
\item[(d)]
$$
{\bf M}(\bar \phi_{ij}(t,x))\leq {\bf M}(\Omega_i(t,x))+\frac 1j\quad\mbox{for all $(t,x)\in I_0(n+1,s_{ij})_0$};
$$
\item[(e)]$\bar \phi_{ij}([0],x)=\Omega_i(0,x)=\bar\Omega_i(x)$ for all  $x\in  I(n,s_{ij})_0.$
\end{itemize}
From Lemma \ref{flat+mass=f}, properties (c), and (d), we get
\begin{equation*}
\lim_{j\to\infty}\sup\{{\bf F}(\bar \phi_{ij}(t,x),\Omega_i(t,x)):(t,x)\in I_0(n+1,s_{ij})_0\}=0.
\end{equation*}
Hence, using \eqref{omega.phi} and a diagonal sequence argument, we can find $\{\bar \phi_i=\bar \phi_{ij(i)}\}$ such that 
\begin{equation}\label{phi_ij.omega_ij}
\lim_{i\to\infty}\sup\{{\bf F}(\bar \phi_{i}(t,x),\Omega_i(t,x)):(t,x)\in I_0(n+1,s_{ij})_0\}=0
\end{equation}
and
\begin{equation}\label{phi_ij.phi}
\lim_{i\to\infty}\sup\{{\bf F}(\bar\phi_{i}(t,x),\Phi(x)): t  \in I(1,s_{ij})_0, x\in I_0(n,s_{ij})_0\}=0.
\end{equation}
We define $\hat \phi_i: I(1,s_{ij(i)})_0\times I(n,s_{ij(i)})_0\rightarrow   \mathcal{Z}_2(M)$ to be equal to zero on $$ I(1,s_{ij(i)})_0 \times (T(n,s_{ij(i)})_0\cup B(n,s_{ij(i)})_0),$$ and equal to $\bar \phi_i$ otherwise. Since  ${\bf f}(\hat \phi_i)$ tends to zero,  we obtain from  \eqref{phi_ij.phi} that $\phi_i=\hat \phi_i([1],\cdot)$ is $n$-homotopic to $\hat \phi_i([0],\cdot)$ in $(\mathcal{Z}_2(M;{\bf M}),\Phi_{|I_0^n})$ with fineness tending to zero.

On the other hand, it follows from  \eqref{massomega} and property e) that $\hat \phi_i([0],\cdot)$ is $n$-homotopic to $\phi^*_i$ in $(\mathcal{Z}_2(M;{\bf M}),\Phi_{|I_0^n})$ with fineness tending to zero. Hence  $S=\{\phi_i\}_{i\in \N}\in \Pi$. From property (a) and \eqref{massomega} we obtain  that $S$ is a critical sequence, i.e., ${\bf L}(S)={\bf L}(\Pi)$.

We are left to show that
 $${\bf C}(S)\subset A_0\cap {\bf C}(S^*).$$
Given $V\in {\bf C}(S)$, there exists  a sequence $\{|\phi_{k_i}(x_i)|\}_{i\in\N}$, $k_i \rightarrow \infty$, that converges to $V$ in the sense of varifolds. It follows from \eqref{phi_ij.omega_ij}   that $|\Omega_{k_i}(1,x_i)|$ also tends to $V$ as varifolds. Moreover, from \eqref{massomega} we see that a subsequence of $|\bar\Omega_{k_i}(x_i)|$ converges as varifolds to an element $W$ of ${\bf K}(S^*)$. Since the map $H$ is continuous in the ${\bf F}$-metric, we have
 $$V=\lim_{i\to\infty}|\Omega_{k_i}(1,x_i)|=\lim_{i\to\infty} |H(1,\bar\Omega_{k_i}(x_i))|=f(h(W),W)_{\#}W.$$ 
 If $V\notin A_0$ then, from property (iii) of $H$, we get
$${\bf L}(\Pi)=||V||(M)=||f(h(W),W)_{\#}W||(M)<||W||(M)\leq {\bf L}(\Pi).$$
This is a contradiction, hence $V\in A_0$. Property (ii) of $H$ implies  that $V=W\in {\bf C}(S^*)$.
  
\end{proof}

%%%%%%%%%%%%%%%%%%%%%%%%%%%%%%%%%%%%%%%%%%%%%%%%%%%%%%%%
%%%%%%%%%%%%%%%%%%%%%%%%%%%%%%%%%%%%%%%%%%%%%%%%%%%%%%%%%%%
\appendix

\section{}\label{F1.section}

Let 
\begin{multline*}
\mathcal{F}_1=\{ S \subset S^3: S \mbox{ is an embedded closed minimal surface of}  \\
\mbox{genus } g(S)\geq 1\}.
\end{multline*}
The goal of this appendix is to prove
\subsection{Theorem}\label{existence.minimizer} \textit{There exists $\Sigma$ in  $\mathcal{F}_1$ such that
$$
{\rm area}(\Sigma) = \inf_{S \in \mathcal{F}_1} {\rm area}(S).
$$}

The proof is largely standard and the method well known among the experts (see, for instance,  \cite[Theorem 2.1]{kuwert-li-schatzle}).
\begin{proof}
Let $\Sigma^{i}\in \mathcal{F}_1$ be a minimizing sequence, i.e., such that $$\lim_{i\to\infty}{\rm area}(\Sigma^{i})=\inf_{S \in \mathcal{F}_1} {\rm area}(S).$$ 
 {Allard Compactness Theorem \cite[Theorem 42.7]{simon}} implies that we can extract a subsequence converging in ${\mathcal V}_2(S^3)$ to an integral stationary varifold $\Sigma.$
Since the Clifford torus has area $2\pi^2$, we have $$||\Sigma||(S^3)=\lim_{i\to\infty}{\rm area}(\Sigma^{i})  \leq 2\pi^2<8\pi(1-\delta)$$ for some $\delta>0$.

\subsection{Lemma}\label{density.smaller.two}\textit{ There is $r_0$ so that
$$\frac{||\Sigma||(B_r(p))}{\pi r^2}\leq 2-\delta\mbox{ for all }r\leq r_0, p\in S^3.$$
}
\begin{proof}
Suppose not. Then we could find  sequences $\{q_i\}_{i\in \N}$, $\{r_i\}_{i\in \N}$ tending to $q\in \Sigma$ and zero, respectively, such that
$$\lim_{i\to\infty} \frac{||\Sigma||(B_{r_i}(q_i))}{\pi {r_i}^2}\geq 2-\delta.$$
The monotonicity formula \cite[Theorem 17.6]{simon} on  a general ambient manifold implies that
\begin{equation*}
\lim_{r\to 0}\frac{||\Sigma||(B_r(q))}{\pi r^2}\geq 2-\delta.
\end{equation*}

Consider the cone $C$ in $\R^4$ defined by
$$C=\mu_{\#}(\Sigma\times\R),\quad\mbox{where }\mu:S^3\times\R\rightarrow \R^4 \,\,\mu(p,r)=rp.$$
Because $\Sigma$ is a stationary varifold in $S^3$, $C$ is a stationary integral varifold in $\R^4$ where,  {denoting by $\omega_3$ the volume of a $3$-ball}, we have
\begin{itemize}
\item $$\frac{||C||(B^4_r(0))}{\omega_3 r^3}=\frac{||\Sigma||(S^3)}{4\pi}\leq 2(1-\delta)\quad\mbox{for all }r>0,$$
\item $$\lim_{r\to 0}\frac{||C||(B^4_r(q))}{\omega_3 r^3}=\lim_{r\to 0}\frac{||\Sigma||(B_r(q))}{\pi r^2}\geq 2-\delta.$$
\end{itemize}
Combining these two facts  with the monotonicity formula we  obtain a contradiction because
$$
2-\delta\leq \lim_{r\to 0}\frac{||C||(B^4_r(q))}{\omega_3 r^3}\leq \lim_{r\to \infty}\frac{||C||(B^4_r(q))}{\omega_3 r^3}=\lim_{r\to \infty}\frac{||C||(B^4_r(0))}{\omega_3 r^3}=2-2\delta.
$$
\end{proof}

\subsection{Lemma}\textit{$\Sigma$ is smooth with multiplicity one.}

\begin{proof}
	From Allard Regularity Theorem \cite[Theorem 24.2]{simon} it suffices to see that
	\begin{equation}\label{densities.one}
	\lim_{r\to 0}\frac{||\Sigma||(B_r(p))}{\pi r^2}=1\quad\mbox{for all }p\in \Sigma.
	\end{equation}
Choose $p\in \Sigma$ and consider, for  every $\lambda\in\R$, the dilation map $\mu^{\lambda}(x)=\lambda x$ defined in $\R^4$. Set
$$\Sigma_j=\mu_{{\#}}^{j}(\Sigma-p), \quad j\in\N $$
which  is a  {varifold in $\R^4$ with generalized mean curvature tending to zero uniformly.} From Allard Compactness  Theorem, we have that a subsequence converges to a stationary varifold $V\subset p^{\bot}$, where $p^{\bot}\subset \R^4$ denotes the hyperplane orthogonal to $p$. Moreover, we must have from scale invariance and Lemma \ref{density.smaller.two}  that, for all $s>0$,
\begin{equation}\label{densities.blow.ups}
\frac{||V||(B^4_s(0))}{\pi s^2}=\lim_{j\to \infty} \frac{||\Sigma^j||(B^4_s(0))}{\pi s^2}=\lim_{r\to 0} \frac{||\Sigma||(B_r(p))}{\pi r^2}\leq 2-\delta,
\end{equation}
and so the  monotonicity formula implies that $V$ is a stationary cone in $p^{\bot}$.

From \cite{allard-almgren} we know that $V$ is a cone over a stationary $1$-varifold $\gamma\subset S^2$, which is a network  consisting of  geodesic segments meeting at triple junctions. If we show that $\gamma$ has no triple junctions then $V$ must be a plane, which has multiplicity one from \eqref{densities.blow.ups}, and so 
 \eqref{densities.one} follows at once. 

Suppose $x_0$ is a triple junction of $\gamma$ and consider the sequence of integral stationary varifolds
$$V_k=\mu_{{\#}}^{k}(V-x_0),\quad k\in \N.$$
From \cite[Theorem A.4]{simon} we know that, after passing to a subsequence, $V_k$ converges to a stationary varifold $U$, which consists of three half-planes  $\{P_1,P_2,P_3\}$ of $p^{\bot}$ meeting along a common line $L$. Note that these half-planes must have multiplicity one from \eqref {densities.blow.ups}. We can extract a diagonal subsequence from
$$\Sigma_{i,j,k}=\mu_{{\#}}^{k}(\mu_{{\#}}^{j}(\Sigma^{(i)}-p)-x_0),\quad i,j,k\in\N$$ 
denoted simply by $\{\Sigma_i\}_{i\in \N}$, where the relevant properties are
\begin{itemize}
\item[a)] $\partial \Sigma_i=0$;
\item[b)] $\Sigma_i$  {has generalized mean curvature tending to zero uniformly};
\item[c)] from Lemma \ref{density.smaller.two} there is $C>0$ such that for every $R$ and $i$ sufficiently large, we have  $$||\Sigma_i||(B^4_s(x))\leq  Cs^2\quad\mbox{for all }x\in B^4_R(0), 0\leq s\leq R.$$
\end{itemize}
From Federer Compactness Theorem we know that $\Sigma_i$ converges to $T\in {\mathcal Z}_2(S^3)$ in the flat topology. We  claim  that we can assign orientations to the half-planes $\{P_1,P_2,P_3\}$ so that $U=T$. This gives a contradiction because, regardless the orientation we assign to each half-plane,  we have $\partial U\neq 0$.

Denote by $L_j$ the set of all points at distance $2^{-j}$ from the line $L$, which is the line of common intersection of the half-panes $P_k$. We have that $U\setminus L_j$ consists of multiplicity one planes and thus, from property b) and Allard Regularity Theorem, we obtain that $\Sigma_i\llcorner(\R^4\setminus L_j)$ converges strongly to $U\llcorner (\R^4\setminus L_j)$ for every $j\in \N$. This induces an orientation on $U$. 

 Consider any $2$-form $\omega$ with support contained in  $B^4_R(0)\subset \R^4$, for some $R$,  and comass $||\omega||\leq 1$. We now  argue that $U(\omega)=T(\omega)$ and this finishes the proof. There is an integer $N$, independent of $j$, such that we can cover $L_j\cap B^4_R(0)$ with  balls  $\{B_k\}_{k=1}^{N2^j}$ of radius $2^{-j}$. Hence, we obtain from property c) that for all $i$ sufficiently large
$$||\Sigma_i||(L_j\cap B^4_R(0))\leq \sum_{k=1}^{N2^{j}}||\Sigma_i||(B_k)\leq CN2^{-j}.$$
The strong convergence property of $\Sigma_i$ outside $L_j$ implies at once that
$$|T(\omega)-U(\omega)|\leq 2CN2^{-j},$$
and thus, making $j\to \infty$, we obtain $U(\omega)=T(\omega)$.

\end{proof}

We are left to argue that the genus of $\Sigma$ must be bigger than zero. Indeed, because $\Sigma$ has multiplicity one, Allard Regularity Theorem implies the sequence $\Sigma_i$ converges strongly to $\Sigma$ and thus its  genus $g(\Sigma)\geq 1$.
\end{proof}

%%%%%%%%%%%%%%%%%%%%%%%%%%%%%%%%%%%%%%%%%%%%%%%%%%%%%%%%
%%%%%%%%%%%%%%%%%%%%%%%%%%%%%%%%%%%%%%%%%%%%%%%%%%%%%%%%%%%
\section{}
\label{conformal.images}

In this appendix we collect some facts about conformal transformations of $\mathbb{R}^4$.  
For each $v \in B^4$, let $F_v :S^3 \rightarrow S^3$ be given by
$$
F_v(x) = \frac{(1-|v|^2)}{|x-v|^2}(x-v) -v.
$$
Given $p,N\in  S^3$ with $\langle p, N\rangle = 0$, we define
\begin{align*}
\Delta(p,N,r) &= S^3 \setminus \Big(B_{r}\big((\cos  r) p + (\sin r) N\big)
 \cup \, B_{r}\big((\cos  r) p - (\sin r) N\big) \Big)\\
&= S^3 \setminus \left(B^4_{\sqrt{2(1-\cos r)}}\big((\cos  r) p + (\sin r) N\big)\right.\\
 &\quad\quad \quad\quad\cup \,\left. B^4_{\sqrt{2(1-\cos r)}}\big((\cos  r) p - (\sin r) N\big) \right).
\end{align*}

\subsection{Proposition}\label{sphere.image} {\em There is  $C_0>0$ and, for each $r\in(0,\pi/4),$  $C_1=C_1(r)>0$ and $\varepsilon_0=\varepsilon_0(r)>0$ such that the following holds.

For every
$$v = (1-s)(\cos t \, p+ \sin t \, N),$$
 with $$p,N \in S^3, \quad\langle p, N\rangle = 0,\quad 0<s\leq \varepsilon_0,\quad\mbox{and}\quad |t| \leq \varepsilon_0,$$ we have
 $$
B^4_{\overline{R}-C_0\sqrt{|(s,t)|}}(\overline{Q}) \cap S^3\subset F_v\left(B^4_{\sqrt{2}}(-N) \cap S^3\right) \subset B^4_{\overline{R}+C_0\sqrt{|(s,t)|}}(\overline{Q}) \cap S^3
$$
and
$$
F_v(\Delta(p,N,r)) \subset B^4_{\overline{R}+C_1\sqrt{|(s,t)|}}(\overline{Q}) \setminus B^4_{\overline{R}-C_1\sqrt{|(s,t)|}}(\overline{Q}),
$$
where
\begin{eqnarray*}
\overline{Q}&=&-\frac{t/s}{\sqrt{1+(t/s)^2}}p-\frac{1}{\sqrt{1+(t/s)^2}}N,\\
\overline{R}&=& \sqrt{2\left(1-\frac{t/s}{\sqrt{1+(t/s)^2}}\right)}.
\end{eqnarray*}
}
\begin{proof}
The next of lemma collects some basic identities, the proof of which are left to the reader.
\subsection{Lemma}\label{intersection.sphere}\textit{
\begin{itemize}
\item[(i)] Let $\tilde{Q} \in \mathbb{R}^4\setminus \{0\}$ and $\tilde{R}\geq 0$  such that  $(1-|\tilde{Q}|)^2 \leq \tilde{R}^2$. Then
$$
B^4_{\tilde{R}}(\tilde{Q}) \cap S^3 = B^4_R\left(\frac{\tilde{Q}}{|\tilde{Q}|}\right) \cap S^3\quad\mbox{where}\quad R=\sqrt{2+\frac{\tilde{R}^2-|\tilde{Q}|^2-1}{|\tilde{Q}|}}.
$$
\item[(ii)]Let $Q\in S^3$. Then
$$
S^3 \setminus B^4_{\sqrt{2(1-\cos \alpha)}}(Q)=\overline{B}^4_{\sqrt{2(1+\cos \alpha)}}(-Q) \cap S^3.
$$
\item[(iii)]Let $h \in \mathbb{R}^4 \setminus \{0\}$, $|h| \leq 1$,
and  $$E =  \{ x \in \mathbb{R}^4: \langle x - h, h \rangle \geq 0\}.$$ Then
$
E \cap S^3 =B^4_{\sqrt{2(1-|h|)}}(h/|h|)\cap S^3.
$
\end{itemize}
}

Let $i:\mathbb{R}^4 \setminus \{0\}\rightarrow \mathbb{R}^4$, $T_w:\mathbb{R}^4  \rightarrow \mathbb{R}^4$,
 and $D_{\lambda}:\mathbb{R}^4  \rightarrow \mathbb{R}^4$ be the conformal transformations given by
$$i(x)=\frac{x}{|x|^2},\quad  T_w(x)=x+w,\quad D_{\lambda}(x)=\lambda x,$$
where $\lambda \in \mathbb{R}$ and $w\in \R^4$. We have
\begin{equation}\label{fv.appendix}
F_v = D_{1-|v|^2}\circ T_{-\frac{v}{1-|v|^2}}\circ i\circ T_{-v}.
\end{equation}

\subsection{Lemma}\label{inversion.hyperplane}\textit{
Let $h \in \mathbb{R}^4 \setminus \{0\}$ and  $E =  \{ x \in \mathbb{R}^4: \langle x - h, h \rangle \geq 0\}.$ Then
$$
i(E) = \overline{B}^4_r(c),\quad\mbox{where}\quad c=\frac{h}{2|h|^2}, \quad r=\frac{1}{2|h|}.
$$
}
\begin{proof}
The lemma follows from the calculation:
\begin{eqnarray*}
|i(x)-c|^2 - r^2 &=&\left|\frac{x}{|x|^2} - \frac{h}{2|h|^2}\right|^2 - \frac{1}{4|h|^2}\\
&=& -\frac{\langle x- h,h\rangle}{|h|^2|x|^2}.
\end{eqnarray*}
\end{proof}

\subsection{Lemma}\label{hyperplane.image}\textit{
Let $h \in \mathbb{R}^4 \setminus \{0\}$, $|h| \leq 1$. If $v \in B^4$, then 
$$F_v\left(B^4_{\sqrt{2(1-|h|)}}\left(\frac{h}{|h|}\right)\cap S^3\right)=B^4_R\left(\frac{Q}{|Q|}\right) \cap S^3,$$
where 
\begin{eqnarray*}
Q &=& (1-|v|^2)h-2(|h|^2-\langle h,v\rangle)v,\\
R &=& \sqrt{2\left(1-\frac{|h|^2(1+|v|^2)-2\langle h,v\rangle}{|Q|}\right)}.
\end{eqnarray*}
}
\begin{proof}
From Lemma \ref{intersection.sphere}(iii) and \eqref{fv.appendix} we have
\begin{eqnarray*}
F_v\left(B^4_{\sqrt{2(1-|h|)}}\left(\frac{h}{|h|}\right)\cap S^3\right) &=&D_{1-|v|^2}\circ T_{-\frac{v}{1-|v|^2}}\circ i\circ T_{-v}(E \cap S^3)\\
&=& \left(D_{1-|v|^2}\circ T_{-\frac{v}{1-|v|^2}}\circ i\circ T_{-v}(E)\right) \cap S^3,
\end{eqnarray*}
where $E =  \{ x \in \mathbb{R}^4: \langle x - h, h \rangle \geq 0\}.$

Suppose $|h|^2-\langle h,v\rangle \neq 0$ and set
$$\sigma=1\quad\mbox{if  }v \in \,E,\quad \sigma=-1\quad\mbox{if }v\notin E,\quad\mbox{and}\quad h_v=\frac{|h|^2-\langle h,v\rangle}{|h|^2}h.$$
Then
\begin{equation}\label{plane.appendix.conformal}
T_{-v}(E)=\{ x \in \mathbb{R}^4: \sigma\langle x - h_v, h_v \rangle \geq 0\}.
\end{equation}

Suppose $|h|^2-\langle h,v\rangle >0$, i.e.,  $v$ is in the interior of $E$.
We have from Lemma \ref{inversion.hyperplane}
$$i(T_{-v}(E))=\overline{B}^4_r(c)$$
where
\begin{equation}\label{candr.appendix}
c=\frac{h}{2(|h|^2-\langle h,v\rangle)},\quad\mbox{and}\quad r=\frac{|h|}{2||h|^2-\langle h,v\rangle|}.
\end{equation}
Therefore
$$
D_{1-|v|^2}\circ T_{-\frac{v}{1-|v|^2}}\circ i\circ T_{-v}(E) =  \overline{B}^4_{(1-|v|^2)r}\big((1-|v|^2)c-v\big)
$$
and we conclude that $F_v(E)=B^4_{\tilde{R}}(\tilde{Q})$, where 
\begin{eqnarray*}
\tilde{Q}&=&\frac{(1-|v|^2)h}{2(|h|^2-\langle h,v\rangle)}-v,\\
\tilde{R}&=&\frac{(1-|v|^2)|h|}{2||h|^2-\langle h,v\rangle|}.
\end{eqnarray*}
It follows from Lemma \ref{intersection.sphere}(i) that 
$$F_v(E)=B^4_{\hat{R}}\left(\frac{\tilde{Q}}{|\tilde{Q}|}\right),\quad\mbox{where}\quad
\hat{R}=\sqrt{2+\frac{\tilde{R}^2-|\tilde{Q}|^2-1}{|\tilde{Q}|}}.
$$

Since
\begin{eqnarray*}
\tilde{R}^2-|\tilde{Q}|^2-1&=&\frac{(1-|v|^2)\langle h,v\rangle}{|h|^2-\langle h,v\rangle} - (1+|v|^2)\\
&=&\frac{2\langle h,v\rangle-|h|^2(1+|v|^2)}{|h|^2-\langle h,v\rangle},
\end{eqnarray*}
we have
$$
\hat{R}=\sqrt{2\left(1- \sigma \, \frac{|h|^2(1+|v|^2)-2\langle h,v\rangle}{|Q|}\right)}.
$$
Lemma \ref{hyperplane.image} follows immediately when $|h|^2-\langle h,v\rangle >0$ because $\frac{Q}{|Q|}=\frac{\tilde{Q}}{|\tilde{Q}|}$ and $R=\hat{R}$. 

Suppose now $|h|^2-\langle h,v\rangle <0$, i.e., $v\notin E$.  From \eqref{plane.appendix.conformal} and Lemma \ref{inversion.hyperplane} we have
$$i(T_{-v}(E))=\R^4\setminus{B}^4_r(c)$$
where $c$ and $r$ are as in \eqref{candr.appendix}. Therefore
$$
D_{1-|v|^2}\circ T_{-\frac{v}{1-|v|^2}}\circ i\circ T_{-v}(E) =  \mathbb{R}^4 \setminus B^4_{(1-|v|^2)r}\big((1-|v|^2)c-v\big).
$$
Thus $$F_v(E)=\mathbb{R}^4 \setminus B^4_{\tilde{R}}(\tilde{Q})=\mathbb{R}^4 \setminus B^4_{\hat{R}}(\frac{\tilde{Q}}{|\tilde{Q}|}),$$ where $\tilde{Q}$, $\tilde{R}$ and $\hat R$  are as above. Since $|h|^2-\langle h,v\rangle >0$, we have $\frac{Q}{|Q|}=-\frac{\tilde{Q}}{|\tilde{Q}|}$ and $R^2+\hat{R}^2=4$. We  apply Lemma \ref{intersection.sphere}(ii) and conclude Lemma \ref{hyperplane.image}.

Finally, if $|h|^2-\langle h,v\rangle =0$  then the result follows from the previous cases by approximation, since the set of all $v$ with 
$|h|^2-\langle h,v\rangle \neq 0$ is everywhere dense in $B^4$.
\end{proof}

Next we compute the conformal image of a geodesic ball in $S^3$.
\subsection{Lemma}\label{hyperplane.image.2}\textit{
Let $x \in S^3$. If  $v \in B^4$, then  $$F_v \left(B^4_{\sqrt{2}}(x) \cap S^3\right) = B^4_R\left(\frac{Q}{|Q|}\right) \cap S^3,$$ where
$$Q = (1-|v|^2)x+2 \langle x,v\rangle v\quad\mbox{and}\quad R = \sqrt{2\left(1+\frac{2\langle x,v\rangle}{|Q|}\right)}.$$
}

\begin{proof}
We apply Lemma \ref{hyperplane.image} with $h_t=tx$ in place of $h$, and let $t$ go to zero.
\end{proof}

We can now prove  the first statement  of Proposition \ref{sphere.image}.

\subsection{Lemma}\textit{There is $C_0$ so that for every
$$v = (1-s)(\cos t \, p+ \sin t \, N),$$
 with $$p,N \in S^3, \quad\langle p, N\rangle = 0,\quad 0<s\leq 1/2,\quad\mbox{and}\quad |t| \leq 1/2,$$ we have
 $$
B^4_{\overline{R}-C_0\sqrt{|(s,t)|}}(\overline{Q}) \cap S^3\subset F_v\left(B^4_{\sqrt{2}}(-N) \cap S^3\right) \subset B^4_{\overline{R}+C_0\sqrt{|(s,t)|}}(\overline{Q}) \cap S^3,
$$
where $\overline Q$ and $\overline R$ are defined in Proposition \ref{sphere.image}.}
\begin{proof}
From Lemma \ref{hyperplane.image.2},  $$F_v\left(B^4_{\sqrt{2}}(-N(p)) \cap S^3\right)=B^4_R\left(-\frac{Q}{|Q|}\right) \cap S^3,$$ where
$$
Q = (1-|v|^2)N+2 \langle N,v\rangle v\quad\mbox{and}\quad
R = \sqrt{2\left(1-\frac{2\langle N,v\rangle}{|Q|}\right)}.
$$
Thus
\begin{align}\label{q.expression}
Q  &= (2s-s^2) N +2(1-s)^2 \sin t(\cos t \, p+ \sin t \, N)\\ \notag
&= (2s-s^2+2(1-s)^2\sin^2t)N+ 2(1-s)^2 \sin t \cos t\, p, \notag
\end{align}
and
$$|Q|^2 = (2s-s^2)^2 + 4(1-s)^2 \sin^2t.$$
Hence we can find $C_1$ so that for all $|s|\leq 1/2$ and $|t|\leq 1/2$ we have
\begin{equation}\label{appendix.bigo.Q}
\frac{s^2+t^2}{C_1}\leq |Q|^2\leq 4(s^2+t^2)(1+C_1|(s,t)|).
\end{equation}
This implies the existence of $C_2$ so that for all $|s|\leq 1/2$ and $|t|\leq 1/2$ we have
\begin{equation}\label{Q.inequality.appendix}
\left|\frac{2}{|Q|}-\frac{1}{\sqrt{s^2+t^2}}\right|\leq C_2.
\end{equation}
From this inequality, \eqref{appendix.bigo.Q}, \eqref{q.expression}, and $s>0$, we obtain constants $C_3$ and $C_4$ so that
\begin{multline*}
\left|\overline{Q}-\left(-\frac{Q}{|Q|}\right)\right|^2 = \left(\frac{2s-s^2+2(1-s)^2\sin^2t}{|Q|}-\frac{s}{\sqrt{s^2+t^2}}\right)^2\\
+ \left(\frac{2(1-s)^2\sin t\cos t}{|Q|}-\frac{t}{\sqrt{s^2+t^2}}\right)^2\\
\leq 2\left(\frac{-s^2+2(1-s)^2\sin^2 t}{|Q|}\right)^2
+2\left(\frac{2(1-s)^2\sin t\cos t-2t}{|Q|}\right)^2\\+
C_3|(s,t)|^2
\leq C_4|(s,t)|^2.
\end{multline*}
From \eqref{Q.inequality.appendix},  \eqref{appendix.bigo.Q}, and $s>0$, we obtain constants $C_5$ and $C_6$ so that
\begin{eqnarray*}
|R^2 -\overline{R}^2| &=& \left|2\left(1-\frac{2\langle N,v\rangle}{|Q|}\right)-2\left(1-\frac{t/s}{\sqrt{1+(t/s)^2}}\right)\right|\\
&=&\left|-\frac{4(1-s)\sin t}{|Q|}+\frac{2t}{\sqrt{s^2+t^2}}\right|\\
&\leq&\left|\frac{4t-4(1-s)\sin t}{|Q|}\right|+C_5|(s,t)|\\
&\leq& C_6|(s,t)|.
\end{eqnarray*}
Hence $|R-\overline{R}|\leq \sqrt{C_6}\sqrt{|(s,t)|}$. 

If we choose $C=\sqrt{C_4}+\sqrt{C_6}$, the result follows by applying the  triangle inequality.
\end{proof}

The next lemma finishes the proof of Proposition \ref{sphere.image}

\subsection{Lemma}\textit{For every $r\in(0,\pi/4)$ there is $C_1=C_1(r)$ and $\varepsilon_0=\varepsilon_0(r)$ so that for every
$$v = (1-s)(\cos t \, p+ \sin t \, N),$$
 with $$p,N \in S^3, \quad\langle p, N\rangle = 0,\quad 0<s\leq \varepsilon_0,\quad\mbox{and}\quad |t| \leq \varepsilon_0,$$ we have
 $$
F_v(\Delta(p,N,r) \cap S^3) \subset B^4_{\overline{R}+C_1\sqrt{|(s,t)|}}(\overline{Q}) \setminus B^4_{\overline{R}-C_1\sqrt{|(s,t)|}}(\overline{Q}),
$$
where $\overline Q$ and $\overline R$ are defined in Proposition \ref{sphere.image}.}
\begin{proof}
Let  $\sigma_i=(-1)^{i+1}$, $i=1,2$. Define $$B_{i}=B^4_{\sqrt{2(1-\cos r)}}\big((\cos  r) p +\sigma_i (\sin r) N\big)\cap S^3$$
and $h_i= (\cos r) (\cos r\, p+ \sigma_i \sin r\, N)$. Then, by Lemma \ref{hyperplane.image},   $$F_v(B_{i})=B^4_{R_i}\left(\frac{Q_i}{|Q_i|}\right)\cap S^3$$ where
\begin{eqnarray*}
Q_i &=& (1-|v|^2)h_i-2(|h_i|^2-\langle h_i,v\rangle)v,\\
R_i &=& \sqrt{2\left(1-\frac{|h_i|^2(1+|v|^2)-2\langle h_i,v\rangle}{|Q_i|}\right)}.
\end{eqnarray*}

Notice that
\begin{eqnarray*}
|h_i|^2-\langle h_i,v\rangle &=& \cos^2 r - (1-s)(\cos^2 r \cos t + \sigma_i \cos r \sin r \sin t)\\
&=& - \sigma_it\cos r \sin r  +s \cos^2 r  + \cos r \, O(|(s,t)|^2)
\end{eqnarray*}
and so
\begin{multline}\label{Qi3.appendix.conformal}
Q_i = (2s-s^2)(\cos r) (\cos r\, p+ \sigma_i \sin r\, N)\\
-2(1-s)(\cos t\,p+\sin t \, N)(|h_i|^2-\langle h_i,v\rangle)\\
%&=& 2\sigma_i \cos r \, \sin r \, t \, p + 2 \sigma_i \cos r \, \sin r \, s \, N+ \cos r \, O(|(s,t)|^2),
= 2\sigma_i t\cos r \, \sin r  \, p + 2 \sigma_i s\cos r \, \sin r  \, N+  O(|(s,t)|^2),
\end{multline}
and
\begin{align}\label{Qi.appendix.conformal}
|Q_i|^2 &= 4\cos^2 r\, \sin^2 r \, (s^2+t^2)+O(|(s,t)|^3).\\ \notag
&= 4\cos^2 r\, \sin^2 r \, (s^2+t^2) \left(1 + \frac{O(|(s,t)|)}{\sin^2 r}\right). \notag
\end{align}
Choose $\varepsilon_0=\varepsilon_0(r)$ so that for all $0<s\leq \varepsilon_0$, $|t|\leq \varepsilon_0$ we have
\begin{equation}\label{Qi2.appendix.conformal}
|Q_i|^2\geq 2\cos^2 r\, \sin^2 r \, (s^2+t^2).
\end{equation}
This inequality and \eqref{Qi.appendix.conformal} implies that for some $C_2=C_2(r)$ we have
\begin{equation}\label{Qi4.appendix.conformal}
\left| \frac{2\cos r\sin r}{|Q_i|}-\frac{1}{\sqrt{s^2+t^2}}\right| \leq C_2
\end{equation}
and therefore, from \eqref{Qi3.appendix.conformal} and \eqref{Qi2.appendix.conformal}, we have
\begin{multline*}
\left|\frac{Q_i}{|Q_i|}-(-\sigma_i \overline{Q})\right|\leq \left|\frac{Q_i}{|Q_i|}-\frac{2\sigma_i\cos r\sin r}{|Q_i|}(tp+sN)\right|+C_3|(s,t)|\\
 \leq \frac{O(|(s,t)^2|)}{|Q_i|}  +C_3|(s,t)|\leq C_4|(s,t)|,
\end{multline*}
for some constants $C_3=C_3(r), C_4=C_4(r)$.

Now
\begin{eqnarray*}
|h_i|^2(1+|v|^2)-2\langle h_i,v\rangle &=& 2(|h_i|^2-\langle h_i,v\rangle)+2\cos^2 r(-2s+s^2) \\
&=&-2\sigma_it \cos r \, \sin  r +  O(|(s,t)|^2)
\end{eqnarray*}
and thus, combining with the expression for $R_i$, \eqref{Qi2.appendix.conformal}, and \eqref{Qi4.appendix.conformal}, we obtain
\begin{multline*}
\left|R_i^2 -2-2\frac{\sigma_it}{\sqrt{t^2+s^2}}\right|\leq \left|R_i^2 -2-2\frac{\sigma_it\cos r\sin r}{|Q_i|}\right|+C_2|(s,t)|\\
\leq \frac{O(|(s,t)|^2)}{|Q_i|}+C_2|(s,t)|\leq C_5|(s,t)|
\end{multline*}
for some $C_5=C_5(r)$. We can then find $C_6=C_6(r)$ such that
$$
|R_1^2-(4-\overline{R}^2)| \leq C_6 |(s,t)|,\quad
|R_2^2-\overline{R}^2| \leq C_6 |(s,t)|,
$$
which means
$$|R_1-\sqrt{(4-\overline{R}^2)}|  \leq \sqrt{C_6} \sqrt{|(s,t)|}\quad\mbox{and}\quad|R_2-\overline{R}| \leq \sqrt{C_6} \sqrt{|(s,t)|}.$$

If we choose $C_1=\sqrt{C_4} + \sqrt{C_6}$, then
\begin{eqnarray*}
B^4_{\sqrt{4-\overline{R}^2}-C_3\sqrt{|(s,t)|}}(-\overline{Q}) \cap S^3&\subset&  F_v(B_{1}),\\
B^4_{\overline{R}-C_3\sqrt{|(s,t)|}}(\overline{Q}) \cap S^3  &\subset& F_v(B_{2}).
\end{eqnarray*}

We conclude that
\begin{multline*}
F_v(\Delta(p,N,r) \cap S^3) = F_v(S^3 \setminus (B_1 \cup B_2))\\
\subset  S^3 \setminus \left(B^4_{\sqrt{4-\overline{R}^2}-C_3\sqrt{|(s,t)|}}(-\overline{Q})
\cup \, B^4_{\overline{R}-C_3\sqrt{|(s,t)|}}(\overline{Q}) \right).
\end{multline*}

The result follows from Lemma \ref{intersection.sphere}(ii).

\end{proof}

\end{proof}

%%%%%%%%%%%%%%%%%%%%%%%%%%%%%%%%%%%%%%%%%%%%%%%%%%%
%%%%%%%%%%%%%%%%%%%%%%%%%%%%%%%%%%%%%

\section{}\label{appendix.map}

Given $m,j\in \N,$ we  construct 
$${\bf r}_m(j):I(m,j+q)_0\rightarrow S(m+1,j)_0\cup T(m+1,j)_0$$
satisfying:
\begin{itemize}
\item $q$ depends on $m$ but not on $j$;
\item if $x,y\in I(m,j+q)_0$  satisfy ${\bf d}(x,y)=1$, then
\begin{equation}\label{prop2restriction}
{\bf d}({\bf r}_m(j)(x),{\bf r}_m(j)(y))\leq m.
\end{equation}
\item if $x\in I_0(m,j+q)_0$, then
\begin{equation}\label{prop1restriction}
{\bf r}_m(j)(x)=({\bf n}(j+q,j)(x),[0]).
\end{equation}
\end{itemize}

Let
$$R_m:I^m\rightarrow (I_0^m\times[0,1]) \cup (I^m\times\{1\})\subset I_0^{m+1}$$
be a Lipschitz homeomorphism such that
 \begin{equation}\label{Rmapboundary}
 R_m(x)=(x,0)\quad\mbox{ for all } x\in I_0^m.
\end{equation}
We choose $q\in \N$ such that
\begin{equation}\label{lipschitz}
|R_m(x)-R_m(y)|\leq 3^{q-2}|x-y| \quad\mbox{ for all } x,y \in I^m.
\end{equation}

Let $K= S(m+1,j)_0\cup T(m+1,j)_0$. Given $x\in I(m,j+q)_0$, we choose  ${\bf r}_m(j)(x)\in K$ such that
$d({\bf r}_m(j)(x),R_m(x))=d(R_m(x),K).$ This choice might not be unique, but if $x\in I_0(m,j+q)_0$ we obtain from \eqref{Rmapboundary} that $${\bf r}_m(j)(x)=({\bf n}(j+q,j)(x),[0]).$$
This shows \eqref{prop1restriction}.
 If $x,y\in I(m,j+q)_0$ satisfy ${\bf d}(x,y)=1$, we get from \eqref{lipschitz} that
$$ |R_m(x)-R_m(y)|\leq 3^{q-2}3^{-(j+q)}=3^{-(j+2)}.$$
This implies that ${\bf r}_m(j)(x)$ and ${\bf r}_m(j)(y)$ must be contained in a common $m$-cell of $I_0(m+1,j)$. Hence  property \eqref{prop2restriction} follows as well.

%\subsection{Remark} In Lemma \ref{hyperplane.image}, 
%\begin{eqnarray*}
%|Q|^2 &=& (1-|v|^2)^2|h|^2 + 4(|h|^2|v|^2 - \langle h,v\rangle)(|h|^2-\langle h,v\rangle)\\
%&=&(1-|v|^2)^2|h|^2 -4(1-|v|^2)|h|^2(|h|^2-\langle h,v\rangle)+ 4(|h|^2 - \langle h,v\rangle)^2\\
%&=&(1-|v|^2)^2|h|^2(1-|h|^2)+\left(|h|^2(1+|v|^2)-2\langle h,v\rangle\right)^2.
%\end{eqnarray*}
%In particular, $|Q|\geq |h|^2(1+|v|^2)-2\langle h,v\rangle$.

%\medskip

%\subsection{Lemma}\label{parametrization}\textit{Let $p,N \in S^3$ with $\langle p,N\rangle=0$,
%$t\in [-\pi,\pi]$, and $k\in [-\infty\,+\infty]$. Let $f:\{p,N(p)\}^\perp \rightarrow S^3$ be given by
%\begin{eqnarray*}
%f(w)  &=&\frac{2\cos\, t}{(1+|w|^2)}\left(\frac{w}{\sqrt{1+k^2}}+\frac{p}{1+k^2}-\frac{kN(p)}{1+k^2}\right) \\
%&&+\frac{2\sin\,t}{(1+|w|^2)}\left(\frac{kw}{\sqrt{1+k^2}}+\frac{kp}{1+k^2}-\frac{k^2N(p)}{1+k^2}\right)\\
%&&-\cos \, t\,p + \sin\, t\,N(p),
%\end{eqnarray*}
%for $w\in \{p,N(p)\}^\perp$. Then $f$ parametrizes ***, where 
%}

%\begin{proof}
%\end{proof}

\bibliographystyle{amsbook}

\end{document}